\documentclass{amsart}
\usepackage{graphicx}

\newtheorem{thm}{Theorem}[section]
\newtheorem{cor}[thm]{Corollary}
\newtheorem{lem}[thm]{Lemma}
\newtheorem{prop}[thm]{Proposition}
\newtheorem{defi}[thm]{Definition}
\newtheorem{rem}[thm]{Remark}
\newtheorem*{thmnonnum}{Theorem}

\newenvironment{dimo}{\noindent{\bf Proof: }}{\hfill$\Box$\medskip}
\numberwithin{equation}{section}
\newenvironment{proofof}[2]{\begin{proof}[Proof of #1 \ref{#2}.]}{\end{proof}}
\newcommand{\beq}{\begin{equation}}
\newcommand{\eeq}{\end{equation}}
\newcommand{\beqs}{\begin{equation*}}
\newcommand{\eeqs}{\end{equation*}}

\newcommand\NN {{\mathbb N}}
\newcommand\ZZ {{\mathbb Z}}
\newcommand\QQ {{\mathbb Q}}
\newcommand\RR {{\mathbb R}}
\newcommand\CC {{\mathbb C}}

\newcommand\HH {{\mathbb H}}

\newcommand\TT {{\mathbb T}}

\newcommand\glduer{{\rm GL(2,\RR)}}
\newcommand\slduez{{\rm SL(2,\ZZ)}}
\newcommand\sldz{{\rm SL(d,\ZZ)}}
\newcommand\slduer{{\rm SL(2,\RR)}}

\newcommand\veech{{\rm SL}}

\newcommand\lungh{{\rm length }}

\newcommand\distanza{{\rm dist }}
\newcommand\angolo{{\rm angle }}
\newcommand\sys{{\rm Sys }}
\newcommand\area{{\rm Area }}

\newcommand\cardinalita{{\rm Card }}

\newcommand\hol {{\rm Hol}}

\newcommand\scatola{{\rm Box }}
\newcommand\periodici{{\rm Per }}

\newcommand\al{\alpha}
\newcommand\be{\beta}
\newcommand\ga{\gamma}

\newcommand\gatilde{\tilde{\gamma}}
\newcommand\de{\delta}

\newcommand\la{\lambda}

\newcommand\te{\theta}

\newcommand\Ga{\Gamma}
\newcommand\De{\Delta}

\newcommand\Th{\Theta}

\newcommand\cA{{\mathcal{A}  }}

\newcommand\cC{{\mathcal{C}  }}
\newcommand\cD{{\mathcal{D}  }}
\newcommand\cE{{\mathcal{E}  }}
\newcommand\cF{{\mathcal{F}  }}
\newcommand\cG{{\mathcal{G}  }}
\newcommand\cH{{\mathcal{H}  }}
\newcommand\cI{{\mathcal{I}  }}

\newcommand\cK{{\mathcal{K}  }}
\newcommand\cL{{\mathcal{L}  }}
\newcommand\cM{{\mathcal{M}  }}
\newcommand\cN{{\mathcal{N}  }}

\newcommand\cP{{\mathcal{P}  }}

\newcommand\cR{{\mathcal{R}  }}

\newcommand\cU{{\mathcal{U}  }}
\newcommand\cV{{\mathcal{V}  }}

\begin{document}

\title[Lagrange Spectra in Teichm\"uller Dynamics via renormalization]{Lagrange Spectra in Teichm\"uller Dynamics via renormalization}

\author{Pascal Hubert, Luca Marchese, Corinna Ulcigrai}

\address{LATP, Centre de Math\'ematiques et Informatique (CMI),
Universit\'e Aix-Marseille, 39 rue Joliot Curie, 13453 Marseille Cedex 13 France.}

\email{hubert@cmi.univ-mrs.fr}

\address{LAGA, Universit\'e Paris 13, Avenue Jean-Baptiste Cl\'ement, 93430 Villetaneuse, France.}

\email{marchese@math.univ-paris13.fr}

\address{School of Mathematics, University of Bristol, University Walk, Clifton, Bristol, BS8 1TW, United Kingdom.}

\email{corinna.ulcigrai@bristol.ac.uk}



%\date{}%

\begin{abstract}
We introduce Lagrange Spectra of closed-invariant loci for the action of $\slduer$ on the moduli space of translation surfaces, generalizing the classical Lagrange Spectrum, and we analyze them with renormalization techniques. A formula for the values in such spectra is established in terms of the Rauzy-Veech induction and it is used to show that any invariant locus has closed Lagrange spectrum and values corresponding to pseudo-Anosov elements are dense. Moreover we show that Lagrange spectra of arithmetic Teichm\"uller discs contain an Hall's ray, giving an explicit bound for it via a second formula which uses the classical continued fraction algorithm. In addition, we show the equivalence of several definitions of bounded Teichm\"uller geodesics and bounded type interval exchange transformations and we prove quantitative estimates on excursions to the boundary of moduli space in terms of norms of positive matrices in the Rauzy-Veech induction.
\end{abstract}

% ----------------------------------------------------------------

\maketitle

\tableofcontents

\section{Introduction}

The classical Lagrange Spectrum is the set $\cL$ of values $$
L(\al):=
\limsup_{q,p\to\infty}
\frac{1}{q|q\al-p|}
\textrm{ for }
\al\in\RR.
$$
The quantity $L(\al)$ has the following interpretation in  Diophantine Approximation. Recall that by Dirichlet theorem, given any irrational $\alpha \in \mathbb{R}$, there exists infinitely many integers $p,q$, with $q \neq 0$, such that $\left|\alpha - \frac{p}{q}\right|< \frac{1}{q^2}$; we have that
$$
L(\alpha) = \{ \sup c >0 \quad \textrm{such \ that\ } \left|\alpha - \frac{p}{q}\right|< \frac{1}{ c q^2}  \quad \textrm{for \ infinitely \ many } p,q \in \mathbb{Z}, q\neq 0 \}.
$$
The set $\cL$ has been studied in depth for many decades, see for example the survey by Cusick and Flahive \cite{cusick} and the references therein. In particular, some of the basic properties are summarized below.
\begin{enumerate}
\item
$\cL$ is a closed subset of the real line (Cusick 1975).
\item
The values $L(\be)$ for $\be$ quadratic irrational are dense in $\cL$ (Cusick 1975).
\item
$\cL$ contains an \emph{Hall's ray}, that is a positive half-line (Hall 1947).
\end{enumerate}
It is also known that $\cL$ begins with a discrete sequence (Markoff 1879), which coincides with the \emph{Markoff Spectrum} and whose smaller term is $\sqrt{5}$, which is called \emph{Hurwitz constant} (Hurwitz 1891). More recently, the structure of $\cL$ in between the discrete part and the Hall ray has been investigated by G. Moreira in \cite{gugu}. The \emph{continued fraction algorithm} plays a crucial role in the proof of all these results on $\cL$, thanks to the following beautiful formula for $L(\al)$ (Perron 1921). Let $\al=a_0+[a_1,a_2,\dots]$ be the continued fraction expansion of $\al$. Then we have
\beq\label{eqformulaclassica}
L(\al)=
\limsup_{n\to\infty}
\big(
[a_{n-1},\dots,a_{1}]+a_{n}+[a_{n+1},a_{n+2},\dots]
\big).
\eeq
It is clear from this formula that $L(\al)$ is finite if and only if $\al$ is of bounded type, that is the sequence $(a_n)_{n\in\NN}$ is bounded. It is well known that the continued fraction expansion is related to the coding of the geodesic flow on the modular surface (see \cite{series1, series0}). In this geometric context there is a way to establish a correspondence between bounded type real numbers $\al$ and bounded geodesic rays $(g^{\al}_t)_{t>0}$. Any value $L(\al)$ corresponds to the following geometric quantity
$$
\limsup_{t\to\infty}
\textrm{heigth}(g^{\al}_t),
$$
which gives the asymptotic depth of penetration of the ray $g^{\al}_t$ into the cusp of the modular surface, where $\textrm{heigth}(\cdot)$ denotes the hyperbolic height.

\medskip

Many generalizations of Lagrange and Markhoff Spectra have been studied by several authors. We do not attempt here to summarize the huge literature, but we just mention a partial list of examples. A first natural generalization of the modular surface are quotients of the upper half plane by Hecke groups and more in general by Fuchsian groups with cusps. Investigations of the Hurwitz constant and the
Markoff spectrum for Hecke and triangle groups were carried out in \cite{series2, haasseries, Vulakh:Triangle} and for Fuchsian
groups with cusps in \cite{Vulakh:Fuchsian}. Moreover, the existence of a Hall ray for penetration Markoff Spectra for any Riemann surface with cusps is shown in  \cite{schmidtsheingorn}. Discrete groups acting on higher dimensional hyperbolic spaces were
also studied, in particular the Hurwitz constant and the Lagrange Spectrum are studied in \cite{vulakh3, Vulakh:Farey}, and the closure of the spectrum is  special case of \cite{maucourant} and \cite{paulin1}. Lagrange and
Markhoff spectra for quotients by Bianchi groups have also  been studied since they yield number
theoretical applications for the approximation of a complex number by
elements of an imaginary quadratic number field, see
\cite{vulakh5} \cite{maucourant}. In \cite{maucourant}, Maucourant proves more in general the closure of Lagrange Spectra for manifolds with negative sectional curvature. Other geometric generalizations of the Lagrange Spectra in negative curvature, as spiraling spectra, where defined and investigated in various works by Paulin and Parkonnen in \cite{paulin2, paulin3} and by Paulin and Hersonsky in \cite{hersonskypulin}. In particular, they prove the existence of the Hall ray for all these spectra. Finally, very recently, S. Ferenczi studied in depth the Lagrange Spectrum of interval exchange transformations with three intervals \cite{ferenczi}.

\medskip

In this paper, we generalize Lagrange Spectra to the setting of interval exchange transformations and the \emph{Teichm\"uller flow} on moduli spaces of translation surfaces, which generalize respectively rotations on the circle and the geodesic flow on the modular surface. Using renormalization techniques, we prove generalizations of the classical results (1), (2) and (3) listed above for the geodesic flow on closed and invariant loci for the action of $\slduer$ in the moduli space. Our main tools are two explicit formulas to compute the values in the generalized Lagrange Spectra via two generalizations of the continued fraction algorithm (see Theorem  \ref{renormalizedformula} and Theorem \ref{formulasquaretiled}). The first is the so-called \emph{Rauzy-Veech induction} for interval exchange transformations, which provides a coding of the Teichm\"uller flow on full connected components of strata of moduli spaces.  The second is an extension of the classical continued fraction to \emph{square-tiled surfaces} and allows to code the restriction of the Teichm\"uller flow to their (closed) orbit under $\slduer$.

\subsection{Translation surfaces}

A \emph{translation surface} $X$ is a compact surface of genus $g$, with a flat metric and conical singularities whose angles are multiples of $2\pi$. Alternatively, a translation surface $X$ is a datum $(S,w)$, where $S$ is a compact Riemann surface of genus $g$ and $w$ is an holomorphic 1-form on $S$ with zeros of orders $k_1,\dots,k_r$ at points $p_1,\dots,p_r$. The two definitions are equivalent, one direction of the equivalence being evident, since $w\otimes\overline{w}$ is a flat metric on $S$ with cone angle $2(k_i+1)\cdot\pi$ at each $p_i$. Conversely, the data $(S,w)$ can be recovered from $X$ as follows. Outside the conical singularities $p_1,\dots,p_r$, an atlas $S$ of flat charts can be choosen identifying flat neighborhoods of $X$ with open sets in $\CC$, that is introducing a local coordinate $z$ on each flat neighborhood. Since the angles at $p_1,\dots,p_r$ are multiples of $2\pi$, then the atlas $S$ can be chosen so that $z'=z+const$ whenever any two such charts $z$ and $z'$ overlap, thus it is well-defined also an holomorphic 1-form $w$ whose expression in these charts is $dz$. On a small neighborhood of a conical singularity $p_i$ the 1-form $w$ has a local primitive $\zeta_i:=\int w$ and the latter satisfies $d\zeta_i=z^{k_i}dz$ for any flat chart $z$ overlapping with $\zeta_i$, where $2(k_i+1)\cdot\pi$ is the conical angle at $p_i$. Therefore the atlas $S$ extends holomorphically to each $p_i$ and the same does $w$, admitting a zero of order $k_i$.

For any translation surface $X$ the associated 1-form $w$ induces an area form $(\sqrt{-1}/2)w\wedge\overline{w}$, different from zero outside $p_1,\dots,p_r$, so that the euclidian area of $X$ is given by
$
\area(X)=
(\sqrt{-1}/2)
\int_X w\wedge\overline{w}
$.
The 1-form $w$ induces also a pair of parallel vector fields $\partial_x$ and $\partial_y$ by $w(\partial_x)=1$ and $w(\partial_y)=\sqrt{-1}$. These vector fields are not complete, since they have $(k_i+1)$ entering trajectories and $(k_i+1)$ outgoing trajectories at the conical singularity $p_i$. The \emph{vertical flow} $\phi^t$ of $X$ is the integral flow of $\partial_y$.

\subsubsection{Strata and action of $\slduer$}\label{s1ss1sss1}

For a translation surface $X$ the genus and the orders of zeroes satisfy the relation $k_1+\dots+k_r=2g-2$. Fox fixed integers $k_1,\dots,k_r$ satisfying the last relation, denote $\cH(k_1,\dots,k_r)$ the corresponding \emph{stratum} of the \emph{moduli space} of translation surfaces, that is the set of translation surfaces whose associated 1-form $w$ has $r$ zeroes with orders $k_1,\dots,k_r$. It is a complex orbifold with complex dimension $2g+r-1$. Consider a translation surface $X=(S,w)$ in the stratum $\cH(k_1,\dots,k_r)$ and $A\in\slduer$. A new translation surface
$
A\cdot X=(A_\ast S,A_\ast w)
$
is defined, where the 1-form $A_\ast w$ is the composition of $w$ with $A$ and $A_\ast S$ is the complex atlas for which $A_\ast w$ is holomorphic. Therefore the group $\slduer$ acts on $\cH(k_1,\dots,k_r)$. The \emph{Teichm\"uller flow} $\cF_t$ is the action of the diagonal subgroup, that is
$$
\cF_{t}=
\left(
\begin{array}{cc}
e^{t} & 0 \\
0 & e^{-t}
\end{array}
\right).
$$
The Teichm\"uller flow preserves the hypersurface $\cH^{(1)}(k_1,\dots, k_r) \subset \cH(k_1,\dots, k_r)$ consisting of all $X \in \cH(k_1,\dots,k_r)$ with $\area(X)=1$.

Other relevant subgroups are
$
\cR_\te:=
\begin{pmatrix}
\cos\te & -\sin\te \\
\sin\te & \cos\te
\end{pmatrix}
$
and the \emph{horocyclic flow}
$
\cU_s:=
\begin{pmatrix}
1 & 0 \\
s & 1
\end{pmatrix}
$.

\subsubsection{Bounded geodesics in moduli space}

A \emph{saddle connection} is a geodesic segment $\ga$ for the flat metric of $X$
starting and ending in two conical singularities and not containing any other
conical singularity in its interior. The set $\hol(X)$ of \emph{periods} of $X$
is the set of complex numbers $v:=\int_\ga w$, where $\ga$ is a saddle
connection for $X$ and $w$ is the holomorphic one-form.  Compact subsets in the
stratum are characterized as follows by the \emph{systole} function
$X\mapsto\sys(X)$, defined as the length of the shortest saddle connection of the translation surface $X$, that is
$$
\sys(X):=\min\{|v|;v\in\hol(X)\}.
$$
According to the so-called \emph{Mahler criterion}, a sequence $X_n$ is bounded in $\cH(k_1,\dots,k_r)$ if and only if there exists some $\delta>0$ such that $\sys(X_n)\geq\delta$ for any $n$. We study bounded positive orbits for the action of the Techm\"uller flow on strata. A dynamical estimate of the asymptotic maximal excursion of the positive orbit $\cF_t(X)$ is given by
\beq\label{eq:definitions}
s(X):=\liminf_{t\to\infty}
\frac{1}{\area(X)}\sys(\cF_t\cdot X).
\eeq

For any period $v$ of a translation surface $X$ we set $\area(v)=|\Re(v)|\cdot|\Im(v)|$, where $\Re(v)$ and $\Im(v)$ denote respectively the real and imaginary part of $v \in \mathbb{C}$ (geometrically, $\area(v)$ is the Euclidean area of a rectangle which has $v$ as diagonal). In terms of the flat geometry of $X$, we obtain an estimate of the asymptotic maximal excursion of $\cF_t(X)$ setting
\beq\label{eq:definitiona}
a(X):=\liminf_{|\Im(v)|\to\infty}
\frac{\area(v)}{\area(X)}.
\eeq
For any translation surface $X$ the two quantities introduced above are related by $s(X)=\sqrt{2a(X)}$ (see Proposition \ref{prop:vorobetsidentity} below). In particular, the forward Teichm{\"u}ller geodesics $\cF_t(X)$ is boundeded in moduli space iff, equivalently,  $s(X)>0$ or $a(X)>0$.

\subsection{Interval Exchange Transformations}

An alphabet is a finite set $\cA$ with $d\geq2$ elements. An IET (which is a shortening for \emph{interval
exchange transformation}), is a map $T$ from an interval
$I$ to itself such that $I$ admits two partitions in subintervals
$$
\cP_t:=\{I^t_\al;\al\in\cA\}
\textrm{ and }
\cP_b:=\{I^b_\be;\be\in\cA\}
$$
and for any $\al$ in $\cA$ the restriction of $T$ to $I^t_\al$ is the translation onto $I^b_\al$. For any $\al$ in $\cA$ the intervals $I^t_\al$ and $I^b_\al$ have the same length, that we denote $\la_\al$. We call \emph{length datum} of $T$ the vector $\la$ in $\RR_+^\cA$ whose $\al$-entry is $\la_\al$ for any $\al$ in $\cA$. We call \emph{combinatorial datum} of $T$ the pair of bijections $\pi=(\pi^t,\pi^b)$ from $\cA$ to $\{1,\dots,d\}$ such that for any $\al$ in $\cA$, if we count starting from the left, $I^t_\al$ is in the $\pi^t(\al)$-th position in $\cP_t$ and $I^b_\al$ is in the $\pi^b(\al)$-th position in $\cP_b$. On the interval $I:=(0,\sum_{\chi\in\cA}\la_\chi)$ the data $(\pi,\la)$
determine uniquely $T$. We say that the combinatorial datum $\pi$ is \emph{admissible} if there is no proper subset $\cA'\subset\cA$ with $k<d$ elements such that $\pi^t(\cA')=\pi^b(\cA')=\{1,\dots,k\}$. The parameter space of all IETs with combinatorial datum $\pi$ is $\{\pi\}\times\RR^\cA_+$. IETs have been largely studied, we refer for example to reader to the lecture notes \cite{viana, yoccoz, zorichdue}.

IETs are strictly linked to translation surfaces and to the
Teichm\"uller flow on their moduli space. More precisely, any
translation surface $X$ has an unitary constant vector field whose \emph{first
return map} to a transverse segment $I$ in $X$ is an IET. In \cite{veech}, W. Veech gave a
combinatorial presentation of a translation surface in relation to IETs, known as \emph{zippered rectangles construction}. If $T$ is an IET defined by the data $(\pi,\la)$, where the combinatorial datum $\pi$ is admissible, then there exists a \emph{suspension
datum} $\tau$, that is a vector $\tau\in\RR^\cA$ such that the triple of data $(\pi,\la,\tau)$ defines a translation surface $X$, whose first return map is $T$. The details of Veech's zippered rectangles construction are summarized in \S \ref{backgroundss1}. It is possible to see that if $\hol(X)$ does not contain pure imaginary elements, then there exist data $\pi,\la,\tau$ such that $X=X(\pi,\la,\tau)$. For this last condition to hold it is sufficient to have
$$
\liminf_{t\to\infty}
\sys(\cF_t\cdot X)<+\infty,
$$
which is obviously satisfied by bounded geodesic rays.

\subsubsection{Bounded Type IETs}\label{s1ss2sss1}

Let $X$ be a translation surface and consider an horizontal segment $I$ embedded in $X$. Let $T$ be the IET induced as first return to $I$ of the vertical flow of $X$. Following Boshernitzan, for any positive integer $n$ we define $\cE_{n}(T)$
as the minimum of the distance between two singularities of the $n$-th
iterated $T^{n}$ of $T$. \emph{Boshernitzan's criterion} says that if $T$ is a minimal IET which is not uniquely ergodic, then
$
\lim_{n\to\infty}n\cdot\cE_n(T)=0
$
(see \cite{boshernitzan} and \cite{veechboshernitzancriterion}). The last criterion suggests a third way to estimate the asymptotic maximal excursion of the orbit of $X$, which is to define
\beq\label{eq:definitione}
\cE(T):=
\liminf_{n\to\infty}
\frac{n\cdot\cE_{n}(T)}{|I|}.
\eeq
It turns out that for any $X$ and $T$ related as explained above we always have $\cE(T)=a(X)$ (see  Proposition \ref{prop:vorobetsidentity} below).  In particular, the positive orbit $\cF_t(X)$ is bounded in moduli space if and only if $\cE(T)>0$. In analogy with the classical correspondence between bounded geodesics on the modular surface and numbers $\alpha$ of bounded type, we will say that  an IET $T$ is on \emph{bounded type} if $\cE(T)>0$. In Section \ref{sec:finitetime} we will state another characterization of bounded type IETs in terms of the Rauzy-Veech induction which generalize the definition of $\alpha$ of bounded type in terms of bounded entries in the continued fraction expansion (see Corollary \ref{cor:boundedtypeIETs}).

\subsection{Lagrange spectra of invariant loci}

The notion of bounded type translation surfaces and IETs corresponds to the positivity of the quantities $s(X)$, $a(X)$ and $\cE(T)$ introduced in respectively in Equations \eqref{eq:definitions}, \eqref{eq:definitiona} and \eqref{eq:definitione}. We have actually just one notion, since according to Proposition \ref{prop:vorobetsidentity} below these three asymptotic quantities are in fact the same. Furthermore bounded type IETs can be defined in terms of norms of positive matrices in the so-called \emph{Kontsevich-Zorich cocycle} (see \S~\ref{sec:finitetime}) and according to Theorem \ref{thm1s4} in \S~\ref{sec:finitetime} this still leads to the same notion. Theorem \ref{thm1s4} is a complement to Proposition \ref{prop:vorobetsidentity} and establishes a fourth relation involving the norm of positive matrices and $a(X)$, or equivalently $s(X)$ or $\cE(T)$. Proposition \ref{prop:vorobetsidentity} is due to Y. Vorobets (see \cite{vorobets}). We will give an alternative proof in \S \ref{section2}.

\begin{prop}[Vorobets' identity]\label{prop:vorobetsidentity}
Consider data $(\pi,\la,\tau)$, let $X$ be the underlying translation surface and $T$ be the IET corresponding to $(\pi,\la)$. We have
$$
\cE(T)=a(X)=\frac{s^{2}(X)}{2}.
$$
\end{prop}

The group $\slduer$ acts on strata of translation surfaces. Recently, A. Eskin and M. Mirzakhani \cite{eskinmirza} proved that any closed and invariant subset for such action is an affine sub-orbifold of the corresponding stratum, that is a nice moduli space itself (we refer to \cite{eskinmirza} for definitions). Therefore it is natural to study the asymptotic quantities above for any such locus.

\begin{defi}\label{def:Lagrange}
Let $\cI$ be an \emph{invariant locus}, that is a subset of some stratum, which is closed and invariant under
$\slduer$. We define the \emph{Lagrange Spectrum} of the invariant locus $\cI$
as
$$
\cL(\cI):=\{a^{-1}(X);X\in\cI\},
$$
where $a^{-1}(X)$ denotes the inverse $1/a(X)$.
\end{defi}

An interesting question is to compute the \emph{Hurwitz constant} of the spectra, that is the smallest element of any spectrum $\cL(\cI)$. We can establish the following immediate lower bound (the proof is given in \S \ref{appendixAuxilSs1}).

\begin{lem}\label{lem:lowerbound}
Let $\cI$ be an invariant locus contained in the stratum $\cH(k_1,\dots, k_r)$ of translation surfaces of genus $g$, where $2g-2=k_1+\dots+k_r$. We have
$$
\cL(\cI)
\subset
\left[
{\pi\cdot\frac{2g+r-2}{2} },+\infty
\right].
$$
\end{lem}

\subsubsection{Lagrange spectra of Veech surfaces}

Focusing on closed orbits of $\slduer$ it is possible to see that Definition \ref{def:Lagrange} generalizes the classical definition of Lagrange spectrum.

The \emph{Veech group} of a translation surface $X$ is the stabilizer of $X$ under the action of $\slduer$, and it is denoted $\veech(X)$. We say that a translation surface $X$ is a \emph{Veech surface} if $\veech(X)$ has finite co-volume in $\slduer$. It is well-known (see \cite{smillieweiss}) that $X$ is a Veech surface if and only if its orbit $\slduer\cdot X$ is closed in the stratum, that is it is a closed-invariant locus $\cI_X$ for the action of $\slduer$, and in particular it can be identified with the unitary tangent bundle of the upper half plane quotiented by $\veech(X)$.

\medskip

With the notation of \S \ref{s1ss1sss1}, recall that for any $A$ in $\slduer$ there are $\te$ in
$[0,2\pi)$ and $t$ and $s$ in $\RR$ such that
$A=\cR_\te\cdot\cF_t\cdot\cU_s$. The map $X\mapsto s(X)$ is obviously invariant under $\cF_t$. Moreover for any $t$ and $s$ we have
$
\cF_t\circ\cU_s=\cU_{se^{-2t}}\circ\cF^t
$,
hence $X\mapsto s(X)$ is also invariant under $\cU_s$, since the function $\sys(\cdot)$ is continuous. Recalling Proposition \ref{prop:vorobetsidentity} we can state the following Lemma.

\begin{lem}[Symmetries]\label{lem:symmetries}
Let $X$ be any translation surface. For any $s$
and $t$ in $\RR$ we have
$$
a(\cF_t\cdot X)=
a(\cU_s\cdot X)=
a(X).
$$
\end{lem}

Observe that for any $\te\in[0,2\pi)$ we have the extra symmetry
$
a(\cR_{\te+\pi}\cdot X)=
a(\cR_\te\cdot X)
$,
indeed for any $\te$ we
obviously have
$
\hol(\cR_{\te+\pi}\cdot X)=
-\hol(\cR_{\te}\cdot X)
$.
It follows that if $X$ is a \emph{Veech surface} and $\cI_X:=\slduer\cdot X$ is its closed orbit under $\slduer$, then Lemma \ref{lem:symmetries} implies
$$
\cL(\cI_X)=
\{a^{-1}(\cR_\te\cdot X);
-\frac{\pi}{2}<\te\leq \frac{\pi}{2}\}.
$$
For convenience we pass to the variable
$
\al=\tan(\te)\in(-\infty,+\infty]
$
and we fix the convention that the horizontal direction on the Veech surface $X$ corresponds to $\al=+\infty=\tan(\pi/2)$ and that it is a periodic direction. Moreover we will write simply $\cL(X)$ instead of $\cL(\cI_X)$. With these convention, the Lagrange Spectrum of a Veech surface $X$ is parametrized by a function $L_X:(-\infty,+\infty)\to\cL(X)$, that we call
\emph{standard parametrization}, defined by

\beq\label{eq:standardparametrization}
L_X(\al):=
\frac{1}{a(\cR_{\arctan(\al)}\cdot X)}.
\eeq

We can now explain how we can recover the classical Lagrange spectrum from these definitions. Let $\TT^2$ be the \emph{standard torus}, that is $\TT^2=\RR^{2}/\ZZ^{2}$. Its Veech group $\veech(\TT^2)$ is $\slduez$ and the stratum $\cH(0)$, that is the moduli space of flat tori, coincides with
$
\slduer\cdot\TT^2=\slduer/\slduez
$.
Therefore in genus one the whole parameter space is the only invariant locus and the only Lagrange
spectrum is $\cL(\TT^2)$. In term of Equation \eqref{eq:standardparametrization} it is easy to verify that (a more general formula is established by Lemma \ref{lem2s5ss1})
$$
L_{\TT^2}(\al)=
\limsup_{q,p\to\infty}
\frac{1}{q|q\al-p|}.
$$
Thus $\cL(\TT^2)$ coincides with $\cL$ and Definition \ref{def:Lagrange} is a generalization of the classical definition of Lagrange spectrum.

\subsection{Lagrange spectra via the Rauzy-Veech induction}

The \emph{Rauzy-Veech induction} is  a generalization for translation surfaces of the continued fraction expansion. It is an invertible map acting on data $(\pi,\la,\tau)$ coming from the zippered rectangles construction, whose iteration gives a sequence of data $(\pi^{(r)},\la^{(r)},\tau^{(r)})_{r\in\ZZ}$. The induction was defined  by W. Veech in \cite{veech}, as \emph{natural extension} of a two-to-one induction map acting on IETs, introduced by G. Rauzy in \cite{rauzy}. The induction at the level of IETs  has an acceleration with a finite and ergodic smooth invariant measure, which was discovered in \cite{zorichuno} by A. Zorich. Details on the Rauzy and Rauzy-Veech induction can be found in \S \ref{backgroundss2} and \S \ref{backgroundss3}.

Recall that a real number $\be$ is quadratic irrational if and only if its continued fraction expansion is eventually periodic. In higher genus this notion is naturally generalized by translation surfaces whose data $(\pi,\la,\tau)$ have periodic Rauzy-Veech induction. Therefore, for any invariant locus $\cI$ we are lead to consider the set $PA(\cI)$ of those translation surfaces $X$ in $\cI$ whose data $(\pi,\la,\tau)$ have periodic Rauzy-Veech induction. Equivalently, such $X$ correspond to periodic points for the restriction of the Teichm\"uller flow $\cF_t$ on $\cI$, or to translation surfaces admitting a \emph{pseudo-Anosov diffeomorphism}.

\subsubsection{Closure of Lagrange spectra and density of periodic values}

One of the main tools of this paper is the following formula for the values in Lagrange spectra $\cL(\cI)$ in terms of the Rauzy-Veech induction (the formula is established by Theorem \ref{renormalizedformula}). If $X$ is a translation surface corresponding to combinatorial, length and suspension data $(\pi,\la,\tau)$ and  $(\pi^{(r)},\la^{(r)},\tau^{(r)})_{r\in\ZZ}$ is its sequence of iterates under Rauzy-Veech induction, we have
$$
\frac{1}{a(X)}=
\limsup_{r\to\infty}
\frac{1}{w(\pi^{(r)},\la^{(r)},\tau^{(r)})},
$$
where $w$ is a continuous function of the data $(\pi,\la,\tau)$ which is explicitely defined in \ref{s3ss1sss2} (see Definition \ref{def:w}). The formula can be seen as a generalization of the classical formula (\ref{eqformulaclassica}) to the Rauzy-Veech induction.  We prove the following theorem.

\begin{thm}\label{thm:closure}
Let $\cI$ be any invariant locus contained in some stratum $\cH$ of the moduli space
of translation surfaces. Then
\begin{enumerate}
\item
The Lagrange Spectrum $\cL(\cI)$ is a closed subset of the real line.
\item
The values $a^{-1}(X)$ for $X$ in $PA(\cI)$ are dense in $\cL(\cI)$.
\end{enumerate}
\end{thm}

\subsubsection{Relations with other results in the literature}

In \cite{maucourant}, F. Maucourant treats the closure of Lagrange Spectra and density of values coming from periodic orbits for the geodesic flow on the unitary tangent bundle of a non-compact Riemaniann manifold $M$ with finite volume and sectional curvature not grater then $-1$. Maucourant's result obviously provides a proof of the special case of Theorem \ref{thm:closure} for the Lagrange spectrum $\cL(X)$ of a Veech surface $X$. Nevertheless it cannot be applied in the general case to any spectrum $\cL(\cI)$, indeed strata and their sub-loci are not the tangent bundle of the moduli space $\cM_g$ of Riemann surfaces of genus $g$ (but just components of a stratification) and moreover the Teichm\"uller metric on $\cM_g$ has vanishing curvature.

Maucourant's proof is essentially based on the classical closing Lemma for negatively curved manifolds together with geometric approximation arguments, which also appear in Paulin-Parkonnen \cite{paulin1}. Very recently, Eskin, Mirzhakhani and Rafi proved in \cite{EMR} a generalization of the classical closing Lemma for strata of quadratic differentials (see also \cite{ursula} for another proof of the closing Lemma on strata and \cite{Wright} where the closing lemma is formulated in each invariant locus). Combined with the geometric approach in the spirit of Maucurant/Paulin-Parkonnen \cite{maucourant,paulin1}, the closing Lemma for strata could be used to  provide a different proof of properties (1) and (2) of the classical Lagrange spectrum for the Lagrange spectrum of any invariant locus.

\medskip

Our proof of Theorem \ref{thm:closure} is of combinatorial flavour and is a generalization of the proof given by  Cusick in \cite{cusick} in terms of the Rauzy-Veech induction. Nevertheless, at the combinatorial level, we also essentially prove a closing and shadowing Lemma for the Teichm\"uler flow on closed-invariant sub-loci $\cI$ of strata (see Appendix \ref{loci}). It is remarkable that, even though Rauzy-Veech induction is in principle used to study the whole connected component of a stratum, it follows \emph{a posteriori} that our construction provides closed orbits which live in an invariant locus $\cI$. This follows in virtue of the \emph{local product structure} of invariant loci, recently proved in the already mentioned work by Eskin and Mirzhakhani \cite{eskinmirza}. Thus, our techniques can for example be used to produce pseudo-Anosov diffeomorphism in a given locus by using Rauzy-Veech algorithm.

\subsection{Hall's ray for square-tiled surfaces}

A \emph{square-tiled surface}, also said \emph{origami}, is a translation surface $X$ tiled by copies of the square $[0,1]^2$. Equivalently, $X$ is square-tiled if there exists a ramified covering $\rho:X\to \RR^2/\ZZ^2$, unramified outside $0\in \RR^2/\ZZ^2$ and such that $\rho^\ast(dz)$ is the holomorphic 1-form of $X$. Square-tiled surfaces are Veech's surfaces. Some more detail can be found in  \S \ref{appendixss2}.

In Section \ref{section5} we prove the following extension of the classical formula (\ref{eqformulaclassica}) to square-tiles surfaces (see Theorem \ref{formulasquaretiled}). We need to consider \emph{reduced} origamis, where an origami $X$ is reduced if the lattice spanned by $\hol(X)$ is $\ZZ^2$ (see \S \ref{s5ss0sss1}). Any arithmetic closed Teichm\"uller disc contains reduced origamis. If $X$ is a reduced origami with $N$ squares and $\al=a_0+[a_1,a_2,\dots]$ denotes the continued fraction expansion of $\al$, then (under a technical condition on $\al$ for which we refer to Theorem \ref{formulasquaretiled}) we have
$$
L_X(\al)=
N\cdot\limsup_{n}
\frac
{[a_{n},\dots,a_{1}]+a_{n+1}+[a_{n+2},a_{n+3}\dots]}
{m_X^2(p_n/q_n)},
$$
where $p_n/q_n$ denotes the $n$-th convergent $[a_1,\dots,a_n]$ of $\al$ and the factor $m_X(p_n/q_n)$ is the \emph{multiplicity} of $p_n/q_n$, which is defined in \S \ref{s5ss1sss2}. As an application of this formula, we prove the following result.

\begin{thm}\label{thm:Hall}
Let $X$ be any reduced \emph{square-tiled surface}. Then there exists $r(X)>0$, whose explicit value appears in the statement of Theorem \ref{thm2s5ss2}, such that $\cL(X)$ contains the half-line $[r(X),+\infty)$.
\end{thm}

Obviously, an inclusion between invariant loci $\cI\subset \cI'$ implies the inclusion $\cL(\cI)\subset\cL(\cI')$ for the corresponding Lagrange spectra, thus we have the following criterion for the existence of Hall's ray.

\begin{cor}[Criterion for Hall's ray]\label{cor:criterionhallray}
If an invariant locus $\cI$ contains a square-tiled surface $X$, then its Lagrange Spectrum $\cL(\cI)$ contains an Hall's ray.
\end{cor}

\subsubsection{Relation with other results in the literature}

Let us remark that Teichm\"uller disks of square tiled surfaces are hyperbolic surfaces with constant negative curvature and cusps (see Veech \cite{veechveechcurves}). For any such surface Sheirngorn and Schmidt proved the existence of a Hall ray for the Markoff spectrum of penetration of geodesics into cusps (see \cite{schmidtsheingorn}). One can show that this implies the existence of Hall ray for the set of values
$$
\left\{
\sup_{t>0}\frac{2}{\sys^2(\cF_t Y)}
\textrm{ for }
Y \in \slduer\cdot X
\right\}
$$
where $X$ is a Veech surface. We stress that such result does not imply that there is a Hall ray for the Lagrange spectrum $\cL(X)$. Finally let us mention  that in \cite{paulin2} Paulin et Parkonnen showed the existence of Hall rays for Lagrange spectra of hyperbolic manifolds of negative curvature, but only in dimension $n\geq 3$. Moreover S. A. Roma\~{n}a announced the proof, obtained in his Ph-D thesis, of the existence of open intervals (but not Hall's ray) in the Lagrange Spectrum for a generic small $C^k$-perturbation of an hyperbolic metric on a punctured surface, where $k\geq 2$ (see \cite{ibarra}).

\subsection{Bounded type IETs and penetration estimates}

Another contribution of this paper is to present the proof of a series of equivalent characterizations of bounded-type IETs and translation surfaces whose forward Teichm\"uller orbit is bounded.

First of all, in addition to the equivalences stated by Vorobets' identity (Proposition \ref{prop:vorobetsidentity}), in \S \ref{sec:finitetime} we also provide two other equivalent characterizations. The first is a characterization of bounded type IETs as IETs such that the matrices of a suitable acceleration of the Rauzy-Veech induction have bounded norms (see Corollary \ref{cor:boundedtypeIETs} in Section \ref{sec:finitetime}). This generalizes the classical definition of bounded type real numbers as those $\al\in\RR$ whose continued fraction entries $a_n$ are bounded. The second equivalence characterizes bounded type IETs as those whose orbit under the Rauzy-Veech map stays in a compact set. This is analogous to the fact that orbits of bounded type real numbers under the \emph{Gauss map} stay away from zero.

Finally we establish estimates on excursions to the boundary of moduli space in terms of norms of positive matrices in the Rauzy-Veech induction (see Proposition \ref{prop1s4ss1} and Proposition \ref{prop2s4ss1}). This generalizes the classical relation between the bound on $a_n$ and the depth of penetration to the cusp of the corresponding geodesic. We believe that these quantative  estimates (whose proofs is rather techical and postponed to the Appendix \ref{s7}) will be useful in future applications, since they show that Rauzy-Veech induction, which is a  very well-studied \emph{combinatorial tool} (see the surveys \cite{yoccoz, viana, zorichdue}) can be actually used to deduce properties about the \emph{geometry} (specifically, the depth of penetration to infinity) of a Teichmueller geodesic.

\subsection*{Structure of the paper}

The rest of the paper is arranged as follows. In  \S  \ref{backgroundRazyVeech} we include some technical background material. In particular, we define  the Rauzy-Veech induction for IETs and for translation surfaces (in \S \ref{backgroundss2} and \S \ref{backgroundss3} respectively) and  recall the Veech zippered rectangles construction in \S \ref{backgroundss1}.

In  \S  \ref{section2} we prove Vorobet's identity (Proposition \ref{prop:vorobetsidentity}). We also introduce the notions of reduced periods for a translation surface (\S~\ref{sec:reducedperiods}) and reduced triples for an IET (\S~\ref{sec:reducedtriples}) and explain their connection via the zippered rectangles contruction (\S~\ref{sec:periods_vs_triples}).

In \S \ref{s3} we prove the formula which enables us to compute the
Lagrange spectrum $\cL(\cI)$ of any invariant locus via the Rauzy-Veech algorithm  (Theorem \ref{renormalizedformula}) and then we prove Theorem \ref{thm:closure} (closure of $\cL(\cI)$ and density of values of periodic orbits).
First, in \S \ref{s3ss1(formulastrata)} we prove  Theorem \ref{renormalizedformula}.
In \S \ref{sec:finitetime} we state quantitative relation between elements $a^{-1}(X)$ in the Lagrange spectra $\cL(\cI)$ and the maximal asymptotic size of positive matrices occurring in the Rauzy-Veech expansion which encodes $X$ (see Theorem \ref{thm1s4}). The proof of Theorem \ref{thm1s4} is based on estimations in finite (uniformly bounded) time, namely Proposition \ref{prop1s4ss1} and Proposition \ref{prop2s4ss1}, whose proof is postponed to the Appendix \ref{s7}. As consequence of the estimations in \S \ref{sec:finitetime}, we prove the qualitative equivalence between several notion of bounded type IET's and translation surfaces. In
\ref{sec:subshifts} we define a sub-shift associated to a given Rauzy class, whose alphabet is the set of \emph{positive paths} in the Rauzy diagram associated to the connected component $\cC$ which contains $\cI$. The language of sub-shifts provides a convenient formalism to prove Theorem \ref{thm:closure}. The proof of Theorem \ref{thm:closure} is given in \S~\ref{subsection:closureanddensity}.

In \S \ref{section5} we prove Theorem \ref{thm:Hall}. The main tool is the formula appearing in Theorem \ref{formulasquaretiled}, which enables us to compute the
Lagrange spectrum $\cL(X)$ of a reduced square-tiled surface $X$ in terms of a \emph{skew-product} over the classical continued fraction. Applying the formula, and following the classical argument given by Hall, the existence of the Hall ray for $\cL(X)$ is established in Theorem \ref{thm2s5ss2}.

In the Appendix \ref{s7} we give the proof of Proposition \ref{prop1s4ss1} and Proposition \ref{prop2s4ss1}, which provide quantitative finite time relations between the norms of the Rauzy-Veech matrices  and two geometrical quantities, namely $w(\pi,\la,\tau)$ appearing in the formula of Proposition \ref{renormalizedformula} and the \emph{distortion} $\De(T)$ of the IET associated to $(\pi,\la,\tau)$ (see \S \ref{s4ss1}).

Finally, in Appendix \ref{loci} we describe the structure of invariant loci in the parameter space of Rauzy-Veech induction and give the proof of the combinatorial closing and shadowing lemmas for invariant loci (Propositions \ref{prop:appendixclosinglemma} and \ref{prop:appendixshadowinglemma}) which are used in the proof of Theorem \ref{thm:closure}.

\subsection*{Acknowledgements}

We would like to thank U.~Hamenst{\"a}dt, E.~Lanneau, S. Leli\`evre, C.~Matheus, F.~Paulin and J.-C.~Yoccoz for useful discussions. We would also like to thank the Hausdorff Institute for Mathematics for the hospitality during the program \emph{Geometry and Dynamics of Teichm\"uller Spaces}, where this project started. Some of the research visits which made this collaboration possible were supported by the EPSRC Grant EP/I019030/1 and the ANR Project GeoDyM. P. Hubert is partially supported by the ANR Project GeoDyM. C. Ulcigrai is partially supported by EPSRC Grant EP/I019030/1. C. Ulcigrai would also like to thank the FIM Institute for Mathematical Research at ETH Zurich where part of this work was completed.

\section{Background: zippered rectangles and Rauzy-Veech induction}\label{backgroundRazyVeech}

Fix an alphabet $\cA$ with $d\geq 2$ letters and call $\mathfrak{S}(\cA)$ the set of all admissible combinatorial data $\pi$ over $\cA$. If $T$ is an IET determined by the combinatorial-length data $(\pi,\la)$, we assume that it acts on the interval $I=(0,\sum_\chi\la_\chi)$. If $\al$ and $\be$ are letters in $\cA$ with $\pi^t(\al)>1$ and $\pi^b(\be)>1$ we set
respectively
$$
u^t_\al:=
\sum_{\pi^t(\chi)<\pi^t(\al)}\la_\chi
\textrm{ and }
u^b_\be:=
\sum_{\pi^b(\chi)<\pi^b(\be)}\la_\chi.
$$
For any such $\al$ and $\be$ the points $u^t_\al$ and $u^b_\be$ are the left endpoints of the subintervals $I^t_\al$ and $I^b_\be$. In general $T$ is not continuous at $u^t_\al$ and $T^{-1}$ is not continuous at $u^b_\be$. We will often consider the sup norm of length data, that is $\|\la\|:=\sum_\chi\la_\chi$.

\subsection{Rauzy induction for IETs}\label{backgroundss2}

In this paragraph we give a brief survey of the basic properties of the Rauzy induction for IETs, following \cite{agy, mmy, yoccoz}. For a comprehensive introduction, we refer the reader to the lecture notes \cite{yoccoz}.

\subsubsection{Rauzy map and Rauzy elementary operations}\label{backgroundss2sss1}

Let $T$ be an IET defined by combinatorial-length data $(\pi,\la)$, acting on the interval $I=(0,\sum_\chi\la_\chi)$. Consider the variable $\epsilon\in\{t,b\}$, where the letter $t$ stands for \emph{top} and the letter $b$ for \emph{bottom}. If $\epsilon=t$ we set $1-\epsilon:=b$ and if $\epsilon=b$ we set $1-\epsilon:=t$. Let $\al_t$ and $\al_b$ be the letters in $\cA$ such that respectively $\pi^t(\al_t)=d$ and $\pi^b(\al_b)=d$. The rightmost singularity of $T$ is therefore $u^t_{\al_t}$ and the rightmost singularity of $T^{-1}$ is $u^b_{\al_b}$. When the condition
\begin{equation}\label{eq1backgroundss2}
u^t_{\al_t}\not= u^b_{\al_b}
\end{equation}
is satisfied, we say that $T$ is of \emph{type}
$\epsilon\in \{t,b\}$ if
$$
u^\epsilon_{\al_\epsilon}<
u^{1-\epsilon}_{\al_{1-\epsilon}}.
$$
Consider the subinterval
$
\widetilde{I}:=
I\cap (0,u^{1-\epsilon}_{\al_{1-\epsilon}})
$
of $I$ and define
$
\widetilde{T}:\widetilde{I}\to\widetilde{I}
$
as the first return map of $T$ to $\widetilde{I}$, which is obviously an IET. The combinatorial datum
$
\widetilde{\pi}=
(\widetilde{\pi}^t,\widetilde{\pi}^b)
$
of $\widetilde{T}$ is given by:
\begin{equation}\label{eq2backgroundss2}
\begin{array}{lllll}
\widetilde{\pi}^\epsilon(\al)=
\pi^\epsilon(\al)
\textrm{ for any } \al \in \cA\\
\widetilde{\pi}^{1-\epsilon}(\al)=
\pi^{1-\epsilon}(\al)
\textrm{ if }
\pi^{1-\epsilon}(\al)\leq \pi^{1-\epsilon}(\al_{\epsilon})\\
\widetilde{\pi}^{1-\epsilon}(\al_{1-\epsilon})=
\pi^{1-\epsilon}(\al_{\epsilon})+1\\
\widetilde{\pi}^{1-\epsilon}(\al)=
\pi^{1-\epsilon}(\al)+1
\textrm{ if }
\pi^{1-\epsilon}(\al_{\epsilon})<\pi^{1-\epsilon}(\al)< d.
\end{array}
\end{equation}
The length datum $\widetilde{\la}$ of $\widetilde{T}$ is given by:
\begin{equation}\label{eq3backgroundss2}
\begin{array}{ll}
\widetilde{\la}_{\al}=
\la_{\al}
\textrm{ if }
\al\neq\al_{\epsilon}\\
\widetilde{\la}_{\al_{\epsilon}}=
\la_{\al_{\epsilon}}-\la_{\al_{1-\epsilon}}.
\end{array}
\end{equation}
When $T=(\pi,\la)$ satisfies condition (\ref{eq1backgroundss2}), Equations (\ref{eq2backgroundss2}) and (\ref{eq3backgroundss2}) define a map $T\mapsto Q(T):=\widetilde{T}$, called \emph{Rauzy map}. Introduce the operations $R^t$ and $R^b$ from $\mathfrak{S}(\cA)$ to itself setting $R^\epsilon(\pi):=\widetilde{\pi}$, where $\epsilon$ and $\pi$ are respectively the type and the combinatorial datum of $T$ and $\widetilde{\pi}$ is the combinatorial datum of $\widetilde{T}$. It is easy to check that if $\pi$ is admissible then both $R^{t}(\pi)$ and $R^{b}(\pi)$ are admissible. The maps $R^t$ and $R^b$ are called the \emph{Rauzy elementary operations}.

A \emph{Rauzy class} is a minimal non-empty subset $\cR$ of $\mathfrak{S}(\cA)$ which is invariant under $R^t$ and $R^b$. A \emph{Rauzy diagram} is a connected oriented graph $\cD$ whose vertexes are the elements of $\cR$ and whose oriented arcs, or \emph{arrows}, correspond to Rauzy elementary operations
$
\pi\mapsto R^\epsilon(\pi)
$.
An arrow corresponding to $R^t$ is called a \emph{top} arrow and we say that $\al_t$ is its winner and $\al_b$ is its loser. Conversely an arrow corresponding to $R^b$ is called a \emph{bottom} arrow and we say that $\al_t$ is its loser and $\al_b$ is its winner.

A concatenation of $r$ compatible arrows $\ga_1,\dots,\ga_r$ in a Rauzy diagram is called a \emph{Rauzy path} and is denoted $\ga=\ga_1\ast\dots\ast\ga_r$. The set of all Rauzy paths connecting elements of $\cR$ is denoted $\Pi
(\cR)$. If a path $\ga$ is concatenation of $r$ simple arrows, we say that $\ga$ has length $r$. Length one paths are arrows, elements of $\cR$ are identified with trivial (that is length-zero) paths.

\subsubsection{Linear action and iterations of the Rauzy map}\label{backgroundss2sss2}

Let $\{e_\xi\}_{\xi\in\cA}$ be the canonical basis of
$\RR^\cA$. For any Rauzy class $\cR$ over $\cA$ and any path $\ga\in\Pi(\cR)$ define a linear map $B_\ga\in\sldz$ as follows. If $\ga$ is trivial then $B_\ga:=id$. If $\ga$ is an arrow with winner $\al$ and loser $\be$ set
$$
B_\ga e_\al=e_\al+e_\be
\textrm{ and }
B_\ga e_\xi=e_\xi
\textrm{ for }
\xi\not=\al.
$$
Then extend the definition to paths in $\Pi(\cR)$ so that $B_{\ga_1\ga_2}=B_{\ga_2}B_{\ga_1}$. For a combinatorial datum $\pi\in\cR$ and a Rauzy path $\ga\in\Pi(\cR)$ starting at $\pi$ define the simplicial sub-cone $\De_\ga\subset\RR^\cA_+$
by
$$
\De_\ga=^t\!B_\ga(\RR_{+}^{\cA}),
$$
where  $^t\!B_\ga$ is the transpose matrix of $B_\ga$. Finally, denoting by  $\vec{1}$ the vector of $\NN^\cA$ whose entries are all equal to $1$, for any $\ga\in\Pi(\cR)$ define a vector $q^\ga\in\NN^\cA$ by
$$
q^\ga:=B_\ga\vec{1}.
$$
If $T=(\pi,\la)$ is an IET in $\cR\times\RR^\cA_+ $ admitting $r$ elementary steps $\ga_1,\dots,\ga_r$ of the Rauzy map $Q$ and $\ga=\ga_1\ast\dots\ast\ga_r$ is the corresponding Rauzy path, we also write $q^{(r)}$ instead of $q^\ga$. For such $T$, denote $T^{(r)}:=Q^{r}(T)$ the $r$-th step of the induction, whose combinatorial-length data are $(\pi^{(r)},\la^{(r)})$ . The interval where $T^{(r)}$ acts is
$
I^{(r)}=\big(0,\sum_\chi\la^{(r)}_\chi\big)
$,
which is a sub-interval of $I$ with the same left endpoint. For $\al\in\cA$ call $u^{(r),t}_\al$ and $u^{(r),t}_\al$ the singularity respectively of $T^{(r)}$ and of $(T^{(r)})^{-1}$ corresponding to the letter $\al$. The length datum of $T^{(r)}$ is given by the following Lemma (see \cite{yoccoz} for a proof).

\begin{lem}\label{lem1backgroundss2}
Fix combinatorial data $\pi$ and $\pi'$ in $\cR$ and let $\ga\in\Pi(\cR)$ be a path starting at $\pi$ and ending in $\pi'$ with length $r$. Then the $r$-th iterate of $Q$ is a linear isomorphism defined on $\pi\times\De_\ga$ with values in $\pi'\times\RR^\cA_+$ and the length datum $\la^{(r)}$ of $T^{(r)}$ is given by the formula
$$
\la^{(r)}=^t\!B_\ga^{-1}\la.
$$
\end{lem}

Condition (\ref{eq1backgroundss2}) corresponds to an union of hyperplanes in $\cR\times\RR^\cA_+$. Similarly, for any $r$ the iterate $Q^r$ is defined outside a finite union of linear spaces, therefore IETs admitting infinitely many iterations of the map $Q$ form a set with full Lebesgue measure. The complement of such set corresponds to those $T$ such that $T^{(r)}=Q^{r}(T)$ eventually does not satisfy condition (\ref{eq1backgroundss2}), that is to those $T$ such that the algorithm stops. The following combinatorial characterization holds (see \cite{yoccoz} for a proof).

\begin{lem}
$T$ admits infinitely many steps of the Rauzy induction if any only if it does
not have connections.
\end{lem}

\subsubsection{Return times}\label{backgroundss2sss3}

Let $\pi\in\cR$ and $\ga\in\Pi(\cR)$ be a path of length $r$. Let $T$ be an IET in $\{\pi\}\times\De_\ga$ and $T^{(r)}$ be the $r$-th step of the Rauzy induction applied to $T$, whose combinatorial-length data are $(\pi^{(r)},\la^{(r)})$. The sub-interval of $I^{(r)}$ where $T^{(r)}$ (respectively $(T^{(r)})^{-1}$) acts as a translation are
$
I^{(r),t}_{\chi}:=
(u^{(r),t}_\chi,u^{(r),t}_\chi+\la^{(r)}_\chi)
$
(respectively
$
I^{(r),b}_{\chi}:=
(u^{(r),b}_\chi,u^{(r),b}_\chi+\la^{(r)}_\chi)
$)
for $\chi\in\cA$. For such $\chi$ and $r$ let $R_r(\chi)$ be the minimal positive integer $k$ such that $T^k(I^{(r),t}_\chi)$ belongs to $I^{(r)}$. Observe that $R_r(\chi)$ is also equal to the minimal positive integer $k$ such that $T^{-k}(I^{(r),b}_\chi)$ belongs to $I^{(r)}$. For a matrix $A$ let $[A]_{\al\be}$ be the entry of $A$ in row $\al$ and column $\be$. A proof of the following Lemma can be found in \cite{yoccoz}.

\begin{lem}\label{lem3backgroundss2}
Let $T$ be an IET admitting $r$ elementary steps $T\mapsto T^{(r)}$ of the Rauzy map, represented by the path $\ga\in\Pi(\cR)$. For any $\al\in\cA$ we have
$$
R_r(\al)=q^{(r)}_\al=
\sum_{\be\in\cA}[B_\ga]_{\al\be},
$$
that is $q^{(r)}_\al$ is the return time to $I^{(r)}$ of the sub-interval $I^{(r),t}_\al$. More precisely, for any $\al,\be,\in\cA$ the entry $[B_\ga]_{\al\be}$ of $B_{\ga}$ is equal to the cardinality of the following two sets:
\begin{enumerate}
\item
The integers $k$ with $0\leq k<R_r(\al)$ such that $T^k(I^{(r),t}_\al)\subset  I^t_{\be}$.
\item
The integers $k$ with $0\leq k<R_r(\al)$ such that
$T^{-k}(I^{(r),b}_\al)\subset I^b_\be$.
\end{enumerate}
\end{lem}

\subsection{Zippered rectangles construction}\label{backgroundss1}

Here we describe Veech's \emph{zippered rectangles construction}, following \cite{mmy}. Fix $\pi\in\mathfrak{S}(\cA)$ and introduce the labelling $\al(1),\dots,\al(d)$ of the letters in $\cA$, according to their order in $\pi^t$. Similarly, introduce the labelling $\be(1),\dots,\be(d)$ of the letters in $\cA$ according to their order in $\pi^b$.

A \emph{suspension datum} for $\pi$ is a vector
$\tau$ in $\RR^\cA$ such that for any $1\leq k\leq d-1$
%%%CHECK it was d but should be d-1 right?(Corinna)
%%%right!(Luca)
we have
$$
\sum_{j\leq k} \tau_{\al(j)}>0
\textrm{ and }
\sum_{j\leq k} \tau_{\be(j)}<0.
$$
Let $\Theta_\pi$ be the open sub-cone of $\RR^\cA$ of suspension data $\tau$ for $\pi$. Observe that the vector $\tau$ with coordinates $\tau_\xi:=\pi^b(\xi)-\pi^t(\xi)$ satisfies the inequalities above, hence $\Theta_\pi$ is not empty.

Consider $\la\in\RR_+^\cA$ and $\tau\in\Th_\pi$ and define the complex vector
$
\zeta=\lambda +i\tau \in \CC^\cA
$,
then for any $1\leq k\leq d-1$ introduce the complex numbers
$$
\xi_{\al(k)}^t:=\sum_{j\leq k}\zeta_{\al(j)}
\textrm{ and }
\xi_{\be(k)}^b:=\sum_{j\leq k}\zeta_{\be(j)}.
$$

\subsubsection{The height function $h$}\label{backgroundss1sss1}

Let $T$ be the IET defined by the data $(\pi,\la)$, acting on the interval $I=(0,\sum_\al\la_\al)$. Observe that the singularities $u^t_\be$ and $u^b_\al$  of $T$ and $T^{-1}$ satisfy $u_\al^t=\Re(\xi_\al^t)$ and $u^b_\be=\Re(\xi_\be^b)$. The condition $\tau\in\Th_\pi$ is equivalent to $\Im(\xi_{\al(k)}^t)>0$ and $\Im(\xi_{\be(k)}^b)<0$ for any $1\leq k\leq d-1$. Moreover admissibility of $\pi$ implies $\al(1)\not=\be(1)$ and $\al(d)\not=\be(d)$, hence the vector $h=h(\pi,\tau)$   defined by setting for $1\leq k\leq d$
$$
h_{\al(k)}:=
\sum_{j\leq k} \tau_{\al(j)}-
\sum_{j\leq k} \tau_{\be(j)}
$$
is a positive vector, that is $h_\al>0$ for any $\al$.
It is useful to consider $h$ as a piecewise constant function $h:I\to\RR_{+}$, with constant value $h_\al$ on each sub-interval
$I^t_\al$ of $I$.
%We call it a \emph{roof function} for $T$. The reason is that the \emph{suspension flow} $\phi^t$ over $T$ under the roof function $h$ will be identified with the vertical flow of the translation surface $X=X(\pi,\la,\tau)$ that is going to be constructed in \S \ref{backgroundss1sss2}.

Observe that for any $k$ with $1<k\leq d$ the suspension conditions imply
$
h_{\al(k)}\geq \tau_{\al(1)}+\dots+\tau_{\al(k-1)}>0
$
and
$
-h_{\al(k)}\leq \tau_{\be(1)}+\dots+\tau_{\be(k-1)}<0
$.
Combining this remark with the admissibility of $\pi$ it is not difficult to get the following estimate for the height function $h$.

\begin{lem}\label{lem1backgroundss1}
We have
$$
\|\tau\|_\infty\leq \|h\|_\infty.
$$
\end{lem}

\subsubsection{Construction of a translation surface}\label{backgroundss1sss2}

We construct a translation surface $X=X(\pi,\la,\tau)$ whose vertical flow is the suspension flow over $T$ with roof
function $h$, that is the unit vertical flow $\stackrel{\cdot}{y}=1$ with respect to the equivalence relation $(x, y)\sim (T(x), y - h(x))$.  Embed $I$ in the complex plane, identifying it with $(0,\sum_\al\la_\al)\times \{0\}$ and accordingly denote elements in $\CC\cap I$ just with their real coordinate. Define $2d$ open rectangles in $\CC$ setting
$$
R_\al^t:=
(u^t_\al,u^t_\al+\la_\al)\times(0,h_\al)
\textrm{ and }
R_\be^b:=
(u^b_\be,u^b_\be+\la_\be)\times(-h_\be,0).
$$
Consider the \emph{translation datum} $\te$, that is the vector in $\CC^\cA$ defined by
$
\theta_\al:=\xi^b_\al-\xi^t_\al
$
for any $\al$. In order to get a surface we "zip" together these rectangles by glueing their boundaries according to the identifications described below.

\begin{enumerate}
\item
For each $\al$ the rectangle $R_\al^t$ is equivalent to the
rectangle $R_\al^b$ via the translation by $\te_\al$.
\item
For each $k>1$ we paste together $R^t_{\al(k-1)}$ and
$R^t_{\al(k)}$
along the common vertical open segment connecting $u^t_{\al(k)}$ to $\xi^t_{\al(k)}$.
\item
For each $k$ with $k>1$ we paste together $R^b_{\be(k-1)}$ and
$R^b_{\be(k)}$ along the vertical open segment connecting $\xi^b_{\be(k)}$ to $u^b_{\be(k)}$.
\item
For any $\al$ we paste $R_\al^t$ to $I$ along its lower
horizontal boundary segment
$
(u^t_\al,u^t_\al+\la_\al)\times \{0\}
$.
\item
For any $\be$ we paste $R_\be^b$ to $I$ along its upper
horizontal boundary segment
$
(u^b_\be,u^b_\be+\la_\be)\times\{0\}
$.
\item
If $\sum_\al\tau_\al\geq 0$ identify the vertical segment
connecting $u^t_{\be(d)+1}+ih_{\be(d)}$
to $\xi^t_{\be(d)+1}$ and the vertical segment connecting
$\sum_\al\la_\al$ to $\sum_\al\zeta_\al$. This is possible because
$h_{\be(d)}=
\Im(\xi^t_{\be(d)})-\Im(\xi^b_{\be(d)})=
\Im(\xi^t_{\be(d)+1})-\sum_\al\tau_\al\leq
\Im(\xi^t_{\be(d)+1})
$.
If $\sum_\al\tau_\al<0$ identify the vertical segment
connecting $\xi^b_{\al(d)+1}$ to $u^b_{\al(d)+1}-ih_{\al(d)}$ and the vertical segment connecting $\sum_\al\zeta_\al$ to
$\sum_\al\la_\al$. This is possible because
$
-h_{\al(d)}=
-\sum_\al\tau_\al+\Im(\xi^b_{\pi(d)+1})>
\Im(\xi^b_{\al(d)+1})
$.
\item
Finally we add the points $\xi_{\al(k)}^t$ and $\xi_{\be(k)}^b$ for $1\leq k\leq d-1$.
\end{enumerate}

\begin{prop}\label{prop1backgroundss1}
For any connected component $\cC$ of any stratum of translation surfaces there exists a Rauzy class $\cR$ such that the following holds.
\begin{enumerate}
\item
For any $\pi\in\cR$, any $\la\in\RR^\cA_+$ and any $\tau\in\Th_\pi$ the surface $X(\pi,\la,\tau)$ given by the zippered rectangles construction belongs to $\cC$.
\item
If $X$ is a translation surface in $\cC$ such that do not exist data $\pi\in\cR$, $\la\in\RR^\cA_+$ and $\tau\in\Th_\pi$ with $X=X(\pi,\la,\tau)$, then $X$ has a vertical and an horizontal saddle connection. Therefore the zippered rectangles construction fills a subset of $\cC$ of co-dimension 2.
\item
In particular, if $s(X)>0$ (recall the definition \ref{eq:definitions}) then there exist data $\pi$, $\la$ and $\tau$ such that $X=X(\pi,\la,\tau)$.
\end{enumerate}
\end{prop}

\subsection{Rauzy-Veech induction}\label{backgroundss3}

Consider combinatorial-length suspension data $(\pi,\la,\tau)$. When the data $(\pi,\la)$ satisfy condition (\ref{eq2backgroundss2}), we define a new suspension datum $\widetilde{\tau}$ extending Equation (\ref{eq3backgroundss2}) to $\tau$ by
\beq\label{eq1backgroundss3}
\begin{array}{ll}
\widetilde{\tau}_{\al}=\tau_{\al}
\textrm{ if }
\al\neq\al_{\epsilon}\\
\widetilde{\tau}_{\al_{\epsilon}}=
\tau_{\al_{\epsilon}}-\tau_{\al_{1-\epsilon}},
\end{array}
\eeq
where the letters $\al_{\epsilon}$ and $\al_{1-\epsilon}$ are defined at the beginning of \S \ref{backgroundss2sss1}. Equation (\ref{eq1backgroundss3}), together with Equations (\ref{eq2backgroundss2}) and (\ref{eq3backgroundss2}), define a map
$$
(\pi,\la,\tau)\mapsto
\widehat{Q}(\pi,\la,\tau):=
(\widetilde{\pi},\widetilde{\la},\widetilde{\tau}),
$$
called \emph{Rauzy-Veech map}. Let $\ga$ be a Rauzy path starting at $\pi$ and ending in $\pi'$ and let
$B_{\ga}$ be the associated matrix. Define the cone
$
\Th_{\ga}:=^{t}B^{-1}_{\ga}\Th_{\pi}
$.

\begin{lem}
If $\ga$ starts at $\pi$ and ends in $\pi'$, then
$\Th_{\ga}$ is a subcone of $\Th_{\pi'}$.
\end{lem}

\begin{dimo}
It is enough to prove the Lemma for a simple arrow. In this case, if $\widetilde{\tau}$ is obtained from a suspension datum $\tau$ in $\Th_\pi$ via Equation (\ref{eq1backgroundss3}), then is easy to check
that if $\widetilde{\tau}$ is a suspension datum for $\pi'$.
\end{dimo}

The iteration of the map $\widehat{Q}$ satisfies the same formalism as the Rauzy map $Q$. More precisely, for any Rauzy path $\ga$ as above, we have a linear homeomorphism $\widehat{Q}_{\ga}$ from
$
\{\pi\}\times\De_\ga\times\Th_\pi
$
with values in
$\{\pi'\}\times\RR_+^\cA\times\Th_\ga
$
defined by
$$
\widehat{Q}_\ga(\pi,\la,\tau):=
(\pi',^tB^{-1}_\ga\la,^tB^{-1}_\ga\tau).
$$

\begin{lem}
Let $\ga$ be a Rauzy path starting at $\pi$ and ending in $\pi'$. Consider
$\la$ in $\De_{\ga}$ and $\tau$ in $\Th_{\pi}$, then
set
$\la':=^{t}B^{-1}_{\ga}\la$
and
$\tau':=^{t}B^{-1}_{\ga}\tau$.
Then the zippered rectangles construction applies to the data $(\pi',\la',\tau')$, moreover
the translation surfaces $X(\pi,\la,\tau)$ and $X(\pi',\la',\tau')$ define the
same element in the moduli space.
\end{lem}

\section{On Vorobet's identity}\label{section2}

\subsection{Reduced periods for translation surfaces}\label{sec:reducedperiods}

Fix a translation surface $X$ and consider a saddle connection $\ga$ for $X$. The \emph{standard orientation} of $\ga$ is the orientation induced by any smooth parametrization $t\mapsto\ga(t)$ of $\ga$, where $t\in[0,1]$, such that for any $t$ we have
$$
\angolo
\bigg(
\frac{d}{dt}\ga(t),\partial_y
\bigg)<
\frac{\pi}{2}.
$$
In this case we say that $\ga$ \emph{starts at} $p_{out}$ and \emph{ends in} $p_{in}$ if we have respectively $\ga(0)=p_{out}$ and $\ga(1)=p_{in}$. Observe that if $\ga$ has the standard orientation, then its period $v=\int_\ga w$ belongs to the upper half plane $\HH$. For such a period, let $R_v$ be the rectangle in $\HH$ whose diagonal is $v$.

\begin{defi} We say that a period $v\in\hol(X)$ is a \emph{reduced period} if there exists an immersion $\rho:R_v\to X$ isometric with respect to the flat metric of $X$ such that $t\mapsto\rho(tv)$ for $0\leq y \leq 1$ is a standard parametrization of $\ga$ and $\rho(R_v)$ does not contain other conical singularities in its interior (the endpoints of $\ga$ of course belong to the boundary of $\rho(R_v)$).
\end{defi}
Note that in general $\rho$ cannot be an embedding, since it may not be injective.

\begin{lem}\label{lem1s2ss2}
For any $X$, the value $a(X)$ can be computed taking the $\liminf$ just on reduced periods.
\end{lem}

\begin{dimo}
Just observe that if a period $v\in\hol(X)$ is not reduced then there exists an other reduced period $v'$ in $\hol(X)$ with $|\Re(v')|<|\Re(v)|$ and
$|\Im(v')|<|\Im(v)|$, and hence $\area(v')<\area(v)$. Then the Lemma follows because $\hol(X)$ is a discrete subset of the plane.
\end{dimo}

\subsection{Reduced triples for IETs}\label{sec:reducedtriples}
We now define the notion of reduced triples for an IET, which is the combinatorial counterpart at the level of IETs of the notion of reduced periods given in S\ref{sec:reducedperiods} above.
Let $T$ be an IET with combinatorial and length data $(\pi,\la)$, acting on the interval $I=(0,\sum_\chi\la_\chi)$. Recall that if $\al$ and $\be$ are two letters in $\cA$ with $\pi^t(\al)>1$ and $\pi^b(\be)>1$ the points $u^t_\al$ and $u^b_\be$ denote the left endpoint of the subintervals $I^t_\al$ and $I^b_\be$ respectively and in general $T$ is not continuous at $u^t_\al$ and $T^{-1}$ is not continuous at $u^b_\be$. A \emph{connection} for $T$ is a triple $(\be,\al,n)$ with $\pi^b(\be)>1$, $\pi^t(\al)>1$ and $n\in\NN$ such that $T^n(u^b_\be)=u^t_\al$. In particular one can check that if $T$ has no connections then $\pi$ is admissible. Consider a triple $(\be,\al,n)$ with $n$ in $\NN$, $\pi^b(\be)>1$ and
$\pi^t(\al)>1$ and suppose that it is not a connection for $T$. Then we denote $I(\be,\al,n)$ the open subinterval of $I$
whose endpoints are $T^n(u^b_\be)$ and $u^t_\al$. %%%, for both the two possible reciprocal orders between them.

\begin{defi}\label{defreducedtriple}
We say that $(\be,\al,n)$ is a \emph{reduced triple} for $T$
if it is not a connection for $T$ and moreover for any $k$ in $\{0,\dots,n\}$ the pre-image $T^{-k}(I(\be,\al,n))$ does not contain any singularity of $T$ or of $T^{-1}$.
\end{defi}

\subsubsection{An alternative way to compute Lagrange Spectra}\label{sec:reducedtriplesLagrange}
%\subsubsection{An alternative way to compute Lagrange Spectra}

We introduce the following asymptotic quantity, which gives an alternative way to compute Lagrange Spectra.
$$
l(T):=
\frac{1}{\|\la\|}\cdot
\liminf_{n\to\infty}
n\cdot
\min_{\be,\al}|T^{n}u^{b}_{\be}-u^{t}_{\al}|.
$$

\begin{prop}\label{prop1s2ss1}
Let $T$ be an IET. We have $l(T)=\cE(T)$.
\end{prop}

The following Lemmas are useful in the proof of Proposition \ref{prop1s2ss1}.

\begin{lem}\label{lem1s2ss1}
Let $(\be,\al,n)$ be a triple which is not reduced for
$T$. Then there exists a triple
$(\be',\al',m)$
with $m< n$ which is
reduced for $T$ and such that
$$
|T^m(u^b_{\be'})-u^t_{\al'}|<
|T^n(u^b_\be)-u^t_\al|.
$$
\end{lem}

\begin{dimo}
Let $(\be,\al,n)$ be as in the statement and assume without loss of generality that $T^n(u^b_\be)<u^t_\al$ and that $I(\be,\al,n)$ does not contain other singularities. Since $(\be,\al,n)$ is not reduced for $T$, then there exists $k$ with
$1\leq k\leq n$ such that $T^{-k}I(\be,\al,n)$
contains either a singularity for $T$ or a singularity for $T^{-1}$. Let $k$ be minimal satisfying the property above. Minimality implies that $T^{-k}I(\be,\al,n)$ contains a singularity of $T$ (observe that singularities of $T^{-1}$ are the image of singularities of $T$). Calling $u^t_\chi$ such singularity, we have
$
T^{n-k}(u^b_\be)<u^t_\chi<T^{-k}(u^t_\al)
$,
hence, since $T^k$ is an isometry on $I(\beta, \alpha, n)$, the triple $(\be,\chi,n-k)$ satisfies
$
|T^{n-k}(u^b_\be)-u^t_\chi|< | T^{n-k}(u^b_\be)- T^{-k}(u^t_\al)| =
|T^{n}(u^b_\be)-u^t_\al|
$.
Moreover we can assume that the interval $(T^{n-k}(u^b_\be),u^{t}_{\chi})$
does not contains other singularities. If the triple $(\be,\chi,n-k)$ is reduced
then the Lemma is proved, otherwise we replace $(\be,\al,n)$ by $(\be,\chi,n-k)$
and we repeat the procedure above.
Such procedure admits at most $n$ steps, hence repeating it finitely many
times we end up with a triple
$(\be',\al',m)$ as in the statement.
\end{dimo}

The following Lemma is the analogous for reduced triples of Lemma \ref{lem1s2ss2} for reduced periods.
\begin{lem}\label{lem2s2ss1}
Let $T$ be any IET. Then $l(T)$ can be computed taking the $\liminf$ on reduced triples only.
\end{lem}

\begin{dimo}
Observe that for triples $(\be,\al,n)$ with $n$ bounded the values
$
|T^n(u^b_\be)-u^t_\al|
$
are bounded from below by a positive constant. Thus the Lemma is an immediate consequence of Lemma \ref{lem1s2ss1}.
\end{dimo}

Recall that a singularity of $T^n$ has the form $T^{-i}(u^t_\al)$, where $\pi^t(\al)>1$ and $0\leq i< n$. Recall also that $u^b_\al=T(u^t_\al)$. Therefore the following Lemma is obvious.

\begin{lem}\label{lem3s2ss1}
For any $n$ there exists $m=m(n)$ with $m<n$ and a triple $(\be,\al,m)$ reduced
for $T$ such that
$$
\cE_{n}(T)=
|T^{m}(u^{b}_{\be})-u^{t}_{\al}|.
$$
\end{lem}

\subsubsection{Proof of Proposition \ref{prop1s2ss1}}

Fix $N>0$ and $\epsilon>0$. Consider $n$ in $\NN$ such that
$n\cE_{n}(T)<\cE(T)+\epsilon$ and $m(n)>N$, where $m=m(n)$ is given by Lemma \ref{lem3s2ss1}. Let $(\be,\al,m)$ be the triple reduced for $T$ associated to
$n$ by the same Lemma. We have
$$
m\cdot|T^{m}(u^{b}_{\be})-u^{t}_{\al}|=
\frac{m}{n}n\cdot\cE_{n}(T)\leq
\cE(T)+\epsilon,
$$
hence it follows that $l(T)\leq \cE(T)$, since $m$ is arbitrarily large. To prove the reverse inequality fix again $N>0$ and $\epsilon>0$. Consider a
triple $(\be,\al,n)$ reduced for $T$ such that $n>N$ and
$
n\cdot
|T^{n}(u^{b}_{\be})-u^{t}_{\al}|
<l(T)+\epsilon
$.
Observe that such triple exists according to Lemma \ref{lem2s2ss1}. Since $(\be,\al,n)$ is reduced then $T^{-(n+1)}$ is continuous on $I(\be,\al,n)$.
Thus, recalling that $u^t_\be=T^{-1}(u^b_\be)$ we have
$$
n\cdot\cE_{n+1}(T)\leq
n\cdot|u^{t}_{\be}-T^{-n-1}(u^{t}_{\al})|= %%%CHECK this was an -m-1 but it should be -n-1 right?
n\cdot
|T^{n}(u^{b}_{\be})-u^{t}_{\al}|\leq
l(T)+\epsilon.
$$
It follows that $\cE(T)\leq l(T)$ and the Proposition is proved.

\subsection{Correspondence between triples and periods}\label{sec:periods_vs_triples}

In this paragraph we use Veech's zippered rectangles construction, following the notation introduced in \S \ref{backgroundss1}, to explain the connection between periods (see \S\ref{sec:reducedperiods}) and triples (see \S\ref{sec:reducedtriples}). Consider combinatorial, length and suspension data $(\pi,\la,\tau)$, let $X$ be the underlying translation surface and $\phi^t$ be the vertical flow of $X$. Let $T$ be the IET corresponding to $(\pi,\la)$ and acting on the interval $I=(0,\sum_\al\la_\al)$, which is naturally embedded in $X$ along the horizontal direction. Remark that in this setting $\phi^t$ corresponds to the suspension flow over $T$ under the roof function $h=h(\pi,\tau)$ defined in \S \ref{backgroundss1}.

The singularities of $T$ lie on the vertical negative time
orbits  of the conical singularities of
$X$ and can be obtained considering the first intersections of these orbits with $I$.  In particular, if $\pi^t(\al)>1$, the vertical segment connecting $u^t_\al$ to the conical singularity where its orbit ends has length $\Im(\xi^t_\al)$, where the complex number $\xi^t_\al$ is defined at the beginning of \S \ref{backgroundss1}. Similarly, the singularities of $T^{-1}$ can be obtained as first intersection with $I$ of the vertical orbits in positive time of the conical singularities of $X$, thus, if $\pi^b(\be)>1$, the vertical segment connecting $u^b_\be$ to the conical singularity where its orbit starts has length $-\Im(\xi^b_\be)$.

\begin{lem}\label{lem2s2ss2}
Consider data $(\pi,\la,\tau)$ as above and let $X$ and $T$ be the corresponding translation surface and
IET respectively. The following holds.
\begin{enumerate}
\item
For any triple $(\be,\al,n)$ reduced for $T$ and such that
$
|T^n(u^{b}_{\be})-u^{t}_{\al}|<
\min_{\chi\in\cA}\la_\chi
$
there exists a saddle connection
$\ga$ whose period $v$ is reduced and satisfies
$
|\Re(v)|=|T^{n}(u^{b}_{\be})-u^{t}_{\al}|
$
with $n=\sharp(I\cap\ga)$.
\item
For any saddle connection $\ga$ whose period $v$ is reduced and satisfies
$
|\Im(v)|>
\max_{\al,\be}
\Im(\xi^t_\al-\xi^b_\be)
$
there exists a triple $(\be,\al,n)$ such
that
$
|\Re(v)|=|T^n(u^b_\be)-u^t_\al|
$
and
$
n=\sharp(I\cap\ga)
$.
\end{enumerate}
\end{lem}

\begin{dimo}
We first prove part (1). Consider a triple $(\be,\al,n)$ as in the statement. Let $V_{out}$ be the outgoing vertical half-line whose first intersection with $I$ is $u^b_\be$ and call $p_{out}$ its starting point. Similarly, let $V_{in}$ be
the ingoing vertical half-line whose first
intersection with $I$ is $u^t_\al$ and call $p_{in}$ its ending point. Recall that the points $\xi_\be^b$ and $\xi_\al^t$ defined in \S \ref{backgroundss1} are identified to $p_{out}$ and $p_{in}$ respectively. Our aim is to define a saddle connection $\ga$ starting in $p_{out}$ and ending in $p_{in}$ and satisfying the properties in the statement. Assume without loss of generality that $T^n(u^b_\be)<u^t_\al$, the other case being the same, then set
$
\delta:=
u^{t}_{\al}-T^{n}(u^{b}_{\be})
$
and observe that by assumption we have $\de<\la_\chi$ for any $\chi$ in $\cA$. In particular, $T^n(u^b_\be)$ belongs to the subinterval $I^t_A$, where $A$ is the letter with $\pi^t(A)=\pi^t(\al)-1$. Equation (\ref{eq2s4ss3controloncords}) implies that $\xi^b_\be$ belongs to the boundary of the rectangle $R^b_\be$. Consider the horizontal segment in $R^b_\be$ with length $\de$ and left endpoint $\xi^b_\be$ and let $H$ be the image of the same segment under the isometric embedding of $R^b_\be$ in $X$. In particular $H$ has length $\de$ and its left endpoint is $p_{out}$. Since $(\be,\al,n)$ is reduced for $T$, then $T^{-k}$ is defined on the subinterval
$I(\be,\al,n)$ of $I$ for $0\leq k\leq n$. Recalling that $T$ is the first return of $\phi_t$ to $I$, we have instants $0<t(0)<\dots<t(n)$ such that
$
\phi_{t(k)}(H)=
T^{k-n}\big(I(\be,\al,n)\big)
$
for $0\leq k\leq n$. Let $s$ in $(0,\delta)$ be
the length coordinate on $H$ and call
$p_{s}$ the point of $H$ at distance $s$ from $p_{out}$. Since $(\be,\al,n)$ is
reduced, then $\phi_{t}(p_{s})$ is not a conical singularity for $s$ in $(0,\delta)$ and $t$ in
$(0,t(n)+h_A)$. For $s$ and $t$ in the same intervals, set
$$
\rho(s,t):=\phi_{t}(p_{s}),
$$
which is an isometric immersion from
$
(0,\delta)\times(0,t(n)+h_A)
$
to $X$ which avoids singularities. Observe that $\rho(0,0)=p_{out}$. If
$
\Im(\xi^t_\al)\leq h_A
$
then we set $a:=t(n)+\Im(\xi^t_\al)$ and
$
\ga(t):=\rho(t\delta,ta)
$,
which defines a saddle connection,
since $\rho(\delta,a)=p_{in}$. Otherwise, if $\Im(\xi^t_\al)>h_A$ then we must have $\pi^b(\al)=d$ and $0<\sum_\chi\tau_\chi<h_B$ where $B$ is the letter such that $\pi^t(B)=d$. Moreover the upper side of the rectangle $R^b_A$ is identified to the lower side of $R^t_B$, where the right endpoints coincide, therefore $\rho(s,t)$ extends to
$
(0,\delta)\times(0,t(n)+h_A+h_B)
$.
We have
$
\Im(\xi^t_\al)=
\sum_\xi\tau_\xi+h_A
$,
thus setting again $a:=t(n)+\Im(\xi^t_\al)$ and
$
\ga(t):=\rho(t\delta,ta)
$
we get a saddle connection, because $\rho(\delta,a)=p_{in}$. In both cases above, the period
$v=\int_{\ga}w$ obviously satisfies the required conditions.

Now we prove part (2). Let $\ga$ be a saddle connection whose period $v$ is as in the statement. Assume that $\Re(v)>0$, the other case being the same. Since $v$
is reduced, we have an isometric immersion
$
\rho:(0,\Re(v))\times(0,\Im(v))\to X
$
which avoids the conical singularities and such that
$\rho(tv)=\ga(t)$. Since
$
\Im(v)>\max_{\al,\be}
\Im(\xi^t_\al-\xi^b_\be)
$
then $\ga$ crosses $I$, thus there are instants
$
0<t(1)<\dots<t(n)<\Im(v)
$
such that
$
\rho
\big((0,\Re(v))\times\{t(k)\}\big)
\subset I
$
for $0\leq k\leq n$. The left and right vertical sides of the rectangle
$
(0,\Re(v))\times(0,\Im(v))
$
correspond to segment of vertical half-trajectories in $X$ respectively starting at $p_{out}$ and ending in $p_{in}$. Therefore there exist letters $\be$ and $\al$ such that for $0\leq k\leq n$ we have
$$
\rho\big((0,\Re(v))\times\{t(k)\}\big)=
\big(T^{k}(u^b_\be),T^{k-n}(u^t_\al)\big).
$$
The required identity follows considering $k=n$ in the equation above.
\end{dimo}

\begin{rem}In general, the triple $(\be,\al,n)$ in part (2) of Lemma \ref{lem2s2ss2} fails to be reduced just because
$
\big(T^n(u^b_\be),u^t_\al\big)
$
may contain a singularity for $T$ whose positive vertical orbit reaches a conical singularity in time greater than $\Im(v)-t(n)$ and similarly
$
\big(u^b_\be,T^{-n}(u^t_\al)\big)
$
may contain a singularity for $T^{-1}$ whose negative vertical orbit reaches a conical singularity in time greater than $t(0)$.
\end{rem}

\subsection{Proof of Vorobets identity}

Consider data $(\pi,\la,\tau)$, let $X$ be the underlying translation surface and $T$ be the IET corresponding to $(\pi,\la)$. Here we prove Proposition \ref{prop:vorobetsidentity}, that is
$$
\cE(T)=a(X)=\frac{s^{2}(X)}{2}.
$$
We first show separately the second equality.

\begin{lem}\label{lem1s2ss3}
For any translation surface $X$ we have
$$
a(X)=\frac{s(X)^2}{2}.
$$
\end{lem}

\begin{dimo}
Let $v_{n}$ be a sequence of periods in $\hol(X)$ such that $
\area(v_{n})\to a(X)\cdot\area(X)
$.
Since any $v_{n}$ belongs to a fixed vertical cone, then
$
\angolo(v_{n},\sqrt{-1})\to 0
$,
that is the angle with the vertical direction goes to zero. If $v$ is a period of $X$, corresponding to the saddle connection $\ga$, then for any $t$ we denote $|v|_t$ the length of $\ga$ with respect to the flat metric of $\cF_t\cdot X$. Consider the sequence of
positive instants $\tau_{n}$ such that
$
\area(v_{n})=|v_{n}|^{2}_{\tau_{n}}/2
$.
Since $\angolo(v_{n},\sqrt{-1})\to 0$ then $\tau_{n}\to\infty$, therefore we have
$$
\frac{s^{2}(X)}{2}=
\liminf_{t\to\infty}
\frac{\sys^{2}(\cF_t\cdot X)}{2\cdot\area(X)}\leq
\liminf_{n\to\infty}
\frac{\sys^{2}(\cF_{\tau_n}\cdot X)}{2\cdot\area(X)}\leq
\lim_{n\to\infty}
\frac{|v_{n}|^{2}_{\tau_n}}{2\cdot\area(X)}=
\lim_{n\to\infty}
\frac{\area(v_{n})}{\area(X)}=
a(X).
$$
On the other hand, consider positive instants $t_n$ with
$t_n\to+\infty$ and such that
$
\sys(\cF_{t_n}\cdot X)\to s(X)
$.
Consider a sequence of periods $v_n$ in $\hol(X)$ such that
$
|v_n|_{t_n}=
\sys(\cF_{t_n}\cdot X)
$.
Finally, for any such $v_n$, let $\tau_n$ be the unique instant such that
$
\area(v_n)=|v_n|^{2}_{\tau_n}/2
$.
Since $t_n\to+\infty$ then
$
\angolo(v_n,\sqrt{-1})\to 0
$,
otherwise
$\area(v_n)\to\infty$. This implies that the sequence of periods
$v_n$ eventually belongs to any fixed vertical cone. Hence we have
$$
\frac{s^2(X)}{2}=
\lim_{n\to\infty}
\frac{\sys^2(\cF_{t_n}\cdot X)}{2}=
\lim_{n\to\infty}
\frac{|v_n|^2_{t_n}}{2}\geq
\liminf_{n\to\infty}
\frac{|v_n|^2_{\tau_n}}{2}=
\liminf_{n\to\infty}\area(v_n)\geq
a(X).
$$
\end{dimo}

Now we prove the second equality, that is $\cE(T)=a(X)$. The proof is much simpler for those $X$ whose vertical flow is not uniquely ergodic. For such
$X$, Masur's criterion (\cite{masur})
implies that the projection of $\cF_t\cdot X$ to $\cM_g$ diverges, hence $\cF_t \cdot X$ diverges in its stratum too, therefore we have $s(X)=0$, which is equivalent to $a(X)=0$ according to Lemma \ref{lem1s2ss3}. On the other hand, since $T$ is not uniquely ergodic, Boshernitzan's criterion (see \S \ref{s1ss2sss1}) implies $\cE(T)=0$. Thus, we can prove the identity assuming that $X$ has uniquely ergodic vertical flow. Recall that if the vertical flow $\phi^t$ of $X$ is uniquely ergodic, then for any $p\in X$ we have
\begin{equation}\label{eq1s2ss3}
\lim_{t\to+\infty}
\frac
{\sharp\{s\in(0,t);\phi^s(p)\in I\}}
{t}
=
\frac{|I|}{\area(X)}.
\end{equation}
Equation (\ref{eq1s2ss3}) applies in particular to any trajectory $\phi^t(p)$ emanating from a conical singularity $p$. In this case, let $\ga$ be a saddle connection starting at $p$ whose period $v=\int_\ga w$ is reduced. Let $R_v$ be the rectangle in $\HH$ whose diagonal is $v$ and $\rho:R_v\to X$ be the corresponding isometric immersion. Then the paths
$
t\mapsto\rho(0,t\Im(v))
$
and
$
t\mapsto\rho(t\Re(v),t\Im(v))
$
with $t\in [0,1]$ have the same number of intersections with the transversal segment $I$. Therefore, for any $\epsilon>0$, if $\ga$ is a saddle connection whose period
$v$ is reduced and $|\Im(v)|$ is big enough we have
$$
(1-\epsilon)\cdot\frac{|I|}{\area(X)}
<
\frac{\sharp(\ga\cap I)}{|\Im(v)|}<
(1+\epsilon)\cdot\frac{|I|}{\area(X)}.
$$

Let $(v_{k})_{k\in\NN}$ be a sequence of reduced periods such that
$
\area(v_{k})\to a(X)\cdot\area(X)
$.
For any $k$ let $\ga_{k}$ be the saddle connection
corresponding to $v_{k}$ and $(\be(k),\al(k),n(k))$ be the sequence of triples
given by part (2) of Lemma \ref{lem2s2ss2}. For $k$ big enough $|\Im(v_{k})|$ is arbitrarily
large, hence we have
$$
\frac{\area(v_{k})}{\area(X)}=
\frac{|\Re(v_{k})|\cdot|\Im(v_{k})|}{\area(X)}\geq
\frac{\sharp(\ga_{k}\cap I)\cdot|\Re(v_{k})|}
{|I|\cdot(1+\epsilon)}=
\frac{n(k)\cdot|T^{n(k)}u^{b}_{\be(k)}-u^{t}_{\al(k)}|}
{|I|\cdot(1+\epsilon)}\geq
\frac{l(T)}{(1+\epsilon)^{2}}.
$$

On the other hand, let $(\be,\al,n)$ be a reduced triple
such that
$
n\cdot
|T^{n}u^{b}_{\be}-u^{t}_{\al}|
\leq
(1+\epsilon)l(T)
$
and $n$ is big enough. Let $\ga$ be the saddle connection given by part (1) of Lemma \ref{lem2s2ss2}. For the
associated period $v$ we have
$$
\frac{n\cdot|T^{n}u^{b}_{\be}-u^{t}_{\al}|}{|I|}=
\frac{\sharp(I\cap\ga)\cdot|\Re(v)|}{|I|}\geq
\frac{(1-\epsilon)|\Im(v)|\cdot|\Re(v)|}{\area(X)}\geq
(1-\epsilon)^{2}a(X).
$$
Proposition \ref{prop:vorobetsidentity} is proved.

\section{Closure of spectra and density of values of periodic elements}\label{s3}

Fix a connected component $\cC$ of some stratum of translation surfaces let $\cI$ be a closed invariant locus for the action of $\slduer$ on $\cC$. Let $\cR$ be a Rauzy class representing the surfaces in $\cC$, as it is explained in Proposition \ref{prop1backgroundss1}. In this section we establish a formula (Theorem \ref{renormalizedformula}) which allows to compute values of Lagrange spectra $\cL(\cI)$ via the \emph{Rauzy-Veech induction}, that is in terms of the map $\widehat{Q}$ defined in \S \ref{backgroundss3}. Then we use the formula to give a combinatorial proof of Theorem \ref{thm:closure}.

\medskip

We stress that the map $\widehat{Q}$ is defined on data $(\pi,\la,\tau)$ such that the pair $(\pi,\la)$ satisfies condition \eqref{eq1backgroundss2}. Therefore it can be iterated infinitely many times on data $(\pi,\la,\tau)$ such that the corresponding translation surface $X(\pi,\la,\tau)$ does not have vertical saddle connections. The latter is a full-measure subset of parameter space, and moreover it contains all data $(\pi,\la,\tau)$ of bounded type. The iteration of $\widehat{Q}$ produces a sequence
$
(\pi^{(r)},\la^{(r)},\tau^{(r)})_{r\in\NN}
$
of combinatorial-length-suspension data. For any fixed \emph{Rauzy time} $r$, the $r$-th iterate $\widehat{Q}^r$ is a piecewise linear map. Linear pieces are in bijection with finite \emph{Rauzy paths}. For any such path $\ga$, starting and ending at combinatorial data $\pi$ and $\pi'$ respectively, there are non-empty open sub-cones $\De_\ga\subset\RR_+^\cA$ and $\Th_\ga\subset\Th_{\pi'}$ such that $\widehat{Q}$ is a linear homeomorphism from $\{\pi\}\times\De_\ga\times\Th_\pi$ to $\{\pi'\}\times\RR_+^\cA\times\Th_\ga$.

\subsubsection{Notation for Rauzy paths}\label{subsubsection:notationrauzypaths}

Let $T$ be an IET corresponding to data $(\pi,\la)$. If $T$ admits $r$ steps of the map $Q$, we denote $\ga(0,r)$ the unique finite Rauzy path of combinatorial length $r$ such that
$
T\in\{\pi\}\times\De_{\ga(0,r)}
$.
Similarly, fix data $(\pi,\la,\tau)$ admitting infinitely many steps of the map $\widehat{Q}$, both in the future and in the past. For any pair of Rauzy times $r_-<0<r_+$ we call $\ga(r_-,r_+)$ the concatenation of the Rauzy path $\ga(r_-,0)$ ending in $\pi$, with the path $\ga(0,r_+)$ starting at $\pi$, such that
$$
(\pi,\la,\tau)\in
\{\pi\}\times\De_{\ga(0,r_+)}\times\Th_{\ga(r_-,0)}.
$$
Taking the limits of the corresponding cones, the paths $\ga(r_-,r_+)$ can be defined for $r_-=-\infty$ and/or $r_+=+\infty$. The paths $\ga(0,r)$ and $\ga(r_-,r_+)$ obviously depend on $(\pi,\la)$ and $(\pi,\la,\tau)$ respectively. Therefore, when referring to these paths, we will always specify the corresponding data inducing them. In \S \ref{sec:finitetime}, \ref{sec:subshifts} and \ref{subsection:closureanddensity} we will often denote by the letter $\ga$ both finite, half-infinite and bi-infinite Rauzy paths, taking care to avoid ambiguities when necessary.

\subsection{Renormalized formula on strata}\label{s3ss1(formulastrata)}

In order to state and prove the renormalized formula (Theorem \ref{renormalizedformula} below), we first introduce the notion of  \emph{diagonals} in Rauzy-Veech induction.

\subsubsection{Diagonals}\label{s3ss1sss1}

Fix data $(\pi,\la,\tau)$, let $X$ be the underlying translation surface and $T$ be the IET corresponding to $(\pi,\la)$. The interval $I=(0,\sum_{\chi\in\cA}\la_\chi)$ where $T$ acts is
naturally embedded in $X$. For $\chi\in\cA$ denote $e_\chi$ the correspondent element of the standard
basis of $\RR^\cA$. For a pair of letters $\be$ and $\al$ such that
$\pi^b(\be)>1$ and $\pi^t(\al)>1$ consider the elements of $\ZZ^\cA$
defined by
$$
w^b_{\be,\pi}:=
\sum_{\pi^b(\chi)<\pi^b(\be)}e_\chi
\textrm{ and }
w^t_{\al,\pi}:=
\sum_{\pi^t(\chi)<\pi^t(\al)}e_\chi,
$$
then set
$
w_{\be,\al,\pi}:=
w^b_{\be,\pi}-w^t_{\al,\pi}
$.
For any letter $\al$ in $\cA$, the values
$
\xi^t_\al:=
\langle\la,w^t_{\al,\pi}\rangle+
i\langle\tau,w^t_{\al,\pi}\rangle
$
and
$
\xi^b_\al:=
\langle\la,w^b_{\al,\pi}\rangle+
i\langle\tau,w^b_{\al,\pi}\rangle
$
are the coordinates of the conical singularities of $X$ in the zippered rectangle construction induced by $(\pi,\la,\tau)$ (see \S \ref{backgroundss1}).

\begin{defi}\label{def:diagonal}
We say that a period $v$ in $\hol(X)$ is a \emph{diagonal} with respect to $(\pi,\la,\tau)$ if there exist letters $\be,\al$ with $\pi^{b}(\be)>1$ and
$\pi^{t}(\al)>1$ such that
\begin{equation}\label{eq1s3ss1diagonals}
v=
\langle\la,w_{\pi,\be,\al}\rangle+i
\langle\tau,w_{\pi,\be,\al}\rangle.
\end{equation}
\end{defi}

Notice that not any pair of letters $\be,\al$ corresponds to a diagonal, indeed the vector $v$ defined by Equation \eqref{eq1s3ss1diagonals} may not be a period of $X$. Nevertheless we have the following Lemma, whose
simple proof is left to the reader.

\begin{lem}\label{lem1s3ss1}
If the letters $\be$ and $\al$ are such that
$$
|\langle\la,w_{\pi,\be,\al}\rangle|\cdot
|\langle\tau,w_{\pi,\be,\al}\rangle|=
\min_{\be',\al'\in\cA}
|\langle\la,w_{\pi,\be',\al'}\rangle|\cdot
|\langle\tau,w_{\pi,\be',\al'}\rangle|
$$
then Equation \eqref{eq1s3ss1diagonals} defines a reduced period $v$ of $X$, thus a diagonal with respect to $(\pi,\la,\tau)$.
\end{lem}

According to Lemma \ref{lem1s3ss1}, from the sequence of data
$
(\pi^{(r)},\la^{(r)},\tau^{(r)})
$
given by iteration of the Rauzy-Veech algorithm we can produce a sequence of reduced periods of $X$ as diagonals. Conversely, we have the following  criterion to establish whether or not a period is a diagonal with respect some renormalized data
$
(\pi^{(r)},\la^{(r)},\tau^{(r)})
$.

\begin{prop}\label{prop1s3ss1}
Consider a triple $(\be,\al,n)$ reduced for $T$  and such that
$
|T^n(u^{b}_{\be})-u^{t}_{\al}|<
\min_{\chi\in\cA}\la_\chi
$.
Let $v$ be the period associated to $(\be,\al,n)$ by  Lemma \ref{lem2s2ss2}. Then there exists $r$ in
$\NN$ such that $v$ is a diagonal with
respect to $(\pi^{(r)},\la^{(r)},\tau^{(r)})$.
\end{prop}

The Proposition follows immediately as a Corollary of Lemma \ref{lem2s2ss2} and the following Lemma, which  was first shown by  Rauzy and is proved in this form in \cite{luca1} (see Lemma 3.3 in \cite{luca1}).

\begin{lem}\label{lem2s3ss1}
If the triple $(\be,\al,n)$ is reduced for $T=(\pi,\la)$ then there exists $r$
in $\NN$ such that
$$
|T^{n}u^{b}_{\be}-u^{t}_{\al}|=
|u^{(r),t}_{\al}-u^{(r),b}_{\be}|,
$$
where $u^{(r),t}_{\al}$ and $u^{(r),b}_{\be}$ denote the singularities
of $T^{(r)}=(\pi^{(r)},\la^{(r)})$.
\end{lem}

\subsubsection{The renormalized formula}\label{s3ss1sss2}

Consider data $(\pi,\la,\tau)$ and let $X$ be the underlying translation surface. Recall that
$
\area(X)=
\langle
\la,\Omega_\pi\tau
\rangle
$
for a proper antisymmetric matrix $\Omega_\pi$, depending on $\pi$ (see for example \cite{yoccoz}, page 20).

\begin{defi}\label{def:w}
We define a continuous function on $(\pi,\la,\tau)$ setting
$$
w(\pi,\la,\tau):=
\frac{1}{\langle \la,\Omega_{\pi}\tau\rangle}
\min_{\be,\al\in\cA}
|\langle\la,w_{\pi,\be,\al}\rangle|\cdot
|\langle\tau,w_{\pi,\be,\al}\rangle|.
$$
\end{defi}

According to Lemma \ref{lem1s3ss1}, for a translation surface with unitary area, $w(\pi,\la,\tau)$ equals to $\area(v)$ of some period $v$ of $X$ which is diagonal with respect to $(\pi,\la,\tau)$ in the sense of Definition \ref{def:diagonal}.

\begin{thm}\label{renormalizedformula}
Let $X$ be the translation surface without vertical saddle connection and corresponding to the data $(\pi,\la,\tau)$. We have
$$
a(X)=
\liminf_{r\to\infty}
w(\pi^{(r)},\la^{(r)},\tau^{(r)}).
$$
\end{thm}

\begin{dimo}
Let $T$ be the IET corresponding to $(\pi,\la)$. For any $r$ let $\be(r)$ and $\al(r)$ be the letters satisfying the assumptions of Lemma \ref{lem1s3ss1} with respect to the data $(\pi^{(r)},\la^{(r)},\tau^{(r)})$ and let
$$
v_{r}=
\langle \la^{(r)},w_{\be(r),\al(r),\pi^{(r)}}\rangle
+i
\langle \tau^{(r)},w_{\be(r),\al(r),\pi^{(r)}}\rangle
$$
be the corresponding period. We have obviously
$$
\liminf_{r\to\infty}
w(\pi^{(r)},\la^{(r)},\tau^{(r)})=
\liminf_{r\to\infty}
\frac{\area(v_{r})}{\area(X)}\geq
a(X).
$$
In order to prove the Proposition it just remains to establish the reverse inequality. According to Proposition \ref{prop:vorobetsidentity} and Proposition \ref{prop1s2ss1},
it is enough to prove the required inequality replacing $a(X)$ by $l(T)$. Fix $\epsilon>0$ and two positive integers $N>0$ and $R>0$. Consider a triple
$(\be,\al,n)$ reduced for $T$ and such that $n\geq N$ and
$
n\cdot|T^{n}u^{b}_{\be}-u^{t}_{\al}|/|I|<l(T)+\epsilon.
$
Let $\ga$ and $v$ be respectively the saddle connection and the period given by
part (1) of Lemma \ref{lem2s2ss2}. According to Proposition \ref{prop1s3ss1}, $v$ is a diagonal with respect to
$(\pi^{(r)},\la^{(r)},\tau^{(r)})$ for some proper $r$. Moreover, modulo taking a bigger $N$, we can suppose $r>R$. Let $\pi=\pi^{(r)}$ and
let $\be$ and $\al$ be the letters such that
$
v=
\langle \la^{(r)},w_{\be,\al,\pi}\rangle
+i
\langle \tau^{(r)},w_{\be,\al,\pi}\rangle
$.
We have
$$
l(T)+\epsilon>
\frac{n\cdot|T^{n}u^{b}_{\be}-u^{t}_{\al}|}{|I|}=
\frac
{\sharp(\ga\cap I)\cdot|\langle\la^{(r)},w_{\be,\al,\pi}\rangle|}
{|I|}.
$$
Let us first assume that $l(T)>0$, which implies that the vertical flow of $X$ is uniquely
ergodic. Observe that we have
$
|\Im(v)|=
|\langle\tau^{(r)},w_{\be,\al,\pi}\rangle|
$.
Therefore, according to Equation \eqref{eq1s2ss3}, unique ergodicity implies
$$
\frac
{\sharp(\ga\cap I)
\cdot
|\langle\la^{(r)},w_{\be,\al,\pi}\rangle|}
{|I|}
\geq
(1-\epsilon)
\frac
{|\langle\tau^{(r)},w_{\be,\al,\pi}\rangle|
\cdot
|\langle\la^{(r)},w_{\be,\al,\pi}\rangle|}
{\area(X)}
\geq
(1-\epsilon)
w(\pi,\la^{(r)},\tau^{(r)}).
$$
We get
$
a(X)+\epsilon=
l(T)+\epsilon\geq
(1-\epsilon)\cdot
w(\pi^{(r)},\la^{(r)},\tau^{(r)})
$,
hence the required inequality follows because $\epsilon$ is arbitrarily small and $r$ is arbitrarily big.

Now assume that $l(T)=0$. In this case, for
$
L:=\max_{\al,\be}
\Im(\xi^t_\al-\xi^b_\be)
$
we have $|\Im(v)|\leq L\cdot\sharp(\ga\cap I)$, hence
$$
\frac
{\sharp(\ga\cap I)
\cdot
|\langle\la^{(r)},w_{\be,\al,\pi}\rangle|}
{|I|}
\geq
\frac
{|\langle\tau^{(r)},w_{\be,\al,\pi}\rangle|
\cdot
|\langle\la^{(r)},w_{\be,\al,\pi}\rangle|}
{\area(X)\cdot L}
\geq
L^{-1}\cdot w(\pi^{(r)},\la^{(r)},\tau^{r}).
$$
It follows that
$
\epsilon=l(T)+\epsilon
\geq L^{-1}
w(\pi^{(r)},\la^{(r)},\tau^{(r)})
$,
hence we get
$
\liminf_{r\to+\infty}
w(\pi^{(r)},\la^{(r)},\tau^{(r)})
=0
$,
because $\epsilon$ is arbitrarily small and $r$ is arbitrarily big.
\end{dimo}

\subsection{Control of the excursions via the cocycle}\label{sec:finitetime}

Let $\ga$ be any Rauzy path, either finite, bi-infinite, or half-infinite. We say that a finite sub-path $\nu$ of $\ga$ is \emph{minimal positive} if it is positive and it does not have proper positive sub-paths. Denote $\cM(\ga)$ the set of minimal positive sub-paths $\nu$ of $\ga$, then define
\beq\label{eq1sec:finitetime}
N(\ga):=
\sup\{\|B_\nu\|;\nu\in\cM(\ga)\},
\eeq
where the norm of the real $d\times d$ matrix $A$ is
$
\|A\|:=
\max_{\al,\be\in\cA}|[A]_{\al,\be}|
$.

Fix data combinatorial length and suspension data $(\pi,\la,\tau)$ inducing a bi-infinite path $\ga$. Let $X=X(\pi,\la,\tau)$ and $T=(\pi,\la)$ be respectively the corresponding translation surface and IET. For $r\geq 0$ the positive half-infinite path $\ga(r,+\infty)$ just depends on $T$. Define
$$
N_\infty(T):=
\limsup_{r\to\infty}
N\big(\ga(r,+\infty)\big).
$$
The following Theorem \ref{thm1s4} establishes a relation between $N_\infty(T)$ and the function $X\mapsto a(X)$, which gives rise to Lagrange Spectra. Its proof is given in \S \ref{s4ss2sss1} and is based on Proposition \ref{prop1s4ss1} and Proposition \ref{prop2s4ss1} in the following section.

\begin{thm}\label{thm1s4}
Fix combinatorial-length-suspension data $(\pi,\la,\tau)$, let $X$ be the underlying translation surface and $T$ be the IET corresponding to $(\pi,\la)$. We have
$$
\frac{1}{\sqrt[d-1]{d!}}\cdot
N_\infty(T)^{\frac{1}{2(d-1)(2d-3)}}
\leq
\frac{1}{a(X)}
\leq
\frac{d}{2}\cdot
N_\infty(T)^4.
$$
\end{thm}

From Theorem \ref{thm1s4} we can also deduce a characterization of IETs of bounded type (defined in \S~\ref{s1ss2sss1}) in terms of a suitable accelaration of the Rauzy-Veech induction. Let $T$ be the IET given by $(\pi,\la)$ which induces an half-infinite Rauzy path $\ga(0,+\infty)$. Denote $\nu(T)$ the minimal positive Rauzy path of the form $\ga(0,r)$ with $r\in\NN_\ast$. We can define an acceleration $Q_+$ of the Rauzy-Veech induction, that we call \emph{positive acceleration}, by setting $Q_{+}(T) = Q^{\nu(T)}(T)$. Finally set $P(T):=B_{\nu(T)}$ and let $P_n (T) = P(Q_{+}^n(T))$ be the sequence of positive matrices produced by the accelerated induction.

\begin{cor}\label{cor:boundedtypeIETs}
$T$ is an IET of bounded type if and only if it admits infinitely many iterations of the positive acceleration $Q_+$ and the corresponding matrices $P_n(T)$ have norms uniformly bounded for $n \in \NN$.
\end{cor}

After this paper, Kim and Marmi announced a proof of a similar characterization of IETs of bounded type in \cite{kimmarmi}.

\subsubsection{Positive matrices and flat geometry}\label{s4ss1}

The following two Propositions contain  two geometrical estimates of the distortion of data $(\pi,\la,\tau)$ in terms of a finite and bounded number of steps of the Rauzy-Veech induction. Their proof is quite technical and is given in \S \ref{s7}. Fix combinatorial-length-suspension data $(\pi,\la,\tau)$ and let $T$ be the IET corresponding to $(\pi,\la)$.

\begin{defi}\label{def:distorsion} We define the \emph{distortion} of $T$ by
$$
\Delta(T):=
\max_{\be,\al\in\cA}
\frac{\la_{\al}}{\la_{\be}}.
$$
\end{defi}
Recall from \cite{mmy} that a finite Rauzy path $\ga$ is said \emph{complete} if any letter $\al$ in $\cA$ \emph{wins} in $\ga$ (we refer to \S \ref{backgroundss2} for terminology). Then we say that a path $\ga$ is \emph{strongly complete} if it is a concatenation of $d$ complete paths. Let $\ga(-\infty,0)$ be the half-infinite path induced by the data $(\pi,\la,\tau)$. Consider the minimal positive integer $m$ such that the ending part $\ga(-m,0)$ of $\ga(-\infty,0)$ is strongly complete. Then set
$$
m(\pi,\la,\tau):=
\min_{-m\leq r\leq 0}
w(\pi^{(r)},\la^{(r)},\tau^{(r)}).
$$

\begin{prop}\label{prop1s4ss1}
Consider combinatorial-length-suspension data $(\pi,\la,\tau)$ and let
$T$ be the IET corresponding to $(\pi,\la)$. We have
$$
\Delta(T)\leq
\frac{d!}{m(\pi,\la,\tau)^{d-1}}.
$$
\end{prop}

\begin{prop}\label{prop2s4ss1}
Let $\ga_1,\ga_2,\ga_3,\ga_4$ be positive paths that can be
concatenated and let $r,s$ be instants
such that  that $\ga(-s,0)=\ga_1\ast\ga_2$ and $\ga(0,r)=\ga_3\ast\ga_4$. Then
we have
$$
w(\pi,\la,\tau)\geq
\frac{2}{d}\cdot
\frac
{1}
{\|B_{\ga_1}\|_\infty\cdot
\|B_{\ga_2}\|_\infty\cdot
\|B_{\ga_3}\|_\infty\cdot
\|B_{\ga_4}\|_\infty}.
$$
\end{prop}

The proof of Proposition \ref{prop1s4ss1} and Proposition \ref{prop2s4ss1} is given in \S \ref{s7}. In the rest of the section we apply them to prove Theorem \ref{thm1s4}.

\subsubsection{Auxiliary results}

If $\ga$ is a finite Rauzy path of length $r$ and $T$ is an IET defined by data $(\pi,\la)$ in $\De_\ga$, denote $Q_\ga$ the branch of the Rauzy map induced by $\ga$, that is
$
(\pi^{(r)},\la^{(r)})=
Q_\ga(\pi,\la)
$.
If $T^{(r)}$ is the IET corresponding to $(\pi^{(r)},\la^{(r)})$, we also write $T^{(r)}=Q_\ga(T)$.

\medskip

Consider combinatorial-length data $(\pi,\la)$ and let $T$ be the corresponding IET. Let $\ga(0,+\infty)$ be the half-infinite path induced by $(\pi,\la)$ and recall that $P(T):=B_{\nu(T)}$ where $\nu(T)$ is the smallest  $r\in\NN_\ast$ such that the Rauzy path $\ga(0,r)$ is  positive.

\begin{lem}\label{lem1s4ss1}
Let $T$ be any IET. We have
$$
\Delta(T)\leq \|P(T)\|.
$$
\end{lem}

\begin{dimo}
Set $(\pi',\la')=Q_{\nu(T)}(\pi,\la)$ and $P:=P(T)$. We
have $^tP\la'=\la$. Since $P$ is positive with integer entries then for any $\al$ in $\cA$ we have
$
\sum_{\chi\in\cA}\la'_\chi
\leq
\la_{\al}
\leq
\|P\|_{\infty}\cdot\sum_{\chi\in\cA}\la'_\chi
$.
It follows that for any $\be$ and $\al$ we have
$
\la_{\al}\leq \|P\|_{\infty}\cdot\la_{\be}
$,
thus the Lemma is proved.
\end{dimo}

\begin{lem}\label{lem2s4ss1}
Let $\ga$ be a non-complete path and $(\pi,\la)$ be combinatorial-length data in $\De_\ga$. If $T$ is the IET corresponding to the data $(\pi,\la)$ then we have
$$
\|B_\ga\|_\infty\leq
\De\big(T\big)\cdot\De\big(Q_\ga(T)\big).
$$
\end{lem}

\begin{dimo}
Observe that we have $^tB_\ga\la'=\la$. Let $\al$ be a letter which newer wins in $\ga$. The $\al$-column of $B_\ga$ equals to the vector $e_\al$ of the standard basis of $\RR^\cA$, that is $[B_\ga]_{\al,\al}=1$ and $[B_\ga]_{\chi,\al}=0$ for $\chi\not=\al$, thus $\la_\al=\la'_\al$. On the other hand, if $(\be',\be)$ is the pair of letters such that
$
\|B_\ga\|_\infty=[B_\ga]_{\be,\be'}
$,
we have
$$
\la_{\be}=
\sum_{\chi\in\cA}[B_\ga]_{\chi,\be}\la'_\chi>
[B_\ga]_{\be',\be}\la'_{\be'}=
\|B_\ga\|_\infty\la'_{\be'}.
$$
Therefore we get
$$
\la_\al=
\la'_\al<
\De\big(Q_\ga(T)\big)\cdot\la'_{\be'}<
\De\big(Q_\ga(T)\big)
\frac{\la_\be}{\|B_\ga\|_{\infty}}<
\De\big(Q_\ga(T)\big)
\frac
{\De\big(T\big)\cdot\la_\al}
{\|B_\ga\|_{\infty}}.
$$
\end{dimo}

The following Corollary is motivated by the Lemma at page 833 in \cite{mmy}, where it is proved that the concatenation of $2d-3$ maximal non-complete paths is
positive.

\begin{cor}\label{cor1s4ss1}
Consider a finite Rauzy path $\ga$ that is the concatenation
$
\ga=\ga_1\ast\dots\ast\ga_{2d-3}
$
of $2d-3$ maximal non-complete paths $\ga_1,\dots,\ga_{2d-3}$. For any $k$ with $1\leq k\leq 2d-3$ denote
$Q_k:=Q_{\ga_1\ast\dots\ast\ga_k}$. If $T$ is an IET in $\De_\ga$, then we have
$$
\|P(T)\|_\infty\leq
\De(T)\cdot\De^2(Q_1(T))\cdot\dots
\cdot\De^2(Q_{2d-4}(T))\cdot\De(Q_{2d-3}(T)).
$$
\end{cor}

\begin{dimo}
Observe that
$
\|B_\ga\|_\infty\leq
\|B_{\ga_1}\|_\infty
\cdot\dots\cdot
\|B_{\ga_{2d-3}}\|_\infty
$.
Observe also that $T$ belongs to $\De_{\ga_1}$ and moreover, for any $1<k\leq 2d-3$ the IET $Q_{k-1}(T)$ belongs to $\De_{\ga_k}$. Then the Corollary follows applying Lemma \ref{lem2s4ss1} to each factor in the product.
\end{dimo}

\subsubsection{Proof of Theorem \ref{thm1s4}}\label{s4ss2sss1}

Fix combinatorial-length-suspension data $(\pi,\la,\tau)$, let $\ga=\ga(-\infty,+\infty)$ be the corresponding bi-infinite Rauzy path and $X$ be the underlying translation surface. Let $T$ be the IET corresponding to $(\pi,\la)$.

We first prove the second inequality. Fix $M>N_\infty(T)$. For any $r$ big enough, we can find integers $r_-<r<r+$ and positive paths $\ga_{1},\ga_{2},\ga_{3},\ga_{4}$ with $\|B_{\ga_i}\|<M$ for $i=1,\dots,4$ such that
$
\ga_{1}\ast\ga_{2}=\ga(r_-,r-1)
$
and
$
\ga_{3}\ast\ga_{4}=\ga(r,r_+)
$.
For these big values of $r$, Proposition \ref{prop2s4ss1} implies
$
M^{-4}\cdot 2/d\leq
w(\pi^{(r)},\la^{(r)},\tau^{(r)})
$,
therefore we have
$$
\frac{2}{d\cdot M^4}\leq
\liminf_{r\to\infty}
w(\pi^{(r)},\la^{(r)},\tau^{(r)}).
$$
According to the formula in Proposition \ref{renormalizedformula} we have $M^{-4}\cdot 2/d\leq a(X)$, and since this last estimate holds for any $M>N_\infty(T)$, then the second inequality in Theorem \ref{thm1s4} follows.

Now we prove the first inequality. Recall the formula in Proposition \ref{renormalizedformula} and observe that we have
$$
a(X)=
\liminf_{r\to+\infty}
w(\pi^{(r)},\la^{(r)},\tau^{(r)})=
\liminf_{r\to+\infty}
m(\pi^{(r)},\la^{(r)},\tau^{(r)}),
$$
therefore Proposition \ref{prop1s4ss1} implies
$
\limsup_{r\to\infty}
\De(T^{(r)})\leq
d!\cdot a(X)^{-(d-1)}
$.
Observe that with our notation we have
$
N_\infty(T)=
\limsup_{r\to\infty}
\|P(T^{(r)})\|_\infty
$,
thus
$
N_\infty(T)\leq
\limsup_{r\to+\infty}
\Delta(T^{(r)})^{2(2d-3)}
$,
according to Corollary \ref{cor1s4ss1}. Combining the last two results we get
$$
N_\infty(T)\leq
(d!)^{2(2d-3)}
\cdot
a(X)^{-2(d-1)(2d-3)}.
$$
The first inequality in Theorem \ref{thm1s4} follows. $\square$

\subsection{Subshifts in the Rauzy class}\label{sec:subshifts}

In this subsection we consider bi-infinite Rauzy paths $\ga$ such that the quantity $N(\ga)$ introduced in Equation \eqref{eq1sec:finitetime} is finite and we show that these paths can be encoded by a symbolic sequence in a subshift
of finite type. This coding will be convenient to present the proof of Theorem \ref{thm:closure}.

\subsubsection{Construction of the sub-shift}

For any finite Rauzy path
$\ga$ consider the corresponding matrix $B_{\ga}$ in the \emph{Kontsevich-Zorich
cocycle}, that is the matrix defined in \S \ref{backgroundss2sss2}. Recall that $\ga$ is said to be \emph{positive} if all the entries of the matrix $B_\ga$ are positive. Consider the countable alphabet $\Ga$ of finite and positive Rauzy paths in $\cR$. Let $\Ga^\ZZ$ be the shift space of sequences in the elements of $\Ga$. Let $\Ga^\ZZ_\cR\subset\Ga^\ZZ$ be the sub-shift space consisting of all
sequences $(\ga_k)_{k\in\ZZ}$ in $\Gamma^\ZZ$ such that for all $k\in\ZZ$ the final permutation of $\ga_k$ is the same than the initial
permutation of $\ga_{k+1}$, so that the infinite concatenation
$
\cdots\ast\ga_{k-1}\ast\ga_ k\ast\ga_{k+1}\ast\cdots
$
is a bi-infinite path on $\cR$. For any $M>0$, let $\Ga(M)\subset\Ga$ be the subset of the alphabet of those $\ga\in\Ga$ such that $\|B_\ga\|\leq M$. Clearly $\Ga(M)$ is a finite alphabet. Let $\Ga^\ZZ(M)$ be the sub-shift space consisting of all bi-infinite sequences in the letters of $\Ga(M)$ and let
$$
\Ga^\ZZ_\cR(M):=
\Ga^\ZZ_\cR \cap \Ga^\ZZ(M) .
$$
Then the shift
$
\sigma:\Ga^\ZZ_\cR(M)\to \Ga^\ZZ_\cR(M)
$
is a topological Markov chain, the transition matrix $A$ being defined setting $A_{\ga,\ga'}:=1$ if and only if the
concatenation $\ga\ast\ga'$ is a path on $\cR$, and $A_{\ga,\ga'}:=0$ otherwise, where $\ga,\ga'$ are any pair of letters in $\Ga(M)$. Remark that while the whole shift $\sigma$ is not of finite type, for any $M>0$ the sub-shift
$
\sigma:\Ga^\ZZ_\cR(M)\to \Ga^\ZZ_\cR(M)
$
is of finite type. Recall the function $\ga\mapsto N(\ga)$ defined by Equation \eqref{eq1sec:finitetime}. According to the following Lemma, whose proof is obvious, any bi-infinite path $\ga$ with $N(\ga)<M$ is decomposed into \emph{packets}, each packet $\ga_k$ being the entry of a sequence in $\Ga^\ZZ_\cR(M)$.

\begin{lem}\label{lem:packing}
Let $\ga$ be a bi-infinite path such that $N(\ga)<M$ for some $M>0$. Then there exists a sequence
$
(\ga_k)_{k\in\ZZ}\in\Ga^\ZZ_\cR(M)
$
such that
$$
\ga=\dots\ga_{k-1}\ast\ga_k\ast\ga_{k+1}\dots.
$$
\end{lem}

\subsubsection{Hilbert metric}\label{subsubsection:hilbertmetric}

The Hilbert pseudo-metric on $\RR^2_+$ is
$$
\distanza_{\RR^2_+}(x,y):=
\log\max_{1\leq i,j\leq 2}
\frac{x_ix_j}{x_jy_i},
$$
which is invariant under linear isomorphisms of $\RR^2_+$, that is the group of invertible diagonal matrices. More generally, if $C\subset \RR^\cA\setminus\{0\}$ is an open convex cone whose closure does not contain any
one-dimensional subspace of $\RR^\cA$, one defines a Hilbert pseudo-metric on $C$ as follows. If $x$ and $y$
are collinear, then $\distanza_C(x,y)=0$. Otherwise, $C$ intersects the subspace generated by $x$ and $y$ in a
cone isomorphic to $\RR^2_+$. We let
$
\distanza_{C}(x,y):=
\distanza_{\RR^2_+}(\psi x,\psi y)
$,
where $\psi$ is any such isomorphism. When $C=\RR^\cA$ the formula above gives
$$
\distanza_{\RR^\cA}(x,y)=
\log\max_{\be,\al}
\frac{x_\al y_\be}{x_\be y_\al}.
$$

%Consider an open convex cone $C$ in $\RR^\cA$ which does not contain a line in its closure. Let $x$ and $y$ be two distinct points of $C$. Consider the points $x'$ and $y'$ at which the line passing through $x$ and $y$ intersects the boundary of $C$, where the order of the points is $x',x,y,y'$. The \emph{Hilbert's pseudo metric} on $C$ is
%$$
%\distanza_C(x,y):=
%\log\frac{\|x-y'\|\cdot\|y-x'\|}{\|x-x'\|\cdot\|y-y'\|}.
%$$
%Observe that $x'$ and $y'$ cannot be simultaneously equal to $\infty$, because $C$ does not contains any line. On the other hand, if $C$ contains the half line starting from $x$ and pointing to $y$, then we have $y'=\infty$. In this case we have $\|y-x'\|=\|x-x'\|=\infty$ and the Hilbert's pseudo metric becomes
%$\distanza_C(x,y):=\log\frac{\|x-y'\|}{\|y-y'\|}$.
%The case $x'=\infty$ similar.

It is easy to check that for any two points $x$ and $y$ in $C$ and for any pair of positive real numbers $a$ and $b$ we have
$
\distanza_C(x,y)=
\distanza_C(a\cdot x,b\cdot y)
$,
therefore the Hilbert pseudo-metric induces a well defined metric on the space of
rays $\{tx:t\in\RR_+\}$ contained in $C$, called \emph{Hilbert's metric}.

Let $A$ be linear map such that $A(C)\subset C$. According to a classical result due to G. Birkhoff (see \cite{birkhoff}), if $A(C)$ is bounded in $C$ then $A$ is a contraction of the Hilbert pseudo-metric. More precisely, denoting by $D$ the diameter of $A(C)$, for any $x,y$ in $C$ we have
$$
\distanza_C(A\cdot x,A\cdot y)\leq
\delta\cdot
\distanza_C(x,y),
\textrm{ where }
\delta=\frac{\sqrt{e^D}-1}{\sqrt{e^D}+1}.
$$
The last estimate implies the following Lemma, whose proof is left to the reader.

\begin{lem}\label{lem:contraction}
For any $M>0$ there exists a positive constant $\delta=\delta(M)<1$ such that, for any Rauzy path $\ga:\pi\to\pi'$ with $\|B_\ga\|\leq M$ the following holds
\begin{enumerate}
\item
If $^tB_\ga\RR^\cA_+$ has compact closure in $\RR^\cA_+$ then the map $\la \mapsto ^tB_\ga(\la)$ is a $\delta$-contraction.
\item
If $^tB^{-1}_\ga\Th_\pi$ has compact closure in $\Th_{\pi'}$ then the map $\tau\mapsto ^tB_\ga(\tau)$ is a $\delta$-contraction.
\end{enumerate}
\end{lem}

%\begin{dimo}
%According to Birkhoff's estimate, in both cases it is enough to estimate the diameter of $^tB_\ga\RR^\cA_+$ or $^tB^{-1}_\ga\Th_\pi$ respectively. Recall that for $C=\RR_+^\cA$ the Hilbert's pseudo metric takes the form
%$$
%\distanza_{\RR^\cA_+}(x,y)=
%\log\max_{\be,\al}
%\frac{x_\al y_\be}{x_\be y_\al}.
%$$
%Observe that $^tB_\ga\RR^\cA_+$ has compact closure in $\RR^\cA_+$ if and only if $^{t}B_\ga(e_\al)$ has all its entries positive for any element $e_\al$ of the standard basis of $\RR^\cA$. According to the above expression for the Hilbert's pseudo metric, we have $\diametro(^tB_\ga\RR^\cA_+)\leq\log (\|B_\ga\|^2)$.  On the other hand, consider two extremal points $u,v$ of $\Th_\pi$ and let $\psi$ be a linear isomorphism between the cone $\spanne(^tB^{-1}_\ga(u),^tB^{-1}_\ga(v))\cap \Th_{\pi'}$ and $\RR^2_+$. Let $u':=\psi\big(^tB^{-1}_\ga(u)\big)$ and $v':=\psi\big(^tB^{-1}_\ga(v)\big)$. Since there are finitely many matrices $B$ with integer coefficients, the same argument applies as for $\RR^\cA_+$. The Lemma is proved.
%\end{dimo}

\subsubsection{A continuous function for the sub-shift}

Fix $M>0$ and let $\ga$ be a bi-infinite Rauzy path with $N(\ga)\leq M$, then according to Lemma \ref{lem:packing} decompose it as
$$
\ga=
\dots\ga_{k-1}\ast\ga_k\ast\ga_{k+1}\ast\dots
\textrm{ with }
(\ga_k)_{k\in\ZZ}\in\Ga^\ZZ_\cR(M).
$$
For any $r$ in $\ZZ$ consider the half-infinite path $\ga(r,+\infty)$ defined in \S \ref{subsubsection:notationrauzypaths} and let $\pi^{(r)}$ be the element in $\cR$ where $\ga(r,+\infty)$ starts. Since any $\ga_k$ is positive then $^tB_{\ga_k}\RR^\cA_+$ has compact closure in $\RR^\cA_+$ with respect to the Hilbert's metric. Lemma 4.3 in \cite{agy} implies that the same property holds for the cone of suspension data, modulo considering the concatenation
$
\ga_k\ast\dots\ast\ga_{k+3d-5}
$
of at most $3d-4$ positive paths, where $d$ denotes the cardinality of $\cA$. Therefore, according to Lemma \ref{lem:contraction}, the cone
$
\bigcap^s_r
\big(^{t}B_{\ga(r,s)}\RR^\cA_+\big)
$
shrinks exponentially with respect to the Hilbert's metric as $s\to+\infty$ and similarly
$
\bigcap^{r-1}_{s}
\big(^{t}B^{-1}_{\ga(s,r-1)}\Th_{\pi^{(s)}}\big)
$
shrinks exponentially as $s\to-\infty$. Therefore we have two rays $\widehat{\la}^{(r)}$ in $\RR^\cA_+$ and $\widehat{\tau}^{(r)}$ in $\Th_{\pi^{(r)}}$ defined by
\begin{eqnarray*}
&&\widehat{\la}^{(r)}:=
\bigcap^\infty_{s=r}
\big(^{t}B_{\ga(r,s)}\RR^\cA_+\big)\\
&&\widehat{\tau}^{(r)}:=
\bigcap^{r-1}_{s=-\infty}
\big(^{t}B^{-1}_{\ga(s,r-1)}\Th_{\pi^{(s)}}\big).
\end{eqnarray*}

\begin{defi}\label{def:continuefunction}
Fix a bi-infinite path with $N(\ga)<+\infty$ and any $r\in\ZZ$. Define length and suspension data $\la^{(r)}$ and $\tau^{(r)}$ on the rays $\widehat{\la}^{(r)}$ and $\widehat{\tau}^{(r)}$ respectively via the normalizations
$\|\la^{(r)}\|=1$ and
$
\area(\pi^{(r)},\la^{(r)},\tau^{(r)})=1
$.
Then define
$$
a_r(\ga):=w(\pi^{(r)},\la^{(r)},\tau^{(r)}).
$$
\end{defi}

\begin{lem}[Approximation Lemma]\label{lem:approximation}  For any $\epsilon >0$
and $M>0$ there exists $m=m(\epsilon,M) \in \NN$ such that the following holds. Let $\ga', \ga''$  be two  bi-infinite paths on $\cR$ with $N(\ga')<M$ and $N(\ga'')<M$ and let $r^-< r<
r^+$ be any Rauzy times  such that we can write
\begin{eqnarray*}
&&\ga'(r^-,r)=\ga''(r^-, r) =
\ga_{1}\ast\dots\ast\ga_{m}, \\
&&\ga'(r,r^+)=\ga''(r,r^+)=
\ga_{m+1}\ast\dots\ast\ga_{2m},
\end{eqnarray*}
where  $\ga_1, \dots, \ga_{2m}$ are $2m$ positive paths in $\Ga(M)$.
Then
$$
\left|
\frac{a_r(\ga')}{a_r(\ga'')} -1
\right|
\leq \epsilon .
$$
\end{lem}

\begin{dimo}
For any $r$ let
$
(\pi,\la^{'(r)},\tau^{'(r)})
$
and
$
(\pi,\la^{''(r)},\tau^{''(r)})
$
be the unique data associated respectively to $\ga'$ and to $\ga''$ as in Definition \ref{def:continuefunction}. Proposition \ref{prop2s4ss1} implies
$$
a_r(\ga')\geq\frac{2}{d\cdot M^4}
\textrm{ and }
a_r(\ga'')\geq\frac{2}{d\cdot M^4}.
$$
The condition
$
\ga'(r,r^+)=\ga''(r,r^+)=
\ga_{m+1}\ast\dots\ast\ga_{2m}
$
implies that both $\la^{'(r)}$ and $\la^{''(r)}$ belong to the cone $\De_{\ga_{m+1}\ast\dots\ast\ga_{2m}}$, whose projective diameter is less than $\delta^m$, according to Lemma \ref{lem:contraction}, where $\delta$ is the constant in the Lemma. Similarly, the condition
$
\ga'(r^-,r)=\ga''(r^-, r)=
\ga_{1}\ast\dots\ast\ga_{m}
$
implies that both $\tau^{'(r)}$ and $\tau^{''(r)}$ belong to a subset of $\Th_{\pi^{(r)}}$ with diameter bounded by $\delta^m$. The required $m$ exists because the function
$
(\pi,\la,\tau)\mapsto
w(\pi,\la,\tau)
$
is continuous.
\end{dimo}

Consider two sequences $(\ga'_k)$ and $(\ga''_k)$ in $\Ga^\ZZ_\cR$ and the corresponding bi-infinite Rauzy paths
$
\ga'=
\cdots\ga'_k\ast\ga'_{k+1}\ast\cdots
$
and
$
\ga''=
\cdots\ga''_k\ast\ga''_{k+1}\ast\cdots
$
obtained by bi-infinite concatenation. Recall that any bi-infinite Rauzy path admits such decomposition.

\begin{cor}[Interpolation]\label{cor:interpolation}
For any $\epsilon >0$ and $M>0$, let $m=m(\epsilon, M)$ be the positive integer given by the
Approximation Lemma \ref{lem:approximation}. If there are integers $k_1<k_2$ and finite paths $\sigma_i$ in $\Ga(M)$ with $k_1-m \leq i \leq k_2 +m$ such that
$$
\ga'_i=\ga''_i=\sigma_i
\textrm{ for any }
k_1-m \leq i \leq k_2 +m
$$
then for the instants $r_1<r_2$ such that
$
\ga'(r_1,r_2)=
\ga''(r_1,r_2)=
\sigma_{k_1}\ast\cdots\ast\sigma_{k_2}
$
and for any $r_1 \leq  r< r_2$ we have
$$
\left|
\frac{a_r(\gamma')}{a_r(\gamma'')} -1
\right|\leq\epsilon.
$$
\end{cor}

\begin{dimo}
Let $r_1 \leq  r< r_2$. Consider $k$ with $k_1 \leq k
\leq k_2$ and the two instants $r_k^-$ and $r_k^+$ with
$
r_1\leq r_k^-\leq r<r_k^+\leq r_2
$
such that
$
\ga_k=\ga(r_k^-,r)\ast\ga(r+1,r_k^+)
$.
Thus we have
$$
\ga(r_1,r)=
\sigma_{k_1}\ast\cdots\ast\sigma_{k-1}\ast\ga(r_k^-,r)
\textrm{ and }
\ga(r+1,r_2)=
\ga(r+1,r_k^+)\ast\sigma_{k+1}\cdots\ast\sigma_{k_2}
$$
Observe that
$
\gatilde_{k-1}:=\sigma_{k-1}\ast\ga(r_k^-,r)
$
and
$
\gatilde_{k+1}:=\ga(r+1,r_k^+)\ast\sigma_{k+1}
$
are both positive paths such that
$
\|B_{\gatilde_{k-1}}\|\leq M^2
$
and
$
\|B_{\gatilde_{k+1}}\|\leq M^2
$.
According to the assumptions, we have two Rauzy times $r^-<r<r^+$ such that
\begin{eqnarray*}
&&\ga'(r^-,r)=\ga''(r^-,r)=
\sigma_{k_1-m}\ast\dots\ast\gatilde_{k-1},\\
&&\ga'(r,r^+)=\ga''(r,r^+)=
\gatilde_{k+1}\ast\dots\ast\sigma_{2m},
\end{eqnarray*}
where $m=m(\epsilon, M^2)$ and $\gatilde_{k-1},\gatilde_{k+1} $ and all $\sigma_i$ appearing above belong to $\Ga(M^2)$. The Corollary follows applying the Approximation Lemma \ref{lem:approximation}.
\end{dimo}

\subsubsection{Recurrence for the sub-shift}

In this paragraph we show that for sub-shifts of finite type there always exists instants such that the hypothesis required by Lemma \ref{lem:approximation} and Corollary \ref{cor:interpolation} are satisfied.

\begin{lem}[Cantor diagonal Lemma]\label{lem:Cantor}
Consider a sequence of sequences
$
(\ga_k^{(n)})_{k \in \ZZ}
$
in $\Ga^\ZZ (M)$, indexed by $n\in\NN$. There exists a sequence $\sigma=(\sigma_k)_{k\in\ZZ}$ in $\Ga^\ZZ (M)$ and a subsequence $(n_r)_{r\in\NN}$ such that for any $m\in\NN$, we have
$$
\ga_k^{(n_r)}=\sigma_{k}
\textrm{ for all }
-m \leq k \leq m
\textrm{ and all }
r \geq m.
$$
\end{lem}

In other words, the subsequence
$(\ga_k^{(n_r)})_{k \in \ZZ}$ provided by the
Lemma consists of sequences which eventually agree on arbitrarily large central
blocks with the limiting sequence $(\sigma_k)_{k \in \ZZ}$.
\medskip

\begin{dimo}
Recall that $\Ga(M)$ is a finite set. Since the sequence of central letters
$
(\ga_0^{(n)})_{n\in\NN}
$
has values in $\Ga(M)$, then there exists a letter
$\sigma_0\in \Ga(M)$ which occurs infinitely often, that is there exists an increasing subsequence $(n(0,j))_{j\in\NN}$ of natural numbers such that  $\ga_0^{(n(0,r))}=\sigma_0$ for any $r\in\NN$. Fix $m\geq 0$ and suppose that we have $2m+1$ letters
$
\sigma_{-m},\dots,\sigma_{0},\sigma_1,\dots,\sigma_m
$
in $\Ga(M)$ and a nested family of subsequences
$$
n(m,\cdot)\subset
n(m-1,\cdot)\subset
\dots\subset n(0,\cdot)
$$
such that for any $r$ with $0\leq r\leq m$ and for any $j$ we have
$$
\ga_i{(n(r,j))}=\sigma_i
\textrm{ for all }
-r\leq i \leq r,
$$
where we stress that each $n(r,\cdot)$ is a sequence in the index $j$, that is $n(r,\cdot)=(n(r,j))_{j\in\NN}$. The sequence of pair of letters
$
(\ga_{-m-1}^{(n(m,j))},\ga_{m+1}^{(n(m,j))})_{j\in \NN}
$
has values in $\Ga(M)\times\Ga(M)$, which is a finite set. Therefore there exits a pair of letters $(\sigma_{-m-1},\sigma_{m+1})$ and a subsequence $n(m+1,\cdot)\subset n(m,\cdot)$, where $n(m+1,\cdot)=(n(m+1,j))_{j\in\NN}$, such that
$$
(\ga_{-m-1}^{(n(m+1,j))},\ga_{m+1}^{(n(m+1,j))})=
(\sigma_{-m-1},\sigma_{m+1})
\textrm{ for all }
j\in\NN.
$$
Therefore the sequence $\sigma=(\sigma_j)_{j\in\ZZ}$ is defined inductively, moreover we have also an inductively defined nested sequence of subsequences $(n(m,\cdot))_{m\in\NN}$ with
$
\dots
n(m,\cdot)\subset n(m-1,\cdot)
\dots.
$
The required subsequence is diagonal subsequence, that is $m_m:=n(m,m)$.
\end{dimo}

\begin{cor}[Recurrence]\label{cor:recurrence}
Consider a sequence $(\ga_k)_{k\in\ZZ}$ in
$\Ga^\ZZ(M)$ and any subsequence $(k_{n})_{n\in\NN}$ with $k_n\to\infty$. For any $m\in\NN$ there exist $2m+1 $
letters
$
\sigma_{-m},\dots,\sigma_{0},\sigma_1,\dots,\sigma_m
$
in $\Ga(M)$ such for infinitely many $n\in\NN$ we have
$$
\ga_{k_n+i}=\sigma_i
\textrm{ for all }
-m\leq i \leq m.
$$
In other words, the finite word
$
\sigma_{-m}\cdots\sigma_m
$
appears infinitely often in $(\ga_k)_{k\in\ZZ}$ and moreover it appears centered
at infinitely many of the prescribed indexes $(k_{n})_{n\in\NN}$.
\end{cor}

\begin{dimo}
For any $n$ define a sequence $(\ga^{(n)}_k)_{k\in\ZZ}$ in the index $k$ setting $\ga^{(n)}_k:=\ga_{k_n+k}$. Let $\sigma=(\sigma_k)_{k\in\ZZ}$ be the sequence in $\Ga^\ZZ(M)$ given by Lemma \ref{lem:Cantor} and $(n_r)_{r\in\NN}$ be the subsequence defined in the same Lemma. Fix any non-negative integer $m$. Then, for any $r\geq m$ set $n=n_r$. The Corollary follows observing that for any $i$ with $-m\leq i\leq m$ we have
$$
\ga_{k_n+i}=
\ga^{(n)}_i=
\sigma_i.
$$
\end{dimo}

\subsection{Proof of closure and density of periodic elements}\label{subsection:closureanddensity}

Recall that we consider a fixed connected component $\cC$ of some stratum of translation surfaces and a closed invariant locus $\cI$ for the action of $\slduer$ on $\cC$. Then $\cR$ is a Rauzy class representing the surfaces in $\cC$, as it is explained in Proposition \ref{prop1backgroundss1}. For any combinatorial datum $\pi\in\cR$ let $\cI_\pi$ be the set of data
$
(\la,\tau)\in\RR_+^\cA\times\Th_\pi
$
such that the surface $X(\pi,\la,\tau)$ belongs to $\cI$. The structure of $\cI_\pi$ is described in Appendix \ref{loci}.

\medskip

If $\ga$ is a positive loop at $\pi\in\cR$, denote respectively $\la^\ga$ and $\tau^\ga$ the maximal and the minimal eigenvector of the matrix $^tB_\ga^{-1}$, normalized so that $\|\la^\ga\|=1$ and
$
\area(\pi,\la^\ga,\tau^\ga)=1
$.
Let $X(\pi,\la^\ga,\tau^\ga)$ be the translation surface corresponding to such data, which is a periodic point for $\cF_t$. More precisely, for
$
T:=-\log\big(\|^tB_\ga^{-1}\la^\ga\|\big)
$
we have
$$
\cF_T\cdot X(\pi,\la^\ga,\tau^\ga)=
X(\pi,\la^\ga,\tau^\ga).
$$
If $\cI$ is a closed and $\slduer$-invariant subset of the connected component corresponding to the Rauzy class $\cR$, then define
$$
\periodici(\cR,\cI):=
\left\{
\ga
\textrm{ positive loop at }
\pi\in\cR
\textrm{ with }
X(\pi,\la^\ga,\tau^\ga)\in\cI
\right\}.
$$
In order to simplify the notation, in the following we will refer to elements of
$
\periodici(\cR,\cI)
$
either as positive loops $\ga$, or as bi-infinite paths corresponding to the infinite concatenation
$
\dots\ga\ast\ga\ast\ga\dots
$.
Any $\ga$ in $\periodici(\cR,\cI)$ corresponds to some closed Teichm\"uller geodesic in $\cI$, that is an element of $PA(\cI)$. Not o any element of $PA(\cI)$ is represented by a positive loop in $\periodici(\cR,\cI)$. Nevertheless, for any $X\in PA(\cI)$ whose orbit has period $T$, there exists a positive integer $m$ and some $\ga\in\periodici(\cR,\cI)$ representing the orbit
$
t\mapsto \cF_{mt}\cdot X
$
for $0\leq t\leq T$. In particular, the set of values
$
\left\{a(X);X\in PA(\cI)\right\}
$
can be computed just considering the elements in $\periodici(\cR,\cI)$.

\subsubsection{Reduction to the sub-shift}

Recall that $\cL(\cI)\subset\cL(\cC)$. We introduce a monotone increasing function $\cN:\cL(\cC)\to\RR_+$, which is obviously defined by restriction also on $\cL(\cI)$, given explicitly by
$$
\cN(L):=(d!)^{2(2d-3)}L^{2(d-1)(2d-3)}.
$$
Consider data $(\pi,\la,\tau)$ inducing a bi-infinite path $\ga$ and let $X=X(\pi,\la,\tau)$ be the corresponding translation surface. For $r\in\ZZ$ consider the corresponding paths $\ga(r,+\infty)$ defined in \S~\ref{subsubsection:notationrauzypaths}. The first inequality in Theorem \ref{thm1s4} implies
\beq\label{eq1:subshiftreduction}
\limsup_{r\to\infty}
N\big(\ga(r,+\infty)\big)\leq
\cN\big(a(X)^{-1}\big).
\eeq
If $N(\ga)<\infty$ then we can consider the function $r\mapsto a_r(\ga)$ introduced in Definition \ref{def:continuefunction} and in view of Theorem \ref{renormalizedformula} we have $a(X)=a(\ga)$, where
\beq\label{eq2:subshiftreduction}
a(\ga):=
\liminf_{r\to+\infty}
a_r(\ga).
\eeq
In general, Condition \eqref{eq1:subshiftreduction} just holds for some positive half-infinite segment $\ga(r,+\infty)$ of $\ga$ for $r$ big enough. According to Lemma \ref{lem:continuefunction} below, for any value in $\cL(\cI)$ it is possible to find a bi-infinite path $\ga$ with $N(\ga)<+\infty$, so that Equation \eqref{eq2:subshiftreduction} applies.

\begin{lem}\label{lem:continuefunction}
Consider any $a>0$ such that $a^{-1}\in\cL(\cI)$. Then there exists $\pi\in\cR$ and data $(\la,\tau)\in\cI_\pi$ such that the bi-infinite path $\ga$ induced by $(\pi,\la,\tau)$ satisfies  $$
N(\ga)\leq \cN(a^{-1})
\textrm{ and }
a=a(\ga).
$$
\end{lem}

\begin{dimo}
Set $M:=\cN(a^{-1})$. Consider $\pi\in\cR$ and data $(\tilde{\la},\tilde{\tau})\in\cI_\pi$ such that the surface
$
\tilde{X}:=
X(\pi,\tilde{\la},\tilde{\tau})$ gives
$
a\big(\tilde{X}\big)=a
$.
Let $\tilde{\ga}$ be the bi-infinite path induced by
$
(\pi,\tilde{\la},\tilde{\tau})
$
and observe that according to Equation \eqref{eq1:subshiftreduction} we have
$
N(\tilde{\ga}(r_0,+\infty))<M
$
for some $r_0\in\ZZ$. Modulo applying the Rauzy-Veech induction to the data
$
(\pi,\tilde{\la},\tilde{\tau})
$
we can assume that $r_0<-m$ for some $m>r_1(M)$, where $r_1(M)$ is the integer appearing in the statement of Proposition \ref{prop:appendixclosinglemma}. We have a decomposition
$$
\tilde{\ga}(r_0,+\infty)=
\ga_0\ast\ga_1\ast\cdots\ast\ga_k\ast\cdots
$$
where any $\ga_k$ is a finite positive sub-path of $\tilde{\ga}$ with $N(\ga_k)<M$. Fix $l\in\NN$ and according to Corollary \ref{cor:recurrence} let $(k_n)_{n\in\NN}$ be a subsequence with $k_n\to+\infty$ and $\sigma_{-l},\dots,\sigma_{l}$ be positive paths in $\Ga(M)$ such that for any $n\in\NN$ we have
$$
\ga_{k_n+i}=\sigma_i
\textrm{ for }
|i|\leq l.
$$
Observe that $\ga_{k_{n+1}-1}=\sigma_{-1}$ ends in the element of $\cR$ where $\ga_{k_n}=\sigma_0$ starts, therefore we obtain a positive loop consider the concatenation
$$
\ga^{(n)}:=
\ga_{k_n}\ast\cdots\ast\ga_{k_{n+1}-1}.
$$
Modulo taking a larger $l$ we can assume that the concatenation
$
\sigma_{-l}\ast\cdots\ast\sigma_{l}
$
has length larger than $2m$, in terms of Rauzy elementary operations. Modulo discarding some elementary arrows we can assume that
$
\sigma_{-l}\ast\cdots\ast\sigma_{l}
$
has length exactly $2m$, thus let $\nu$ and $\de$ be the two paths with length $m$ such that we have
$$
\sigma_{-l}\ast\cdots\ast\sigma_{l}=\nu\ast\de.
$$
Proposition \ref{prop:appendixclosinglemma} implies that
$
\ga^{(n)}\in\periodici(\cR,\cI)
$
and abusing the notation we denote with the same name the corresponding bi-infinite path. We have $N(\ga^{(n)})<M$, since $\ga^{(n)}$ is concatenation of sub-paths of $\tilde{\ga}(r_0,+\infty)$. Moreover, modulo shifting the base point of the loop we have
$$
\ga^{(n)}(-m+1,0)=\nu
\textrm{ and }
\ga^{(n)}(1,m)=\de,
$$
and we can assume without loss of generality that the base point of the loop is $\pi$. Modulo increasing $l$ and repeating the construction above, we can assume that $m$ is big enough to apply Proposition \ref{prop:appendixshadowinglemma} (in fact the constant $r_1$ appearing in Propositions \ref{prop:appendixclosinglemma} is actually the same as the one appearing in Proposition \ref{prop:appendixshadowinglemma}). Consider the infinite concatenation
$$
\ga:=
\cdots\ast\ga^{(1)}\ast\ga^{(1)}\ast
\ga^{(2)}\ast\cdots\ast
\ga^{(n)}\ast\ga^{(n+1)}\ast\cdots,
$$
where the past half-infinite part of $\ga$ coincides with the periodic path with period $\ga^{(1)}$. We obviously have $N(\ga)<M$. Moreover Proposition \ref{prop:appendixshadowinglemma} implies that $\ga$ is induced by data $(\pi,\la,\tau)$ with $(\la,\tau)\in\cI_\pi$. Finally observe that there is some $r_2>0$ such that
$
\tilde{\ga}(r_2,+\infty)=\ga(r_2,+\infty)
$.
Therefore there exists an IET $T$ which is a section both for the vertical flow of $X(\pi,\tilde{\la},\tilde{\tau})$ and the vertical flow of $X(\pi,\la,\tau)$. According to Proposition \ref{prop:vorobetsidentity} we have
$$
a(\ga)=
a\big(X(\pi,\la,\tau)\big)=
\cE(T)=
a\big(X(\pi,\tilde{\la},\tilde{\tau})\big)=a.
$$
\end{dimo}

According to Lemma \ref{lem:packing}, any bi-infinite path $\ga$ with $N(\ga)<M$ is decomposed into \emph{packets}, each packet $\ga_k$ being the entry of a sequence in $\Ga^\ZZ_\cR(M)$. Lemma \ref{lem:packingcompatible} provides a more accurate packing, satisfying a compatibility condition with a prescribed sequence of Rauzy times, which is a technical detail needed in the proof of  Propositions \ref{inclusiondensity} and \ref{inclusionclosure} below.

\begin{lem}[Packing Lemma]\label{lem:packingcompatible}
Fix $M>0$ and let $\ga$ be a bi-infinite Rauzy path with $N(\ga)<M$. Consider an increasing sequence of Rauzy times $(r_n)$ such that $\ga(r_n, r_{n+1}-1)$ is positive for any $n$.  Then there exists a sequence
$
(\ga_k)\in\Ga_{\cR}^\ZZ (M^2)
$
and an increasing subsequence of indexes $(k_n)$ such that for any $n$ we have
$$
\ga(r_n,r_{n+1}-1)=
\ga_{k_n}\ast\ga_{k_n+1}\ast\cdots\ast\ga_{k_{n+1}-1}.
$$
\end{lem}

\begin{dimo}
Since $N(\ga)<M$ then for any $n$, any sub-path $\nu$ of $\ga(r_n,r_{n+1})$ with $\|B_\nu\|\geq M$ is positive. If $\|B_{\ga(r_0,r_1-1)}\|< M^2$,
then we set $\ga_{k_0}:=\ga(r_0,r_1-1)$
and
$k_1:=k_0+1$.
Otherwise, we choose Rauzy times
$$
r_0=:s_0<s_1<\dots<s_m:=r_1
$$
such that for any $j$ with $0\leq j< m$ we have
$M\leq\|B_{\ga(s_{j},s_{j+1}-1})\|<M^2$
and set $k_1:=k_0+m$ and $\ga_{k_0+j}:=\ga(s_j,s_{j+1}-1)$. The required sequences
$
(\ga_k)_{k\in\ZZ}\in
\Ga_{\cR}^\ZZ (M^2)
$
and $(k_n)_{n\in\NN}$ are defined replying the same construction iteratively on $n$.
\end{dimo}

\subsubsection{Density of Pseudo-Anosov elements}

In this section we prove the first part of Theorem \ref{thm:closure}.

\begin{prop}\label{inclusiondensity}
We have
$$
\cL(\cI)
\subset
\overline{\{ a^{-1}(X)\textrm{ ; }X \in PA(\cI) \}}.
$$
\end{prop}

\begin{dimo}
Consider $a>0$ such that $a^{-1}$ belongs to $\cL(\cI)$. Then set $M:=\cN(a^{-1})$ and fix $\epsilon>0$. The statement follows defining a periodic path
$
\gatilde\in\periodici(\cR,\cI)
$
such that
$$
(1-\epsilon)^2 a<a(\tilde{\ga})<(1+\epsilon)^2 a.
$$
According to Lemma \ref{lem:continuefunction} there exists $\pi\in\cR$ and data $(\la,\tau)\in\cI_\pi$ such that the bi-infinite path $\ga$ induced by $(\pi,\la,\tau)$ satisfies $N(\ga)<M$ and $a(\ga)=a$. Consider a sequence of instants $(r_{n})_{n\in\NN}$ such that
$$
r_n\to+\infty
\textrm{ and }
a_{r_{n}}(\ga)\to a
\textrm{ for }
n\to+\infty.
$$
Up to extracting a sub-sequence we can assume that $\ga(r_{n},r_{n+1})$ is positive for any $n$. Thus, by the Packing Lemma \ref{lem:packingcompatible} there exists a sequence
$(\ga_{k})_{k\in\ZZ}$ in $\Gamma^{\ZZ}_{\cR}(M^{2})$ and a subsequence $(k_{n})_{n\in\ZZ}$ with $k_n\to+\infty$ for $n\to+\infty$ such that for $n$ big enough we have
$$
\ga(r_n, r_{n+1}) =
\ga_{k_n}\ast
\ga_{k_n+1}\ast
\cdots \ast
\ga_{k_{n+1}-1}.
$$
Fix $l\in\NN$. Corollary \ref{cor:recurrence} implies that there exist $2l+1$ elements
$
\sigma_{-l},\dots,
\sigma_{0},\sigma_1,\dots,
\sigma_l
$
in $\Gamma (M^{2})$ such that, up to extracting a sub-sequence of $(k_n)_{n\in\NN}$, and hence a subsequence of the instants $(r_n)_{n\in\NN}$, we can assume that for
any $n\in\NN$ we have
$$
\ga_{k_n+i}=
\sigma_i
\textrm{ for all }
-l\leq i \leq l.
$$
Without any loss of generality, we can also assume that $k_{n+1}-k_{n}>2l$.

\smallskip

{\it Step 1: Definition of the periodic path
$
\tilde{\ga}\in\periodici(\cR,\cI)
$}. Let $N$ be such that for any $n\geq N$ we have
$$
(1-\epsilon)a<
a_{r_{n}}(\ga)<
(1+\epsilon)a.
$$
Let $(\gatilde_{k})_{k\in\ZZ}$ be the periodic sequence in
$\Gamma^{\ZZ}$ with period
$k_{N+1}-k_{N}$ defined setting
$$
\gatilde_{k_N+k}:=\ga_{k_N+k}
\textrm{ for }
0\leq k<k_{N+1}-k_N.
$$
Observe that $(\gatilde_{k})_{k\in\ZZ}$ belongs to
$\Gamma^{\ZZ}_{\cR}$, and not just to $\Gamma^{\ZZ}$, that is the concatenation is possible. For $k\not=0 \mod k_{N+1}-k_N$ the concatenation $\ga_{k-1}\ast\ga_{k}$ is possible because $(\ga_k)$ is obtained via Lemma \ref{lem:packingcompatible}. For $k=0 \mod k_{N+1}-k_N$ just observe that by the recurrence property $\ga_{k_{N+1}-1}=\sigma_{-1}$ and $\ga_0=\sigma_0$. Therefore let $\gatilde$ be the bi-infinite periodic path on $\cR $ obtained concatenating the
paths in the sequence $(\gatilde_k)_{k\in\ZZ}$, that is
$$
\gatilde:=
\cdots
\gatilde_{k-1}\ast
\gatilde_{k} \ast
\gatilde_{k+1} \cdots.
$$
In particular, the period $\gatilde(r_N,r_{N+1})$ of $\tilde{\ga}$ is a positive loop at some $\pi\in\cR$. Moreover, modulo taking a larger $l$, we can assume that the concatenation
$
\sigma_{-l}\ast\dots\ast\sigma_l
$
has length at least $2m$ in terms of elementary Rauzy operations, where $m$ is big enough to apply Proposition \ref{prop:appendixclosinglemma}. Arguing as in the proof of Lemma \ref{lem:continuefunction} we get that
$
\gatilde(r_N,r_{N+1})\in\periodici(\cR,\cI)
$.

\smallskip

{\it Step 2: Extra equalities between $(\tilde{\ga}_k)_{k\in\ZZ}$ and $(\ga_k)_{k\in\ZZ}$ from recurrence.} For any $i$ with $0\leq i\leq l$ we have
$
\ga_{k_N-i}=
\sigma_i=
\ga_{k_{N+1}-i}=
\gatilde_{k_{N+1}-i}=
\gatilde_{k_N-i}
$,
where the first two equalities follow from the recurrence property of the instants $k_n$, and the last two from the definition of the sequence $(\gatilde_k)_{k\in\NN}$ and its periodicity. The symmetric argument shows that for any $i$ with $0\leq i\leq l$ we also have
$
\ga_{k_{N+1}+i}=
\gatilde_{k_{N+1}+i}
$.
Resuming, the choice of instants $k_n$ implies that the sequences $(\ga_k)_{k\in\NN}$ and $(\gatilde_k)_{k\in\NN}$ coincide not only for $k_N\leq k\leq k_{N+1}$ but also for a larger interval of instants, that is
$$
\tilde{\ga}_{k_N+k}=\ga_{k_N+k}
\textrm{ for }
-l\leq k<k_{N+1}-k_N+l.
$$

\smallskip

{\it Step 3: Areas are well  approximated.} Since the path $\gatilde$ is periodic with period $k_{N+1}-k_{N}$, there exists an instant $r_{min}$ with $r_N \leq r_{min}< r_{N+1}$ such that
$$
a_{r_{min}}(\gatilde)=
\min_{r_N\leq r<r_{N+1}}a_{r}(\gatilde)=
\liminf_{r\in\NN} a_r(\gatilde)=
a(\gatilde).
$$
Modulo taking a bigger $l$, we can assume that $l>m(\epsilon, M^2)$, where $m(\epsilon, M^2)$ is the positive integer appearing in Corollary \ref{cor:interpolation}. The extra equalities between $(\tilde{\ga}_k)_{k\in\ZZ}$ and
$(\ga_k)_{k\in\ZZ}$ proved in Step 2 imply that Corollary \ref{cor:interpolation} applies to the finite segments $\ga(r_N,r_{N+1})$ and $\gatilde(r_N,r_{N+1})$ and since
$
r_N \leq r_{min}< r_{N+1}
$
then we have
$$
a_{r_{min}}(\gatilde)\leq
a_{r_N}(\gatilde)\leq
(1+\epsilon)a_{r_N}(\ga)\leq
(1+\epsilon)^2a,
$$
where the first inequality follows by definition of the instant $r_{min}$, the second by the Interpolation Corollary \ref{cor:interpolation} and the last by the choice of the instant $r_N$. Similarly we have
$$
a_{r_{min}}(\gatilde)\geq
(1-\epsilon)a_{r_m}(\ga)\geq
(1-\epsilon)^2a,
$$
where the first equality follows again from Corollary \ref{cor:interpolation} and the second holds because $a=\liminf_r a_r(\ga)$. The Proposition follows observing that
$
a_{r_{min}}(\gatilde)=a(\gatilde)
$
and therefore
$$
(1-\epsilon)^2 <
\frac{a}{a(\gatilde)}<
(1+\epsilon)^2.
$$
\end{dimo}

\subsubsection{Closure of the spectrum}

In this section we prove the Proposition below, which gives the opposite inclusion to Proposition \ref{inclusiondensity} and allows to conclude the proof of Theorem \ref{thm:closure}.

\begin{prop}\label{inclusionclosure}
We have
$$
\overline{\{ a^{-1}(X)\textrm{ ; }X\in PA(\cI)\}}
\subset \cL(\cI).
$$
\end{prop}

\begin{dimo}
Let $a^{-1}$ be an accumulation point of the values $a^{-1}(X)$ for $X\in PA(\cI)$. Thus there exists a sequence
$
\gatilde^{(n)}\in\periodici(\cR,\cI)
$
such that
$$
a=\lim_{n\to\infty}a(\gatilde^{(n)}).
$$
Let $r_n$ be the period of $\gatilde^{(n)}$ and assume without loss of generality that any finite loop $\gatilde^{(n)}(0, r_n )$ is based at the same $\pi\in\cR$. We can assume also that $\gatilde^{(n)}(0, r_n )$ is a positive loop at $\pi$. Indeed it is a complete loop, otherwise $\gatilde^{(n)}$ is not admissible. Then positivity follows modulo concatenating at most $2d-3$ copies of it, according to \cite{mmy} (the Lemma at page 833).

We show that $a^{-1} \in \cL(\cI)$ by explicitly constructing a bi-infinite path $\ga$ on $\cR$, induced by data $(\pi,\la,\tau)\in\cI_\pi$ and such that $a(\ga)=a$. Observe that since $a>0$ then there exists some $\delta>0$ such that
$
a(\tilde{\gamma}^{(n)})>\delta
$
for any $n$. Therefore according to Equation \ref{eq1:subshiftreduction} there exists some $M>0$ such that
$$
N(\gatilde^{(n)})<M
\textrm{ for any }
n.
$$

\smallskip

{\it Step 1: Symbolic coding.}
For any $n$, according to the Packing Lemma \ref{lem:packingcompatible} applied to $\gatilde^{(n)}$, there exists a positive integer $p_n$ and finite paths
$
\gatilde^{(n)}_0,\dots,\gatilde^{(n)}_{p_n-1}
$
in $\Ga(M^2)$ such that
$$
\gatilde^{(n)}(0,r_n)=
\gatilde^{(n)}_0\ast
\gatilde^{(n)}_1 \ast
\cdots\gatilde^{(n)}_{p_n-1}.
$$
In other words, the periodic bi-infinite path $\gatilde^{(n)}$ is coded by the periodic sequence
$(\gatilde^{(n)}_k )_{k\in \ZZ}$ of period $p_n$ in the sub-shift space
$\Ga^\ZZ_\cR(M^2)$ defined by
$$
\gatilde^{(n)}_k :=
\gatilde^{(n)}_{k \mod p_n}
\textrm{ for any }
k\in\ZZ.
$$
We stress that $r_n$ is the period of $\gatilde^{(n)}$ in terms of Rauzy elementary operations, whereas $p_n$ is the period of its coding
$
(\gatilde^{(n)}_k )_{k\in \ZZ}
$
in terms if the sub-shift time. We assume without loss of generality that $p_n\geq n$ for any $n\in\NN$.

\smallskip

{\it Step 2: Definition of the path $\ga$.}
Apply Lemma \ref{lem:Cantor} to the
sequence of sequences $(\gatilde^{(n)}_k )_{k\in \ZZ}$ in $\Ga^\ZZ_\cR(M^2)$, indexed by $n\in\NN$. The Lemma provides a sequence $(\sigma_k)_{k\in\ZZ}$ in $\Ga^\ZZ_\cR(M^2)$ and a sequence of integers $n(l)$, indexed by $l\in\NN$ and with $n(l)\to+\infty$, such that for any $l$ in $\NN$  we have
$$
\gatilde_k^{(n(j))}=
\sigma_k
\textrm{ for all }
-l \leq k \leq l
\textrm{ and all }
j\geq l.
$$
In order to simplify the notation, we assume that $n(l)=l$, that is we assume that $(\gatilde^{(n)}_k )_{k\in \ZZ}$ satisfies itself the conclusion of Lemma \ref{lem:Cantor}. Thus, resuming, for any $l\in\NN$ we have
$$
\gatilde_k^{(n)}=
\sigma_{k}
\textrm{ for all } -l \leq k \leq l
\textrm{ and all }n \geq l.
$$
Fix $l\in\NN$. We construct a sequence $(\ga_k)_{k\in\ZZ}$ in $\Ga^\ZZ_\cR(M^2)$ by concatenating the
periods of the sequences $(\gatilde^{(n)}_k )_{k\in \ZZ}$ with $l\geq n$ as follows.
\begin{enumerate}
\item
Set $P_l:=p_l$ and for $n>l$ set
$P_n:=\sum_{i=l}^n p_i$, so that $P_n-P_{n-1}=p_n$.
\item
For $k<P_l$ define the entries $\ga_k$ repeating infinitely the period
$
\gatilde^{(l)}_0,\dots,\gatilde^{(l)}_{p_l-1}
$
of $(\gatilde^{(l)}_k )_{k\in\ZZ}$, that is
$$
\ga_k:=
\gatilde^{(l)}_{k\mod P_l}.
$$
\item
For $n>l$ and any $k$ with $P_{n-1}\leq k < P_{n}$ define the entries $\ga_k$ to be equal to the period
$
\gatilde^{(n)}_0,\dots,\gatilde^{(n)}_{p_n-1}
$
of $(\gatilde^{(n)}_k )_{k\in\ZZ}$, that is
$$
\ga_k:=
\gatilde^{(n+1)}_{k-P_{n-1}}.
$$
\end{enumerate}
For any $n\geq l$ we have
$
\ga_{P_n-1}=
\gatilde^{(n)}_{p_n-1} =
\gatilde^{(n)}_{-1}=
\sigma_{-1}
$
and
$
\ga_{P_n}=
\gatilde^{(n+1)}_0=
\sigma_0
$,
since the sequence $(\gatilde^{(n)}_k )_{k\in \ZZ}$ is periodic with period $p_n$. Hence the concatenation
$
\ga_{P_n-1}\ast\ga_{P_n}
=
\sigma_{-1}\ast\sigma_0
$
is admissible, because $(\sigma_k)_{k \in \ZZ}$ belongs to $\Ga^\ZZ_\cR$. Therefore we can define an bi-infinite path $\ga$ on $\cR$ setting
$$
\ga:=
\dots\ga_{k-1}\ast\ga_{k}\ast\ga_{k+1}\ast\dots.
$$
Finally recall that any $\gatilde^{(n)}$ belongs to $\periodici(\cR,\cI)$ and that $N(\gatilde_k^{(n)})$ is uniformly bounded. Assume without loss of generality that $l$ is such that the concatenation
$
\sigma_{-l}\ast\cdots\ast\sigma_{l}
$
has even length $2m$ in terms of Rauzy elementary operations. Then let $\nu$ and $\de$ be the two paths with length $m$ such that we have
$$
\sigma_{-l}\ast\cdots\ast\sigma_{l}=\nu\ast\de.
$$
Modulo taking a bigger $l$ we can assume that $m$ is big enough to apply Proposition \ref{prop:appendixshadowinglemma}. Therefore the path $\ga$ defined above is induced by data $(\pi,\la,\tau)$ with $(\la,\tau)\in\cI_\pi$.

\smallskip

{\it Step 3: Extra equalities between $\ga$ and $(\gatilde^{(n)})$.} Here we consider any $l\in\NN$ and we show that for all $n>l$ we have
$$
\ga_j=\gatilde^{(n)}_{j-P_{n-1}}
\textrm{ for }
P_{n-1}-l\leq j<P_{n}+l.
$$
When $P_{n-1}\leq j<P_{n}$ the required equality follows trivially from the definition of the sequence $(\ga_j)$ in Step 2. For $P_{n}\leq j< P_{n}+l$, setting $j=P_{n}+k$ we have
$$
\ga_j=
\ga_{P_{n}+k}=
\gatilde^{(n+1)}_{k}=
\sigma_k=
\gatilde^{(n)}_k=
\gatilde^{(n)}_{k+p_{n}}=
\gatilde^{(n-1)}_{j-P_{n-1}},
$$
where the first and last equality follow from the change of index, the second form the definition of $\ga_{P_{n}+k}$, the third and the fourth from Lemma \ref{lem:Cantor} (and since $|k|<l$) and finally the fifth holds because $(\gatilde^{(n)}_k)$ is periodic with period $p_{n}$. The proof for $P_{n-1}-l\leq j<P_{n}$ is similar.

\smallskip

{\it Step 4: Approximation properties of $\gamma$.}
For any $n$ let $R_n$ be the length of the path $\gamma(0, P_n)$ in terms of elementary Rauzy operations, so that $R_{n}-R_{n-1}=r_n$ is the period of $\gatilde^{(n)}$. Recall that we have
$
\ga(R_{n-1},R_{n})=
\gatilde^{(n)}(0,R_{n}-R_{n-1})
$,
by definition of $\ga$, and that $N(\gatilde^{(n)})<M$ and $N(\ga)<M^2$. Fix $\epsilon>0$. According to Step 3 the assumptions in Corollary \ref{cor:interpolation} are satisfied, thus there exists $l=l(\epsilon, M^2)$ such that for any $n>l$ and any $R_n \leq  r< R_{n+1}$ we have
$$
\left|
\frac{a_r(\ga)}{a_{r-R_n}(\gatilde^{(n)})} -1
\right|
\leq\epsilon.
$$

\smallskip

{\it Step 5: Computation of $a(\gamma)$.} Each $\gatilde^{(n)}$ has period $R_{n}-R_{n-1}=r_n$ in terms of elementary Rauzy operations, thus there exists an integer $r_{min}(n)$ with $R_{n-1}\leq r_{min}(n)<R_{n}$ such that
$$
a_{r_{min}(n)-R_{n-1}}(\gatilde^{(n)})=
\min_{R_{n-1}\leq  r<R_{n}}
a_{r-R_{n-1}}(\gatilde^{(n)}) =
a(\gatilde^{(n)}).
$$
Since $a(\gatilde^{(n)})\to a$ for $n\to\infty$, then the approximation established at step 4 implies
$$
(1-\epsilon) a\leq
\liminf_n a_{r_{min}(n)}(\ga) \leq
(1+\epsilon) a.
$$
and hence
$
a(\ga)=
\liminf_{r\to\infty}
a_r(\ga)\leq (1+\epsilon)a
$.
On the other hand, for any $n$ and any $r$ with
$
R_{n-1}\leq r<R_{n}
$
the approximation at Step 4 implies
$$
a_r(\ga)
\geq
(1-\epsilon)
a_{r-R_{n-1}}(\gatilde^{(n)})
\geq
(1-\epsilon)
a_{r_{min}(n)-R_{n-1}}(\gatilde^{(n)})=
(1-\epsilon)
a(\gatilde^{(n)}).
$$
Since $a(\gatilde^{(n)})\to a$ for $n\to\infty$, then we get
$
\liminf_{r\to\infty}a_r(\ga)
\geq (1-\epsilon)a
$.
Moreover the two estimations for $\liminf_ra_r(\ga)$ hold for any $\epsilon$, then the Proposition is proved.

\end{dimo}

\subsubsection{End of the proof}

\begin{proofof}{Theorem}{thm:closure}
The inclusion in Proposition \ref{inclusiondensity} proves the density of values of periodic elements, that is part (2) of Theorem \ref{thm:closure}. Proposition
\ref{inclusionclosure} gives the opposite inclusion, therefore
we get
$$
\cL(\cI)=
\overline{\{a^{-1}(X)\textrm{ ; }X\in PA(\cI)\}},
$$
which implies in particular that $\cL(\cI)$ is a closed subset of $\RR$.
\end{proofof}

\section{Hall's ray for Square-tiled surfaces}\label{section5}

\subsection{Background on square-tiled surfaces}\label{appendixss2}

We recall some basic facts on square-tiled surfaces. We closely follow \cite{hubertlelievre}. A \emph{translation covering} is a map $f:X_1\to X_2$ of translation surfaces that
\begin{enumerate}
\item
is topologically a ramified covering;
\item
maps conical singularities of $X_1$ to conical singularities of $X_2$;
\item
is locally a translation in the translation charts of $X_1$ and $X_2$.
\end{enumerate}

A \emph{square-tiled surface}, also called \emph{origami}, is a translation cover of the standard torus $\RR^2/\ZZ^2$ marked at the origin. Such surfaces are tiled by squares and in particular are Veech surfaces. Square-tiled surfaces can be characterized by their Veech group, according to the following Theorem (for a proof see Appendix C in \cite{hubertlelievre}).

\begin{thmnonnum}[Gutkin-Judge]
A translation surface $X$ is square-tiled
if and only if its Veech group $\veech(X)$ shares a finite-index subgroup with $\slduez$.
\end{thmnonnum}

The following Lemma gives a characterization of square-tiled surfaces in terms of the subgroup $\langle\hol(X)\rangle$ of $\RR^2$ spanned by the elements of $\hol(X)$.

\begin{lem}[Lemma 2.1 in \cite{hubertlelievre}]\label{lem1appendixss2}
A translation surface $X$ is square-tiled if and only
if $\langle\hol(X)\rangle$ is a rank two sublattice of $\ZZ^2$.
\end{lem}

\subsubsection{Reduced origamis}\label{s5ss0sss1}

A square-tiled surface $X$ is called reduced if
$\langle\hol(X)\rangle=\ZZ^2$. Note that this is not always the case, for example $X=\ZZ^2/(\ZZ\oplus2\cdot\ZZ)$ satisfies
$
\langle\hol(X)\rangle=
\ZZ\oplus2\cdot\ZZ
$,
moreover many other examples can be constructed as translation covers of such $X$. For reduced square-tiled surfaces we have a very concrete way to represent Veech groups, according to the following Lemma.

\begin{lem}[Lemma 2.3 in \cite{hubertlelievre}]\label{lem3appendixss2}
For any square-tiled surface $X$ we have
$$
\veech(X)<
\veech
\big(
\RR^2/ \langle\hol(X)\rangle
\big).
$$
In particular, if $X$ is reduced, then $\veech(X)<\slduez$.
\end{lem}

The following Lemma shows that being reduced is preserved by the action of $\slduez$.

\begin{lem}[Lemma 2.4 in \cite{hubertlelievre}]\label{lem4appendixss2}
The $\slduez$-orbit of a reduced square-tiled surface with $N$ squares is the set of reduced square-tiled surfaces with $N$ squares in its $\slduer$-orbit.
\end{lem}

\subsubsection{Action of $\slduez$ for reduced origamis}\label{OrigamiBackgroundActionSL(2,Z)}

Consider the set of generators $\{T,V\}$ of $\slduez$, where
$$
T:=
\begin{pmatrix}
1 & 1\\
0 & 1
\end{pmatrix}
\textrm{ and }
V:=
\begin{pmatrix}
1 & 0\\
1 & 1
\end{pmatrix}.
$$
According to Lemma \ref{lem4appendixss2}, the elements of $\slduez\cdot X$ are those reduced origamis in $\slduer\cdot X$ with the same number of squares as $X$, thus they are a finite set, which is identified with
$$
E(X):=\slduez/\veech(X).
$$
The action of $T$ and $V$ passes to the quotient $E(X)$ and is represented by a oriented graph $\cG(X)$ whose vertices are the elements of $E(X)$ and whose oriented edges correspond to the operations $Y\mapsto T\cdot Y$ and $Y\mapsto V\cdot Y$ for $Y\in E(X)$.

\subsubsection{Cusps}

Let $X$ be a reduced origami, thus in particular a Veech surface. Then $\slduer\cdot X$ is isometric to the unitary tangent bundle of $\HH/\veech(X)$. Any cusp of $\HH/\veech(X)$ is isometric to the quotient of a strip
$
\{z\in\HH;|\Re(z)|<\de,\Im(z)>M\}
$
by the translation $z\mapsto z+\de$, where $\de>0$ is called the \emph{width} of the cusp and $M>0$ is any real number large enough. The cusps of $\HH/\veech(X)$ therefore correspond to conjugacy classes under $\veech(X)$ of its \emph{primitive parabolic elements}, that is the elements in $\veech(X)$ with trace equal to $\pm2$, where primitive means not powers of other parabolic elements of $\veech(X)$.

If $X$ is a reduced origami the eigendirections of parabolic elements of $\veech(X)$ are exactly the elements of $\QQ$. Therefore the cusps of $\HH/\veech(X)$ correspond to equivalence  classes for the homographic action of $\veech(X)$ on $\QQ$ (see \S \ref{s5ss1ssscontinuedfraction}). According to the following Lemma (see Lemma 2.5 in \cite{hubertlelievre}) the cusps of $\HH/\veech(X)$ correspond to the orbits of $T$ on the set $E(X)$ defined in \S \ref{OrigamiBackgroundActionSL(2,Z)}.

\begin{lem}[Zorich]\label{lem5appendixss2}
Let $X$ be a reduced origami. Then the cusps of $\HH/\veech(X)$ are in bijection with the $T$-orbits of $\slduez\cdot X$.
\end{lem}

\subsubsection{Continued fraction}\label{s5ss1ssscontinuedfraction}

Let $\al=a_0+[a_1,a_2,\dots]$ be the \emph{continued fraction expansion} of $\al$, where $a_0\in\ZZ$ and
$a_n\in\NN^\ast$. The $n$-th \emph{Gauss approximation} of $\al$ is the rational number
$$
\frac{p_n}{q_n}:=
a_0+[a_1,a_2,\dots,a_n].
$$

At the projective level, the affine action of $\slduez$ on translation surfaces corresponds to the homographic action on co-slopes, that is
$$
A\cdot \al:=
\frac{a\al+b}{c\al+d}
\textrm{ for }
A=
\begin{pmatrix}
a & b\\
c & d
\end{pmatrix}
\textrm{ and }
\al=\tan(\te),
$$
where $-\infty \leq \alpha \leq + \infty$ is the \emph{co-slope} of  a line which forms an angle $-\pi/2\leq \theta \leq \pi/2$  with respect to the vertical direction measured clockwise. In particular, $\theta=0$ correspond to the vertical direction and to the co-slope $0$, while the horizontal direction on $X$ corresponds to the co-slope $\infty$. The following Lemma describes the relation between the continued fraction and the projective action of $\slduez$.

\begin{lem}\label{lem0s5ss1continuedfraction}
If $\al=a_0+[a_1,a_2,\dots]$ then the sequence of co-slopes $p_n/q_n$ of the Gauss approximations of $\al$ is given by
\begin{eqnarray*}
&&
p_n/q_n=
T^{a_0}\circ V^{a_1}\circ\dots\circ
V^{a_{n-1}}\circ T^{a_n}\cdot 0,
\textrm{ for even }
n
\\
&&
p_n/q_n=
T^{a_0}\circ V^{a_1}\circ\dots\circ
T^{a_{n-1}}\circ V^{a_n}\cdot \infty
\textrm{ for odd }
n.
\\
\end{eqnarray*}
\end{lem}

\begin{dimo}
Just recall that the recursive relations satisfied by the convergents (see for example \cite{kin}) show that the  sequence $(p_n,q_n)$  is obtained by setting $(p_{-2},q_{-2})=(0,1)$ and $(p_{-1},q_{-1})=(1,0)$ and then applying for any $k\in\NN$ the recursive relations
$$
\begin{pmatrix} p_{2k-1} & p_{2k}\\
q_{2k-1} & q_{2k}
\end{pmatrix}=
\begin{pmatrix}
p_{2k-1} & p_{2k-2}\\
q_{2k-1} & q_{2k-2}
\end{pmatrix}
\circ T^{a_{2k}}
\textrm{ and }
\begin{pmatrix}
p_{2k+1} & p_{2k}\\
q_{2k+1} & q_{2k}
\end{pmatrix}=
\begin{pmatrix}
p_{2k-1} & p_{2k}\\
q_{2k-1} & q_{2k}
\end{pmatrix}
\circ V^{a_{2k+1}}.
$$
\end{dimo}

The following Lemma is also well-known (see for example Theorem 19, Chapter 6 in \cite{kin}).

\begin{lem}\label{lem1s5ss1continuedfraction}
Fix $\al$ in $\RR$ and consider a pair of integers $(q,p)$.
$$
\textrm{ If }
|q\al-p|<\frac{1}{2q}
\textrm{ then }
(p,q)=(p_n,q_n)
\textrm{ for some }
n\in\NN.
$$
\end{lem}

\subsection{Renormalized formula for reduced square-tiled surfaces}\label{s5ss1}

Let $X$ be a reduced square-tiled surface. According to Lemma \ref{lem3appendixss2} the lattice
$
\langle\hol(X)\rangle
$
spanned by the relative periods of $X$ equals $\ZZ^2$. Therefore, if $p$ is any conical singularity of $X$, the map
$
\zeta\mapsto\int_p^{\zeta}w
$
is defined modulo elements of
$
\langle\hol(X)\rangle=\ZZ^2
$
and we have a well-defined map
$
\rho:X\to
\TT^2=\RR^2/\ZZ^2
$
setting
\beq\label{eq1s5ss1}
\rho(\zeta):=
\int_p^{\zeta}w
\mod \ZZ^2.
\eeq

\begin{lem}\label{lem1s5ss1}
For any origami $X$ there exists a reduced origami $X_{red}$ in $\glduer\cdot X$ and moreover we have
$$
\cL(X)=\cL(X_{red}).
$$
\end{lem}

\begin{dimo}
Let $\{f_1,f_2\}$ be a basis of $\langle\hol(X)\rangle$ and $G$ be the matrix sending it to the standard basis $\{e_1,e_2\}$ of $\ZZ^2$. Then $X_{red}=G\cdot X$ is reduced. Moreover we have $G=\det(G)\cdot G'$ with $G'\in\slduer$, therefore
$$
\cL(X_{red})=
\cL(G'\cdot X)=
\cL(X),
$$
where the first inequality follows because, on the stratum, the function $X\mapsto a(X)$ is invariant under homothetic transformations (since in the definition \eqref{eq:definitiona} of $a(X)$ we renormalize by $\area(X)$) and the second follows because $G'\cdot X$ belongs to $\slduer\cdot X$ (see Lemma \ref{lem:symmetries}).
\end{dimo}

\subsubsection{Multiplicities}\label{s5ss1sss2}

Let $X$ be a reduced origami and $\rho:X\to \TT^2$ be the covering in Equation (\ref{eq1s5ss1}).

\begin{defi}[Multiplicity of a rational direction]\label{def:multiplicity}
If $\ga:[0,1]\to X$ is a saddle
connection for $X$, we define its \emph{multiplicity} $m(\ga)$ as the degree of
the map $t\mapsto \rho\circ\ga(t)$. The \emph{multiplicity of a rational direction with co-slope} $p/q$ \emph{over the surface} $X$ is the minimal
multiplicity among all saddle connections on $X$ with the same co-slope $p/q$, that is the number $m_X(p/q)$ defined by
$$
m_X(p/q):=
\min
\{m(\ga); \quad \gamma \  \textrm{saddle connection with }\ \hol(\ga)\wedge(p+iq )=0\}.
$$
\end{defi}

Recall from \S \ref{s5ss1ssscontinuedfraction} that the action of $\slduez$ on translation surfaces induces the homographic action on co-slopes. The Remark below is a direct consequence of covariance of the multiplicity under $\slduez$, that is
\beq\label{eq:covariancemultiplicity}
m_X(p/q)=
m_{A\cdot X}\left( A\cdot(p/q) \right).
\eeq

\begin{rem}[Multiplicities are cusp invariant]
If the co-slopes $p/q$ and $p'/q'$ are in the same cusp then there exists some $A\in\veech(X)$ such that
$
p'/q'=A\cdot(p/q)
$
and therefore
$$
m_X(p/q)=m_X(p'/q').
$$
\end{rem}

Fix a reduced origami $X$ and consider the parametrization $L_X:\RR\to \cL(X)$ of its Lagrange Spectrum via the function $\al\to L_X(\al)$ defined by Equation (\ref{eq:standardparametrization}). For convenience of notation we also introduce the function $a_X(\al):=L^{-1}_X(\al)$. Let $N=N(X)$ be the number of squares of $X$, so that $\area(X)=N$.

\begin{lem}\label{lem2s5ss1}
Let $X$ be a reduced origami. We have
$$
L_X(\al)=
N\cdot
\limsup_{q,p\to\infty}
\frac{1}
{m_X^2(p/q)\cdot q\cdot|q\al-p|}.
$$
\end{lem}

\begin{dimo}
Fix a co-slope $\al$ in $\RR$ and set $\te=\arctan(\al)$. Remark that $R_\theta$ sends lines  in direction $\theta$ on $X$ (that is lines which form an angle $\theta$ with the vertical line measured clockwise) to vertical lines on $R_\theta X$. For any period $v$ in $\hol(X)$ call $a_{X,\al}(v)$ its area with respect to $R_\te\cdot X$, that
is
$
a_{X,\al}(v):=
|\Re(e^{i\te}v)|\cdot|\Im(e^{i\te}v)|
$.
If $\ga$ is the saddle connection in $X$ corresponding to the period
$v$ in $\hol(X)$, let $\rho(v)$ in $\hol(\TT^2)$ be the period corresponding to the saddle connection $\rho(\ga)$ in $\TT^2$. By definition of multiplicity of a saddle connection, we have
\begin{equation}\label{multiplicityrelation}
a_{X,\al}(v)=
m_X^2(v)\cdot a_{\TT^{2},\al}(\rho(v)).
\end{equation}
Let us show that if $w$ varies in $\hol(\TT^2)$ we have
\begin{equation}\label{sameasymptotic}
\lim_{|\Im(w)|\to\infty}
\frac
{a_{\TT^2,\al}(w)}
{q\cdot|q\al-p|}=1.
\end{equation}
Given $w \in \hol(\TT^2)$, since $w$ is a period of $\mathbb{T}^2$  we have $w = p+iq$, where $p$ and $q$ are co-prime
integers. Since the line in direction $\theta$ has co-slope $\alpha$, it is parametrized by $(\alpha t, t)$ as $t\in\mathbb{R}$. Thus, in particular it contains the point $(q\alpha, q)$. This shows that $q|q\alpha -p|$ is the area of a parallelogram $P$ through the origin which has  $(q\alpha, q)$ and $(p,q)$ as vertices. Since $a_{\TT^2,\al}(w)$ is the area of a rectangle which can be written as the union of $P$ and two triangles of area $\frac{1}{2}|q\alpha -p|\cos \theta \sin \theta $, the asymptotics in \eqref{sameasymptotic} follows. Thus, letting $v$ vary in $\hol(X)$ and using \eqref{multiplicityrelation} and \eqref{sameasymptotic}, we have
$$
a(R_\te\cdot X)=
\liminf_{|\Im(e^{i\te} v)|\to\infty}
\frac{a_{X,\al}(v)}{\area(X)}=
\frac{1}{N}
\liminf_{q,p\to\infty}
m_X^2(p/q)\cdot q\cdot|q\al-p|.
$$
\end{dimo}

\subsubsection{Renormalized formula for origamis}

The following Lemma, already proved by Perron in 1921, is well-known (see for example \cite{gugu} for the proof).

\begin{lem}\label{lem4s5ss1}
For any $\al$ in $\RR$ we have
$$
L_{\TT^2}(\al)=
\limsup_n
\big(
[a_{n-1},\dots,a_0]+a_n+[a_{n+1},a_{n+2},\dots]
\big).
$$
\end{lem}

Lemma \ref{lem4s5ss1} enables to compute the classical Lagrange spectrum via the continued fraction. In this paragraph we show that the same is possible for reduced origamis under the assumption on $\al$ stated in Theorem \ref{formulasquaretiled} below.  To simplify the notation, set
$$
L(A,n):=
[a_{n},\dots,a_1]+a_{n+1}+[a_{n+2},a_{n+3},\dots],
$$
where $A=(a_k)_{k\in\NN^\ast}$ denotes the sequence of entries of the continued fraction expansion of the irrational slope $\al$, that is $\al=a_0+[a_1,a_2,\dots]$. It is an easy computation (see for example Proposition 1.4, page 58 in \cite{gugu}) to check that that for the $n$-th approximation $p_n/q_n=a_0+[a_1,\dots,a_n]$ of $\al$ we have
$$
L(A,n)=\frac{1}{q_n\cdot|q_n\al-p_n|}.
$$
Set also
$$
M_X(\al):=\max_{n\in\NN}
m_X\bigg(\frac{p_n}{q_n}\bigg).
$$

\begin{thm}\label{formulasquaretiled}
Let $X$ be a reduced origami and $N$ be the number of squares of $X$. If $\al\in\RR$ satisfies
$$
L_{\TT^2}(\al)>2M_X(\al)^2
$$
then we have
$$
L_X(\al)=
N\cdot\limsup_{n\to\infty}
\frac
{[a_{n},\dots,a_1]+a_{n+1}+[a_{n+2},a_{n+3}\dots]}
{m_X^2(p_n/q_n)}.
$$
\end{thm}

\begin{dimo}
Set $M:=M_X(\al)$ and fix $\epsilon>0$ such that
$
L_{\TT^2}(\al)>2M^2+2\epsilon
$.
Since Lemma \ref{lem4s5ss1} gives
$
L_{\TT^2}(\al)=
\limsup_{n}L(A,n)
$
then there are arbitrary large integers $n$ such that
$
L(A,n)>
L_{\TT^2}(\al)-\epsilon
$.
For any such $n$, since
$
1\leq m_X(p_n/q_n)\leq M
$,
we have
$$
\frac
{N}
{m_X^2(p_n/q_n)\cdot q_n\cdot|q_n\al-p_n|}=
N
\frac{L(A,n)}{m_X^2(p_n/q_n)}>
2N
\frac{2M^2+\epsilon}{2M^2}>2N.
$$
On the other hand, consider a pair $(q,p)$ such that $p/q$ does not belong to the sequence of Gauss approximations $p_n/q_n$ of $\al$. Lemma \ref{lem1s5ss1continuedfraction} implies $|q\al-p|>1/(2q)$ and therefore,
since $m_X^2(p/q)\geq 1$ for any $q$ and $p$, we get
$$
\frac
{N}
{m_X^2(p/q)\cdot q\cdot|q\al-p|}<2N.
$$
Finally we apply the factorized formula in Lemma \ref{lem2s5ss1}. The last two inequalities
imply
$$
\limsup_{q,p}\frac
{N}
{m_X^2(p/q)\cdot q\cdot|q\al-p|}=
\limsup_{n}
\frac
{N}
{m_X^2(p_n/q_n)\cdot q_n\cdot|q_n\al-p_n|}.
$$
\end{dimo}

\subsection{Hall's Ray for square-tiled surfaces}\label{s5ss2}
We present here below a constructive proof for the existence of an Hall ray for origamis (see Theorem \ref{thm2s5ss2} below) which  is an adaptation of the original Hall's argument for the classical Lagrange spectrum (\cite{hall}, see also Theorem 3, Chapter 4 in \cite{cusick}).

Let $X$ be a reduced origami. The maximal and the minimal multiplicity over all rational direction of $X$ are respectively
$$
M^+(X):=\max_{p/q\in\QQ}m_X(p/q)
\textrm{ and }
M^-(X):=\min_{p/q\in\QQ}m_X(p/q).
$$

\begin{rem}\label{rem:multiplicitiesarenontrivial}
It is easy to find a reduced origami $X$ with $M^+(X)>1$. For example if $X$ is reduced genus $2$ origami with $N(X)=6$ or $N(X)=7$, that is with $6$ or $7$ squares, then we have $M^+(X)=2$. It is possible to check that the Lagrange Spectra of these examples are non-trivial generalizations of the classical Lagrange Spectrum $\cL$, that is
$$
\cL(X)\not=N(X)\cdot \cL.
$$
On the other hand, we do not know any example of  a reduced square-tiled surface $X$ such that $M^-(X)> 1$.
\end{rem}

Consider the set of generators $\{T,V\}$ for $\slduez$. For any finite word $A=(a_1,\dots,a_n)$ with $n$ elements, where
$a_j\in\NN^\ast$ for any $1\leq j\leq n$, define $g(A)\in\slduez$ setting
\begin{eqnarray*}
&&g(A):=
V^{a_1}\circ T^{a_2}\circ\dots\circ V^{a_{n-1}}\circ T^{a_n}
\textrm{ if }
n
\textrm{ is even }\\
&&g(A):=
V^{a_1}\circ T^{a_2}\circ\dots\circ T^{a_{n-1}}\circ V^{a_n}
\textrm{ if }
n
\textrm{ is odd.}
\end{eqnarray*}

The sequence $(a_1,\dots,a_n)$ correspond to the rational slope $[a_1,\dots,a_n]\in(0,1)$. Applying Equation \eqref{eq:covariancemultiplicity} and Lemma \ref{lem0s5ss1continuedfraction} one gets immediately the following Lemma.

\begin{lem}\label{lem0s5ss2}
For any reduced origami $X$ and any rational slope $[a_1,\dots,a_n]$ we have
\begin{eqnarray*}
&&
m_X([a_1,\dots,a_n])=
m_{g^{-1}(a_1,\dots,a_n)\cdot X}(0)
\textrm{ if }
n
\textrm{ is even }\\
&&
m_X([a_1,\dots,a_n])=
m_{g^{-1}(a_1,\dots,a_n)\cdot X}(\infty)
\textrm{ if }
n
\textrm{ is odd. }
\end{eqnarray*}
\end{lem}

\subsubsection{The graph $\cG_{even}(X)$}\label{sec:even_graph}

Consider the finite set $E(X)$ and let $\cardinalita(E(X))$ be its cardinality. Recall that the cusps of $\HH/\veech(X)$ are in bijection with the orbits of $T$ on $E(X)$, according to Lemma \ref{lem5appendixss2}. For any $Y\in\slduez\cdot X$ denote $w(Y)$ the minimal $w\in\NN^\ast$ such that $T^w\in\veech(Y)$. Observe that we have $V=R\circ T^{-1}\circ R^{-1}$, where $R$ is the counterclockwise rotation by an angle $\pi/2$, that is
$$
R:=
\begin{pmatrix}
0 & -1\\
1 & 0
\end{pmatrix}.
$$
It follows that $w=w(R\cdot Y)$ is the minimal $w\in\NN^\ast$ such that $V^w\in\veech(Y)$. Let $p(X)$ be the maximal cyclic order or elements of $E(X)$ under the action of $T$. Equivalently, $p(X)$ is the maximal width of all cusps of $\HH/\veech(X)$.

Recall from \S~\ref{OrigamiBackgroundActionSL(2,Z)} that the action of $T$ and $V$ on $E(X)$ is represented by the edges of the oriented graph $\cG(X)$. In this section it is more useful to consider a different oriented graph, denoted $\cG_{even}(X)$, whose vertices are the elements of $E(X)$ and whose oriented edges are represented by the action of $g(A)$ for \emph{elementary even words}, that is the words of the form
$$
A=(a_1,a_2)
\textrm{ with }
1\leq a_1<p(X)
\textrm{ and }
1\leq a_2<p(X).
$$
The graph $\cG_{even}(X)$ is a finite graph, non necessarily connected. It is easy to verify that the action of elementary even words is \emph{recurrent} on $\cG_{even}(X)$, that is for any $Y_1$ and $Y_2$ in $E(X)$ and any finite sequence of elementary even words $A_1,\dots,A_k$ such that
$
Y_2=g(A_1\ast\dots\ast A_k)\cdot Y_1
$
there exists a finite sequence of elementary even words $B_1,\dots,B_l$ such that
$
Y_1=g(B_1\ast\dots\ast B_l)\cdot 2_1
$ (since if there is path from $Y_2$ to $Y_1$ of odd length, one can concatenate it with a loop of odd length from $Y_1$ to itself).

Define $E_{even}(X)$ as the subset of $E(X)$ of the vertices of $\cG_{even}(X)$ which belong to the same connected component of $X$. Fix an element $X_{min}\in E(X)$ such that the vertical direction has minimal multiplicity, that is
$
m_{X_{min}}(0)=M^-(X)
$.
Observing that $\cL(X)=\cL(X')$ for any $X'\in\slduez \cdot X$, modulo changing the initial origami $X$ inside $\slduez\cdot X$ we can assume that
$$
X_{min}\in E_{even}(X).
$$

\begin{lem}\label{lem1s5ss2}
For any two even length words $A=(a_1,\dots,a_{2r})$ and $B=(b_1,\dots,b_{2s})$ there exists an integer $i\leq \cardinalita(E(X))$ and a finite even length word $D=(d_1,\dots,d_{2i})$ whose entries satisfy
$$
d_j< p(X)
\textrm{ for }
j=1,\dots,2i
$$
such that the co-slope corresponding to the concatenated word $A\ast D \ast B$ satisfies
$$
m_X([A\ast D\ast B])=M^-(X).
$$
\end{lem}

\begin{dimo}
According to Lemma \ref{lem0s5ss2}, it is enough to find a finite word $D$ as in the statement such that the concatenation $A\ast D\ast B$, which has even length $2(r+s+i)$, satisfies
$
g(A\ast D\ast B)\cdot X_{min}=X
$.
Since
$
g(A\ast D\ast B)=g(A)\cdot g(D)\cdot g(B)
$
then the last condition is equivalent to
$$
g(D)\cdot\big(g(B)\cdot X_{min} \big)=g(A)^{-1}\cdot X,
$$
thus the required word $D$ exits because both $g(B)\cdot X_{min}$ and $g(A)^{-1}\cdot X$ belong to the same connected component of $\cG_{even}(X)$. Observe that $\cG_{even}(X)$ has at most as many vertices as $\cG(X)$, thus the bound on $i$ follows. Finally the bound on the entries $d_j$ of $D$ is trivial in view of the considerations at the beginning of this section \ref{sec:even_graph}.
\end{dimo}

\subsubsection{Constructive proof of the existence of Hall's ray}

\begin{thm}\label{thm2s5ss2}
Let $X$ be a reduced origami. Then $\cL(X)$ contains the half-line $[r(X),+\infty)$, where
$$
r(X)=
\frac{N}{M^-(X)^2}\cdot
\max \{2M^+(X)^2+1,7,p(X)+2\}.
$$
\end{thm}

\begin{dimo}
Let $x$ be any real number with
$
x>r(X)\cdot M^-(X)^2\cdot N^{-1}
$,
where $r(X)$ is as in the statement. According to a classical result due to Hall (\cite{hall}, see also Theorem 1, Chapter 4 in \cite{cusick}) there exist $x_0\in\NN$ and two sequences $A=(a_k)_{k\in\NN^\ast}$ and $B=(b_k)_{k\in\NN^\ast}$ with $a_k,b_k$ in $\{1,2,3,4\}$ for any $k\in\NN$ such that
$$
x=
[a_1,a_2,\dots]+
x_0+
[b_1,b_2,\dots].
$$
Observe that the assumption on $x$ and the definition of $r(X)$ imply that
$
x_0\geq 2M^+(X)^2
$.
We define inductively an infinite word $C=(c_n)_{n\in\NN^\ast}$ and a subsequence  $\{n(m)\}_{m\in\NN}$ such that for any $m$ we have
\begin{equation}%narray}
\label{eq1thm2s5ss2}
m_X\big([c_1,\dots,c_{n(m)}]\big)= M^-(X) \quad \mathrm{and} \quad c_{n(m)+1} = x_0.
%\\c_{n(m)+1}&=&x_0 \label{x0}
\end{equation}
%\end{eqnarray}

\emph{Step one.} Set $A:=(x_0,a_1)$ and $B:=(b_2,b_1)$. Let $i(1)\leq \cardinalita(E(X))$ be the positive integer and $D=(d^1_1,\dots,d^1_{2i(1)})$ be the finite word with
$d^1_j\leq p(X)$ for any $1\leq j\leq 2i(1)$ provided by Lemma \ref{lem1s5ss2}. Condition \eqref{eq1thm2s5ss2} is satisfied if we define the initial block $(c_1,\dots,c_{n(1)})$ of $C$ setting
\begin{equation*}
n(1):=4+2i(1)
\quad \textrm{ and }\\
(c_1,\dots,c_{n(1)}, ,c_{n(1)+1}):=
(x_0,a_1,d^1_1,\dots,d^1_{2i(1)},b_1,b_2,x_0).
\end{equation*}

\emph{Inductive step.} Suppose that the integer $n(m)$ is defined and $c_k$ are defined for any
$1\leq k\leq n(m)+1$ in order to satisfy Condition \ref{eq1thm2s5ss2}. Consider the finite words
\begin{eqnarray*}
&&
A:=(c_1,\dots,c_{n(m)},x_0,a_1,\dots,a_{2m+1})\\
&&
B:=(b_{2m+2},\dots,b_1,x_0)
\end{eqnarray*}
Let $i(m+1)\leq \cardinalita(E(X))$ be the positive integer and $D=(d^{m+1}_1,\dots,d^{m+1}_{2i(m+1)})$ be the corresponding finite word with $d^{m+1}_j\leq p(X)$ for any $1\leq j\leq 2i(m+1)$ which are provided by Lemma \ref{lem1s5ss2}. We set
\begin{eqnarray*}
&&
n(m+1):=n(m)+4(m+1)+2i(m+1)
\quad \textrm{ and }\\
&&
(c_{n(m)+1},\dots,c_{n(m+1)+1}):=
(x_0,a_1,\dots,a_{2m+1},d^{m+1}_1,\dots,
d^{m+1}_{2i(m+1)},b_{2m+2},\dots,b_1,x_0).
\end{eqnarray*}
The finite word $(c_1,\dots,c_{n(m+1)+1})$ is equal to the concatenation
$
(c_1,\dots,c_{n(m)})\ast(c_{n(m)+1},\dots,c_{n(m+1)+1})
$
and satisfies Condition \eqref{eq1thm2s5ss2}. The inductive construction of the infinite word $C$ is therefore complete and we set $\al:=[c_1,c_2,\dots]$, which is of course irrational and of bounded type.

\emph{End of the proof.} Observe that
$
L_{\TT^2}(\al)>x_0\geq 2M^+(X)^2
$.
Thus, Theorem \ref{formulasquaretiled} can be applied and it implies
$$
L_X(\al)=
N\cdot
\limsup_{n\to\infty}
\frac
{L(C,n)}
{m_X^2([c_1,\dots,c_n])}.
$$
Observe that for any integer $m$ Condition \eqref{eq1thm2s5ss2} implies
$$
N\cdot
\frac
{L(C,n(m))}
{m_X^2([c_1,\dots,c_{n(m)}])}>
N\cdot\frac{x_0}{M^-(X)^2}.
$$
On the other hand, since for $n\not=n(m)+1$ we have
$
1\leq c_n\leq \max\{4,p(X)-1\},
$
for all $n\not=n(m)$ we have
$$
N\cdot
\frac
{L(C,n)}
{m_X^2([c_1,\dots,c_n])}<
N\cdot\frac{2+\max\{4,p(X)-1\}}{M^-(X)^2}.
$$
The assumption on $x$ and the definition of $r(X)$ imply that
$
x_0\geq \max\{6,p(X)+1\}.
$
Thus,  it follows that
\begin{eqnarray*}
L_X(\al)&=&
N\cdot
\limsup_{m\to\infty}
\frac
{L(C,n(m))}
{m_X^2([c_1,\dots,c_{n(m)}])}=
\frac{N}{M^-(X)^2}\cdot
\limsup_{m\to\infty}L(C,n(m))
\nonumber \\
&=&
\frac{N}{M^-(X)^2}\cdot
\limsup_{m\to\infty}
\big(
[c_{n(m)},\dots,c_1]+
c_{n(m)+1}+
[c_{n(m)+2},c_{n(m)+3},\dots]
\big)
\nonumber \\
&=& \frac{N}{M^-(X)^2}\cdot
\limsup_{m\to\infty}
\big(
[b_1,\dots,b_m\dots]+x_0+[a_1,\dots,a_m,\dots]\big)
\nonumber \\
&=& \frac{N}{M^-(X)^2}\cdot
\lim_{m\to\infty}
\big(
[b_1,\dots,b_m]+x_0+[a_1,\dots,a_m]\big)=
\frac{N\cdot x}{M^-(X)^2}.
\end{eqnarray*}
\end{dimo}

\appendix

\section{Lower bound for spectra}\label{appendixAuxilSs1}

Here we prove Lemma \ref{lem:lowerbound}, that is for an invariant locus $\cI$ contained in the stratum $\cH(k_1,\dots, k_r)$ of translation surfaces of genus $g$, where $2g-2=k_1+\dots+k_r$, we have
$$
\cL(\cI)
\subset
\left[
{\pi\cdot\frac{2g+r-2}{2} },+\infty
\right].
$$

\begin{dimo}
Fix any $X \in \cI$ and, up to rescaling, assume that $\area(X)=1$. Recall that by Vorobets' identity  we have that  $a^{-1}(X) = 2s^{-2}(X)$  (see Proposition \ref{prop:vorobetsidentity}). Thus, the Lemma follows immediately from the definition \ref{eq:definitions} of $s(X)$ if we show that for any $Y$ in $\cH^{(1)}(k_1,\dots, k_r)$ we have
$$
\frac{\sys^2(Y)}{2} \leq \frac{2}{\pi(2g-2+r)}.
$$
The proof of this upper bound  for the systole function follows from the following argument given by Smillie and Weiss in \cite{smillieweiss}\footnote{The bound given in \cite{smillieweiss} is unfortunately incorrect because of a typo in the arguments of the proof, but the proof that we reproduce here is due to them.}. For each conical singularity $p_i$, $1\leq i \leq r$, consider the set $B_i \subset Y$ of points which have flat distance less than $\sys(Y)/2$ from $p_i$. If $p_i$ has conical angle $2\pi (k_i+1)$ then
$
\area(B_i)=(k_i +1)\pi \sys(Y)^2/4
$
and by definition of systole, $B_i$ are pairwise disjoint. Thus we have
$
\sum_{i=1}^r  (k_i +1)\sys(Y)^2 \pi/4\leq 1
$,
which gives the desired upper bound since
$
\sum_{i=1}^r k_i=2g-2
$.
\end{dimo}

\section{Closing and Shadowing Lemma for affine invariant loci}\label{loci}

The main results of this Appendix are Proposition~\ref{prop:appendixclosinglemma} and Proposition~\ref{prop:appendixshadowinglemma}, which provide a combinatorial version respectively of the closing Lemma and of the Shadowing Lemma for closed $\slduer$-invariant loci, adapted to zippered rectangles coordinates. These two results are used in \S~\ref{s3} to prove the closure of the Lagrange spectrum and density of periodic orbits for any such locus (Theorem \ref{thm:closure}). They show that, given a path which describe a zippered rectangle in an invariant locus,  the combinatorial version of the closing and shadowing lemmas used in \S~\ref{s3} produce new paths in the same $\slduer$-invariant locus. In \S~\ref{sec:Rauzy_loci} we first describe the structure of $\slduer$-invariant loci in terms of zippered rectangles coordinates, which follows from Eskin-Mizhakhani's work \cite{eskinmirza}. In \S~\ref{sec:hypboxes} we prove some technical Lemmas, providing nice flow-boxes, or \emph{hyperbolic boxes}, in the zippered rectangles coordinates. Finally in \S~\ref{sec:closing_loci} we prove Proposition~\ref{prop:appendixclosinglemma} and Proposition~\ref{prop:appendixshadowinglemma}

\subsection{Rauzy-Veech induction and affine loci}\label{sec:Rauzy_loci}

Let $\cH(k_1,\dots,k_r)$ be a stratum of translation surfaces. Up to considering a finite covering, we can assume that it is an affine manifold without orbifold points. Let $\cC$ be a connected component of the stratum and $\cR$ be a Rauzy class representing the surfaces in $\cC$. Rauzy paths and combinatorial data are considered in $\cR$. If $\ga$ is a Rauzy path denote $\widehat{Q}_\ga$ the corresponding linear branch of the Rauzy-Veech map. If $\ga$ is a positive loop $\ga$ at $\pi\in\cR$, denote respectively $\la^\ga$ and $\tau^\ga$ the maximal and the minimal eigenvector of $^tB_\ga^{-1}$, normalized so that $\|\la^\ga\|=1$ and
$
\area(\pi,\la^\ga,\tau^\ga)=1
$.
The $\cF_t$-orbit of such data is a closed geodesic in the stratum. In terms of the Rauzy-Veech induction we have
$$
\widehat{Q}_\ga(\pi,\la^\ga,\tau^\ga)=
\|^tB_\ga^{-1}\la^\ga\|\cdot
(\pi,\la^\ga,\tau^\ga).
$$

\subsubsection{Affine invariant loci}

Let $\cI$ ba a closed subset of $\cC$ which is invariant under $\slduer$. For any combinatorial datum $\pi\in\cR$ let $\cI_\pi$ be the set of data
$
(\la,\tau)\in\RR_+^\cA\times\Th_\pi
$
such that the surface $X(\pi,\la,\tau)$ belongs to $\cI$. Fix data $(\la,\tau)\in\cI_\pi$ representing the surface $X(\pi,\la,\tau)$. Let
$
t\mapsto A_t:=
\begin{pmatrix}
a_t & b_t\\
c_t & d_t
\end{pmatrix}
$
be a continuous path in $\slduer$ with $A_0=Id$ and $I\subset \RR$ be an interval containing $0$ such that

\beq\label{eqCondizioneAppartenenzaLuogoInv}
(a_t\la+b_t\tau,c_t\la+d_t\tau)
\in
\RR_+^\cA\times\Th_\pi
\textrm{ for any }
t\in I.
\eeq
Then the data
$
(a_t\la+b_t\tau,c_t\la+d_t\tau)
$
belong to $\cI_\pi$ for any $t\in I$ and represent the surface
$
A_t\cdot X(\pi,\la,\tau)
$.
The celebrated result of Eskin-Mirzakani \cite{eskinmirza} can be stated in terms of Rauzy-Veech induction as follows.

\begin{thm}[Eskin-Mirzakani]\label{thmEskinMirzakaniRV}
For any $\pi\in\cR$ there exists a countable family $\cV_\pi=(V_j)$ of vector subspaces of $\RR^\cA\times\RR^\cA$ of the same dimension as $\cI$ satisfying the properties above.
\begin{enumerate}
\item
We have
$
\cI_\pi=
\big(\bigcup_{V_j\in\cV_\pi} V_j\big)
\cap
\big(\RR_+^\cA\times\Th_\pi\big)
$,
that is any
$
C_j:=V_j\cap \big(\RR_+^\cA\times\Th_\pi\big)
$
is a convex open subcone of $V_j$ of the same dimension as $\cI$, providing local affine coordinates for the locus.
\item
For any $i$ and $j$ we have
$
V_i\cap V_j\cap\big(\RR_+^\cA\times\Th_\pi\big)
=\emptyset
$
if $i\not=j$. Indeed points in a non-trivial intersection $V_i\cap V_j$ would represent elements of $\cI$ not admitting a local affine chart.
\item
For any loop $\ga$ at $\pi$ (not necessarily positive) and any $V_j\in\cV_\pi$ intersecting $\De_\ga\times\Th_\pi$ there exists $V_i\in\cV_\pi$ such that $\widehat{Q}_\ga(V_j)=V_i$, that is the Rauzy-Veech induction represents affine coordinate changes for $\cI$.
\end{enumerate}
\end{thm}

The following Lemma is immediate.

\begin{lem}\label{lemCondizioneAppartenenzaLuogoInv}
Consider $V\in\cV_\pi$ and data
$
(\la,\tau)\in C=V\cap\cI_\pi
$.
Let $t\mapsto A_t$ be a continuous path in $\slduer$ with $A_0=Id$ and $I$ be an interval with $0\in I$ satisfying condition
(\ref{eqCondizioneAppartenenzaLuogoInv}). Then
$$
(a_t\la+b_t\tau,c_t\la+d_t\tau)\in C
\textrm{ for any }
t\in I.
$$
\end{lem}

\begin{rem}\label{remInvarianceVectorSubspaces}
Consider $V\in\cV_\pi$ and the corresponding convex open subcone $C=V\cap\cI_\pi$. The following facts hold.
\begin{enumerate}
\item
If $\ga$ is a (not necessarily positive) loop at $\pi$ and we have data
$
(\la,\tau)\in C\cap(\De_\ga\times\Th_\pi)
$
such that
$
(^tB^{-1}_\ga\la,^tB^{-1}_\ga\tau)\in C
$,
then $\widehat{Q}_\ga V=V$.
\item
In particular, if $C$ contains pseudo-Anosov data $(\pi,\la^\ga,\tau^\ga)$ corresponding to a positive loop $\ga$ at $\pi$ then $\widehat{Q}_\ga V=V$.
\item
For any $(\la,\tau)\in C$ and any $t\in\RR$ we have
$
(e^t\la,e^{-t}\tau)\in\RR_+^\cA\times\Th_\pi
$,
thus Lemma \ref{lemCondizioneAppartenenzaLuogoInv} implies
$
(e^t\la,e^{-t}\tau)\in C
$
for all $t\in\RR$.
\end{enumerate}
\end{rem}

\subsection{Hyperbolic boxes in the space of zippered rectangles}\label{sec:hypboxes}

Let $\ga$ be a finite path. Let $\De^{(1)}_\ga$ be the simplex of those $\la\in\De_\ga$ normalized by $\sum_\chi\la_\chi=1$. Let $\Th^{(1)}_\ga$ be the codimension one hypersurface of those $\tau\in\Th_\ga$ normalized by $\sum_\chi|\tau_\chi|=1$.

\begin{defi}\label{defhyperbolicboxes}
Consider a positive path $\nu$ ending at $\pi$ and a positive path $\de$ starting at $\pi$, so that the concatenation $\nu\ast\de$ is defined. Let $I\subset\RR$ be any interval. An \emph{hyperbolic box} $\scatola_\pi(I,\de,\nu)$ is the image  of $I\times\De^{(1)}_\de\times\Th^{(1)}_\nu$ into the area-one hyperboloid in
$
\{\pi\}\times\De_\de\times\Th_\nu
$
via the map
$$
(t,\la,\tau)\mapsto
\frac{1}{\area(\pi,\la,\tau)}
\cdot
(\pi,e^t\la,e^{-t}\tau),
$$
satisfying the extra property that any orbit of the Rauzy-Veech map meets
$
\scatola_\pi(I,\de,\nu)
$
at most once.
\end{defi}

To simplify notation we always assume that the path $\de$ and $\nu$ defining an hyperbolic box have same length $m$. Recall that for any Rauzy path $\ga$, either finite, half-infinite or bi-infinite, in Equation \eqref{eq1sec:finitetime} we define $N(\ga)$ as the supremum of $\|B_\nu\|$ over all positive sub-paths $\nu$ of $\ga$ not admitting a proper positive sub-path.

\begin{lem}\label{lem1HypBoxes}
For any $M>0$ there exists a positive integer $m_1=m_1(M)$  such that for any integers $r$ and $m$ such that $r>m_1$ and $m_1 \leq m \leq r$ the following holds. If $\ga(-\infty,+\infty)$ is a bi-infinite path with
$
N\big(\ga(-r,+\infty)\big)<M
$
and if we set
$$
\de:=\ga(1,m)
\textrm{ and }
\nu:=\ga(-m+1,0)
$$
then there exists some interval $I\subset\RR$ such that $\scatola_\pi(I,\de,\nu)$ is an hyperbolic box.
\end{lem}

\begin{dimo}
Given $m\in\NN$ and a bi-infinite path $\ga(-\infty,+\infty)$, let us denote by $\de_m:=\ga(1,m)$ and $\nu_m:=\ga(-m+1,0)$. Let $m_0$  be a positive integer such that the path $\de_{m_0}$  is positive. Since $\De_{\de_{m_0}}^{(1)}$ is compact, there exists $0<c<1$ such that for any $\la\in\De^{(1)}_{\de_{m_0}}$ we have $\la_\chi>c\|\la\|=c$ for any $\chi\in\cA$. According to Lemma \ref{lem:contraction}, when we increase $m \in \mathbb{N}$ the supremum of the diameters of $\De^{(1)}_{\de_m}$ and $\Th^{(1)}_{\nu_m}$ over all paths such $N\big(\ga(-m,+\infty)\big)<M$ shrinks exponentially as a function of $m$. Moreover $\area(\pi,\la,\tau)$ is a continuous function of $(\pi,\la,\tau)$, bounded away from zero, therefore there exists $m_1\geq m_0$ such that the condition
$$
N\big(\ga(-r,+\infty)\big)<M
\textrm{ for some }
r>m_1
$$
implies that for any  $m_1 \leq m \leq r$ and any pair of data $(\la,\tau)$ and $(\la',\tau')$ in
$
\De^{(1)}_{\de_{m}} \times\Th^{(1)}_{\nu_{m}}
$
we have
\begin{equation}\label{eq1lem1HypBoxes}
\sqrt{1-c}
<
\frac
{\area(\pi,\la',\tau')}
{\area(\pi,\la,\tau)}<
\frac{1}{\sqrt{1-c}}.
\end{equation}
Thus, if we choose the interval $I$ small enough so that $e^{-|s-t|}>\sqrt{1-c}$ for any $s,t\in I$, let us show that $\scatola_\pi(I,\de,\nu)$ where $\nu=\nu_{m}$ and $\de=\de_{m}$ is an hyperbolic box. Since  by definition $m \geq m_0$, we have that $\de$ is positive, so we just need to show that there cannot be an orbit of $Q$ meeting
$
\scatola_\pi(I,\de,\nu)
$
twice. If by contradiction this is not the case, there exists  data $(\pi,\la,\tau)$ and $(\pi,\la',\tau')$ in
$
\{\pi\}\times \De^{(1)}_\de\times\Th^{(1)}_\nu
$
and $s,t\in I$ and an integer $n\geq 1$ such that
$$
\frac{1}{\area(\pi,\la',\tau')}
(\pi,e^s\la',e^{-s}\tau')=
\frac{1}{\area(\pi,\la,\tau)}
(\pi^{(n)},e^t\la^{(n)},e^{-t}\tau^{(n)}).
$$
Since
$
\De_{\de_{m}}  \subset\De_{\delta_{m_0}}
$
then $\la_\chi>c\|\la\|=c$ for any $\chi\in\cA$ and therefore
$
\|\la^{(n)}\|
\leq
\|\la^{(1)}\|=
1-\min_\chi\la_\chi<1-c
$.
It follows that
$$
e^{s-t}
\frac
{\area(\pi,\la,\tau)}
{\area(\pi,\la',\tau')}<1-c.
$$
Since  $e^{-|s-t|}>\sqrt{1-c}$  by choice of $I$, this contradicts \eqref{eq1lem1HypBoxes}. Thus $\scatola_\pi(I,\de,\nu)$ is a hyperbolic box.
\end{dimo}

\subsubsection{Hyperbolic boxes and affine loci}

Let  $\cI$  be a closed $SL(R, \mathbb{R})$ invariant subset of $\cC$.

\begin{defi}
An hyperbolic box $\scatola_\pi(I,\de,\nu)$ is \emph{adapted} to $\cI$ if there exists an unique $V\in\cV_\pi$ such that
$
V\cap\scatola_\pi(I,\de,\nu)\not=\emptyset
$.
\end{defi}

If $\scatola_\pi(I,\de,\nu)$ is an adapted hyperbolic box, $V\cap\scatola_\pi(I,\de,\nu)$ is a convex disc of the same dimension as $\cI^{(1)}$, that is the image of $\scatola_\pi(I,\de,\nu)$ is a small open set in $\cH^{(1)}(k_1,\dots,k_r)$ which intersects $\cI^{(1)}$ in an open disc.

Consider data $(\pi,\la,\tau)\in\cI_\pi$ with $\area(\pi,\la,\tau)=1$ and let $\ga(-\infty,+\infty)$ be the corresponding bi-infinite path. Observe that if $\de$ and $\nu$ have length $m$ and $\scatola_\pi(I,\de,\nu)$ contains $(\pi,\la,\tau)$ then $\la\in\De_\de$ and $\tau\in\Th_\nu$ and therefore
$$
\ga(1,m)=\de
\textrm{ and }
\ga(-m+1,0)=\nu.
$$

\begin{lem}\label{lem2HypBoxes}
For any $M>0$ there exists $r_1=r_1(M)\in\NN$ such that the following holds. If $(\pi,\la,\tau)$ are data with $(\la,\tau)\in\cI_\pi$ and $\area(\pi,\la,\tau)=1$, that induce a bi-infinite path $\ga(-\infty,+\infty)$ with
$$
N\big(\ga(-r,+\infty)\big)<M
\textrm{ for some }
r\geq r_1
$$
then there exists an hyperbolic box
$
\scatola_\pi(I,\de,\nu)
$
adapted to $\cI$ which contains $(\pi,\la,\tau)$, where $\de$ and $\nu$ have length bounded by $r_1$.
\end{lem}

\begin{dimo}
Fix $M>0$ and let $m_1=m_1(M)$  be the  positive integer given by Lemma \ref{lem1HypBoxes}. Consider $r>m_1$ such that
$
N\big(\ga(-r,+\infty)\big)<M
$.
Then, if we denote as before by $\de_m:=\ga(1,m)$ and by $\nu_m:=\ga(-m+1,0)$, for any $m_1 \leq m < r$ there exists a proper interval $I$ such that
$
(\pi,\la,\tau)\in\scatola_\pi(I,\de_m,\nu_m)
$.
So we just need to show that there exists $r_1 > m_1$ which depends only on $M$ such that if $r >r_1$ than the hyperbolic box $\scatola_\pi(I,\de_{r_1},\nu_{r_1})$ is also adapted to $\cI$.  Since $\ga(1,r_1)$ and $\ga(-r_1+1,0)$ are positive since $r_1 >m_1$ (see the proof of  Lemma \ref{lem1HypBoxes}),
$
\De_{\de_{r_1}}^{(1)} \times \Theta_{\nu_{r_1}}^{(1)} \subset \cK
$
for some   open sub-cone $\cK\subset\RR^\cA_+\times\Th_\pi$
with compact closure with respect to the Hilbert metric of $\RR^\cA_+\times\Th_\pi$ (see \S~\ref{subsubsection:hilbertmetric}), which depends only on $r_1(M)$ and hence only on $M$.

Recall that $\cI^{(1)}$ is an affine sub-manifold of $\cH^{(1)}(k_1,\dots,k_r)$, thus locally it is the zero locus of some analytic function. It follows that there exists just finitely many hyperplanes $V\in\cV_\pi$ such that $V\cap\cK\not=\emptyset$. Hence there exists some $\epsilon>0$ such that the $\epsilon$-neighborhoods of $V\cap\cK$ with respect to the Hilbert metric are disjoint each other. Thus, reasoning as in the proof of Lemma \ref{lem2HypBoxes}, there exists $r_1=r_1(M)>r_0$  such that for any bi-infinite path $\ga(-\infty,+\infty)$ with
$
N\big(\ga(-r_1,+\infty)\big)<M
$
the diameter of
$
\De_{\de_{r_1}}\times\Th_{\nu_{r_1}}$
a smaller than $\epsilon/2$ with respect to the Hilbert metric. In this case, the hyperbolic box $\scatola_\pi(I,\de_{r_1},\nu_{r_1})$ can intersect at most one $V\in\cV_\pi$ and the Lemma is thus proved.
\end{dimo}

\begin{lem}\label{lem3HypBoxes}
Fix $\pi\in\cR$ and $V\in\cV_\pi$ and consider the open subcone $C:=V\cap\cI_\pi$. Let
$
\scatola_\pi(\de,\nu,I)
$
be an hyperbolic box adapted to $\cI$ with non-empty intersection with $C$, where $\de$ and $\nu$ are positive path of length $m$. If $\ga(0,r)$ is a loop at $\pi$ and moreover we have the extra conditions
$$
\ga(r-m+1,r)=\nu
\textrm{ and }
\ga(r+1,r+m)=\de
$$
then
$
(\pi,\la^{(r)},\tau^{(r)})
=
\widehat{Q}_{\ga(0,r)}(\pi,\la,\tau)
\in C
$.
\end{lem}

\begin{dimo}
The assumption in the statement is equivalent to
$\la^{(r)}\in\De_\de$ and $\tau^{(r)}\in\Th_\nu$, therefore there exists $t\in\RR$ such that
$
(\pi,e^{t}\la^{(r)},e^{-t}\tau^{(r)})
$
belongs to $\scatola_\pi(I,\de,\nu)$. On the other hand $(\la^{(r)},\tau^{(r)})$ belongs to $\cI_\pi$, since $\cI_\pi$ is invariant under $\widehat{Q}_{\ga(0,r)}$ and contains $(\la,\tau)$. Moreover also $(e^{t}\la^{(r)},e^{-t}\tau^{(r)})$ belongs to $\cI_\pi$, since condition \eqref{eqCondizioneAppartenenzaLuogoInv} is satisfied, thus we get
$$
(\pi,e^{t}\la^{(r)},e^{-t}\tau^{(r)})
\in
\cI_\pi\cap\scatola_\pi(I,\de,\nu)=C.
$$
Finally part (3) of Remark \ref{remInvarianceVectorSubspaces} implies
$
(\la^{(r)},\tau^{(r)})\in C
$.
\end{dimo}

\subsection{Closing and Shadowing Lemma}\label{sec:closing_loci}

\subsubsection{The main geometrical Lemma}

\begin{lem}\label{lemMaingeomlemma}
Fix $\pi\in\cR$ and $V\in\cV_\pi$ and consider the convex open subcone $C=V\cap\cI_\pi$. Consider a finite path $\eta$ ending at $\pi$ and a finite path $\ga$ starting at $\pi$, so that the concatenation $\eta\ast\ga$ is possible. If $C$ has non-empty intersection both with $\De_\ga\times\Th_\pi$ and with $\RR_+^\cA\times\Th_\eta$ then we have
$$
C\cap(\De_\ga\times\Th_\eta)
\not=
\emptyset.
$$
\end{lem}

\begin{dimo}
Let
$
(\la,\tau)
\in
C\cap(\De_\ga\times\Th_\pi)
$
and
$
(\la',\tau')
\in
C\cap(\RR_+^\cA\times\Th_\eta)
$.
Part (2) of Remark \ref{remInvarianceVectorSubspaces} and convexity of $C$ imply that for any $t\in\RR$ and $s\in(0,1)$ we have
$$
(1-s)\cdot(\la,\tau)
+
s\cdot(e^{-t}\la',e^t\tau')\in C.
$$
Since $\la$ and $\tau$ belong respectively to $\De_\ga$ and to $\Th_\ga$, and the latter are open subcones of $\RR^\cA$, then taking $s$ small enough and $t$ such that $se^t$ is big enough, we have
$$
(1-s)\cdot\la+se^{-t}\cdot\la'\in\De_\ga
\textrm{ and }
(1-s)\cdot\tau+se^t\cdot\tau'\in\Th_\ga,
$$
thus the data
$
(1-s)\cdot(\la,\tau)
+
s\cdot(e^{-t}\la',e^t\tau')
$
belong to
$
C\cap(\De_\ga\times\Th_\ga)
$.
\end{dimo}

\begin{cor}\label{corMainGeomCor}
Fix $\pi\in\cR$ and $V\in\cV_\pi$ and consider the convex open subcone $C=V\cap\cI_\pi$. Let $\ga_1,\dots,\ga_k$ be positive loops at $\pi$ such that for any $i=1,\dots,k$ the following holds
\begin{enumerate}
\item
$C$ has non-empty intersection both with
$
\De_{\ga_i}\times\Th_\pi
$
and with
$
\RR_+^\cA\times\Th_{\ga_i}
$.
\item
We have $\widehat{Q}_{\ga_i}(V)=V$.
\end{enumerate}
Then for any $i=1,\dots,k$ we have
$$
C\cap
(\De_{\ga_{i+1}\ast\dots\ast\ga_k}\times
\Th_{\ga_1\ast\dots\ast\ga_i})
\not=\emptyset.
$$
\end{cor}

\begin{dimo}
Observe that
$
\widehat{Q}_{\ga_1\ast\dots\ast\ga_i}
$
is a linear bijection from
$
\De_{\ga_1\ast\dots\ast\ga_k}\times\Th_\pi
$
onto
$
\De_{\ga_{i+1}\ast\dots\ast\ga_k}
\times
\Th_{\ga_1\ast\dots\ast\ga_i}
$,
and since assumption (2) in the statement is equivalent to
$
\widehat{Q}_{\ga_{1}\ast\dots\ast\ga_k}(V)=V
$,
then it is enough to prove that for any $k>0$ we have
\beq\label{eq1propMainGeomProp}
C\cap
(\De_{\ga_1\ast\dots\ast\ga_k}\times\Th_\pi)
\not=\emptyset.
\eeq
For $k=1$ the required condition corresponds to assumption (1) in the statement. The general case is proved by induction on $k$. Suppose that condition (\ref{eq1propMainGeomProp}) holds for $k$ positive loops $\ga_1,\dots,\ga_k$ and consider an extra loop $\ga_{k+1}$ at $\pi$ satisfying the assumptions in the statement. According to Lemma \ref{lemMaingeomlemma} we have
$$
C\cap
(\De_{\ga_{k+1}}\times\Th_{\ga_1\ast\dots\ast\ga_k})
\not=\emptyset.
$$
Moreover $\widehat{Q}_{\ga_1\ast\dots\ast\ga_k}$ maps
$
\De_{\ga_1\ast\dots\ast\ga_k\ast\ga_{k+1}}\times\Th_\pi
$
onto
$
\De_{\ga_{k+1}}\times\Th_{\ga_1\ast\dots\ast\ga_k}
$
and since
$
\widehat{Q}_{\ga_i}(V)=V
$
for any $i=,\dots,k$ then condition (\ref{eq1propMainGeomProp}) follows.
\end{dimo}

\subsubsection{Closing Lemma}

\begin{prop}\label{prop:appendixclosinglemma}
Fix $\pi\in\cR$. For any $M>0$ there exists a positive integer $r_1=r_1(M)$ such that for any $m>r_1$ the following holds. If $(\pi,\la,\tau)$ are data with $(\la,\tau)\in\cI_\pi$ that induce a bi-infinite path $\ga(-\infty,+\infty)$ with
$
N(\ga(-m,+\infty))<M
$
and for $r>0$ the segment $\ga:=\ga(0,r)$ is a positive loop at $\pi$ satisfying the extra condition
$$
\ga(r-m+1,r+m)=\ga(-m+1,m)
$$
then $(\pi,\la^\ga,\tau^\ga)\in\cI_\pi$, that is the closed geodesic corresponding to $\ga=\ga(0,r)$ is contained in $\cI$.
\end{prop}

\begin{dimo}
We can assume that $\area(\pi,\la,\tau)=1$, indeed the induced bi-infinite path $\ga(-\infty,+\infty)$ is invariant under scalar multiplication
$
(\pi,\la,\tau)\mapsto (\pi,c\cdot\la,c\cdot\tau)
$
for $c\in\RR_+$. Let $r_1=r_1(M)$ be the integer in the statement of Lemma \ref{lem2HypBoxes}. Then $(\pi,\la,\tau)$ is contained in some hyperbolic box $\scatola_\pi(I,\de,\nu)$ adapted to $\cI$, where the path $\de$ and $\nu$ have length bounded by $r_1$. Moreover $\scatola_\pi(I,\de,\nu)$ remains an hyperbolic box adapted to $\cI$ replacing $\de$ and $\nu$ by longer segments, thus we can assume that $\de=\ga(1,m)$ and  $\nu:=\ga(-m+1,0)$. Consider the unique $V\in\cV_\pi$ intersecting $\scatola_\pi(I,\de,\nu)$ and let $C:=V\cap\cI_\pi$ be the corresponding open subcone of $V$. Since $(\pi,\la,\tau)$ belongs to $\cI_\pi$ and to $\scatola_\pi(I,\de,\nu)$, then $(\pi,\la,\tau)\in C$. On the other hand Lemma \ref{lem3HypBoxes} implies
$
(\pi,\la^{(r)},\tau^{(r)})\in C
$.
Finally $Q_\ga$ maps
$
\{\pi\}\times\De_\ga\times\Th_\pi
$
onto
$
\{\pi\}\times\RR_+^\cA\times\Th_\ga
$,
and of course
$
(\pi^{(r)},\la^{(r)},\tau^{(r)})=
Q_\ga(\pi,\la,\tau)
$,
therefore both
$
(\De_\ga\times\Th_\pi)\cap C
$
and
$
(\RR_+^\cA\times\Th_\ga)\cap C
$
are not empty, hence Lemma \ref{lemMaingeomlemma} implies
$$
(\De_\ga\times\Th_\ga)\cap C\not=\emptyset.
$$
Moreover both $(\pi,\la,\tau)$ and $Q_\ga(\pi,\la,\tau)$ belong to $C$, therefore part (1) of Remark \ref{remInvarianceVectorSubspaces} implies
$$
Q_\ga(C)=C.
$$
The last two conditions imply that Corollary \ref{corMainGeomCor} can be applied with $\ga_1=\ga,\dots,\ga_k=\ga$ for arbitrary $k$, thus we get
$$
(\pi,\la^\ga,\tau^\ga)\in
\bigcap_{\lungh(\ga\ast\dots\ast\ga)\to\infty}
\{\pi\}\times
\De_{\ga\ast\dots\ast\ga}
\times
\Th_{\ga\ast\dots\ast\ga}
\subset\cI.
$$
\end{dimo}

\subsubsection{Shadowing Lemma}

\begin{prop}\label{prop:appendixshadowinglemma}
Fix $\pi\in\cR$ and $M>0$. There exists a positive integer $r_1=r_1(M)$ such that for any $m>r_1$ the following holds. Consider a sequence $(\ga_i)_{i\in\ZZ}$ of positive loops at $\pi$ such that $N(\ga_i)<M$ and satisfying the extra conditions:
\begin{enumerate}
\item
we have
$
(\la^{\ga_i},\tau^{\ga_i})\in\cI_\pi
$
for all $i\in\ZZ$
\item
there exists positive paths $\de$ and $\nu$, both of length $m$, such that for any $i\in\ZZ$ we have
$$
\ga_i(-m+1,0)=\nu
\textrm{ and }
\ga_i(1,m)=\de.
$$
\end{enumerate}
Then the bi-infinite path
$
\ga(-\infty,+\infty)=
\dots\ast\ga_{i}\ast\ga_{i+1}\ast\dots
$
is induced by data $(\pi,\la,\tau)$ with $(\la,\tau)\in\cI_\pi$.
\end{prop}

\begin{dimo}
Let $r_1=r_1(M)$ be the constant in the statement of Lemma \ref{lem2HypBoxes}. Since for any $i$ we have
$
(\la^{\ga_i},\tau^{\ga_i})\in\cI_\pi
$
and moreover $\ga_i$ is periodic with $N(\ga_i)<M$, then Lemma \ref{lem2HypBoxes} implies that any
$
(\pi,\la^{\ga_i},\tau^{\ga_i})
$
is contained in some hyperbolic box adapted to $\cI$, whose size is bounded by $r_1$. Assumption (2) implies that all there hyperbolic boxes contain $
\scatola_\pi(\de,\nu,I)
$
and moreover the latter contains all data
$
(\pi,\la^{\ga_i},\tau^{\ga_i})
$
for $i\in\ZZ$. It is also obvious that
$
\scatola_\pi(\de,\nu,I)
$
is adapted to $\cI$. Let $V$ be the unique linear space in $\cV_\pi$ such that the intersection
$
V\cap\scatola_\pi(\de,\nu,I)
$
is not empty and consider its convex open subcone $C:=V\cap\cI_\pi$. Lemma \ref{lem2HypBoxes} implies
$
(\la^{\ga_i},\tau^{\ga_i})\in C
$
for any $i\in\ZZ$. On one hand, it follows that $C$ intersects
$
\De_{\ga_i}\times\Th_{\ga_i}
$
for any $i\in\ZZ$. On the other hand part (2) of Remark \ref{remInvarianceVectorSubspaces} implies
$
Q_{\ga_i}(V)=V
$
for any $i\in\ZZ$. Therefore, for any $k\in\NN$, we can apply Corollary \ref{corMainGeomCor} to the family of positive loops
$
\ga_{-k+1},\dots,\ga_0,\dots,\ga_k
$
and we get
$$
C\cap
(\De_{\ga_1\ast\dots\ast\ga_k}\times
\Th_{\ga_{-k+1}\ast\dots\ast\ga_0})
\not=\emptyset.
$$
Observe that
$
\bigcap_{k\in\NN}
\De_{\ga_1\ast\dots\ast\ga_k}\times
\Th_{\ga_{-k+1}\ast\dots\ast\ga_0}
$
is the half line spanned by the data $(\pi,\la,\tau)$, because $N(\ga_i)$ is uniformly bounded. Therefore we get
$(\pi,\la,\tau)\in C$.
\end{dimo}

\section{Proof of finite-time estimates}\label{s7}

In this section we prove Proposition \ref{prop1s4ss1} and Proposition \ref{prop2s4ss1}. Some technicalities appear, mostly in the proof of Proposition \ref{prop1s4ss2cords}, due to the fact that we need to treat the combinatorics both of length data and suspension data, and the latter are less clearly understood than the former.

Fix an admissible combinatorial datum $\pi$ and recall that it consists in two bijections $\pi^t$ and $\pi^b$ from $\cA$ to $\{1,\dots,d\}$. In order to simplify the notation, we introduce the labelling $\al(1),\dots,\al(d)$ and the labelling $\be(1),\dots,\be(d)$ of the letters of $\cA$, according to their order in $\pi^t$ and $\pi^b$ respectively. With this notation, which is the same as in \S \ref{backgroundss1}, we have
\beq\label{eqlabelling}
\pi=
\begin{pmatrix}
\al(1) & \dots & \al(d) \\
\be(1) & \dots & \be(d)
\end{pmatrix}.
\eeq

Consider length-suspension data $\la$ and $\tau$ for $\pi$ and let $X$ be translation surface corresponding to $(\pi,\la,\tau)$. We fix the normalization
$$
\area(\pi,\la,\tau):=\area(X)=1.
$$

\subsection{Combinatorial control for $(\pi,\la,\tau)$}\label{s4ss2controlfinitetime}

\subsubsection{Distortion of length data and singularities of IETs}\label{s4ss2sss1controlfinitetime}

Let $\ga$ be a finite Rauzy path of length $r$ and let $T$ be an IET in $\De_\ga$. Set $T^{(r)}:=Q^r(T)$ and let $I^{(r)}$ be the sub-interval of $I$ where $T^{(r)}$ acts.

\begin{lem}\label{lem2s4ss2}
If $\ga$ is a positive path then $I^{(r)}$ does not contain any singularity of $T$ or of $T^{-1}$ in its interior.
\end{lem}

\begin{dimo}
We prove the statement for the singularities of $T$, the argument for those of $T^{-1}$ being the same. For any singularity $u^{(r),t}_\al$ of $T^{(r)}$ we have
$
u^{(r),t}_\al=T^{-k}(u^t_\al)
$,
where $k$ is the minimal positive integer such that $T^{-k}(u^t_\al)$ belongs to the interior of $I^{(r)}$.
Suppose that $u^t_\al$ belongs to the interior of $I^{(r)}$ for some $\al$ with $\pi^{t}(\al)>1$. Then $u^t_\al=u^{(r),t}_\al$, that is $u^t_\al$ is a singularity of $T^{(r)}$ too. Therefore
$
I^t_\be\cap T^{k}(I^{(r),t}_\al)=\emptyset
$
for any letter $\be$ with $\pi^t(\be)<\pi^t(\al)$ and for any $k$ with $0\leq k <q^\ga_\al$. Since $[B_\ga]_{\al\be}$ is positive by assumption then we get an absurd, according to Lemma \ref{lem3backgroundss2}.
\end{dimo}

\begin{cor}\label{cor1s4ss2}
If $\ga$ is a positive path then
$$
\min_{\chi\in\cA}
\la^{(r)}_\chi
\leq
\min_{\pi^t(\al)>1,\pi^b(\be)>1}
|u^t_\al-u^b_\be|.
$$
\end{cor}

\begin{dimo}
Consider $\be$ and $\al$ with $\pi^t(\al)>1$ and $\pi^b(\be)>1$. According to Lemma \ref{lem2s4ss2}, consider the minimum $k$ with $1\leq k\leq r$ such that both $u^t_\al$ and $u^b_\be$ do not belong to the interior of $I^{(k)}$. In this case, either $\al$ or $\be$ is the looser of the step $T^{(k-1)}\mapsto T^{(k)}$. According to the case, we get either
$
\la^{(k)}_\al\leq |u^t_\al-u^b_\be|
$
or
$
\la^{(k)}_\al\leq |u^t_\al-u^b_\be|
$.
\end{dimo}

\subsubsection{Distortion of suspension data and cords}\label{s4ss2sss2controlfinitetime}

Fix combinatorial-length-suspension data $(\pi,\la,\tau)$ and let $X$ be the corresponding translation surface. As in \S \ref{backgroundss1}, for any $\chi\in\cA$ we set $\zeta_\chi:=\la_\chi+\sqrt{-1}\cdot\tau_\chi$. A period $v\in\hol(X)$ is a \emph{cord} with respect to $(\pi,\la,\tau)$ if there exist integers $l,m$ with $1\leq l\leq m\leq d$ and such that
$$
\textrm{ either }
v=\sum_{i=l}^{m}\zeta_{\al(i)}
\textrm{ or }
v=\sum_{i=l}^{m}\zeta_{\be(i)}.
$$
When $v$ is a cord, we say that it is a \emph{top cord} if the first condition above is satisfied, and a \emph{bottom cord} if the second condition is satisfied. In both cases, we denote $\cA(v)$ the subset of $\cA$ of those letters appearing in the sum, for which we have $v=\sum_{\al\in\cA(v)}\zeta_\al$. The same translation surface admits many different combinatorial-length-suspension data, thus there are many ways to represent the same cord. When ambiguities may arise, we write $\cA_{(\pi,\la,\tau)}(v)$ instead of just $\cA(v)$.

%In this section we prove that any estimate on the area of any cord with respect to data $(\pi,\la,\tau)$ can be reduced to an estimate on the area of some diagonal (the notion is introduced in \S \ref{s3ss1sss1}) with respect to data $(\pi^{(-r)},\la^{(-r)},\tau^{(-r)})$ in the past orbit of $(\pi,\la,\tau)$, moreover the time $r$ is bounded by the norm of positive matrices of the Kontsevich-Zorich cocycle arising from $(\pi,\la,\tau)$.

\medskip

\textbf{Notation.}
Let $\ga$ be a finite Rauzy path. A path $\ga'$ is an end of $\ga$ if there exists an other path $\ga''$ such that $\ga=\ga''\ast\ga'$. Consider two pairs $(\pi,\la,\tau)$ and $(\pi',\la',\tau')$ of combinatorial-length-suspension data and a finite Rauzy path $\ga$ such that
\beq\label{eq1s4ss2data}
(\pi,\la,\tau)=
\widehat{Q}_\ga(\pi',\la',\tau').
\eeq
For $\chi\in\cA$ write also
$
\zeta'_\chi:=\la'_\chi+\sqrt{-1}\cdot\tau'_\chi
$.
Recall from \S \ref{backgroundss3} the following facts.
\begin{enumerate}
\item
If the path $\ga$ starts at $\pi'$ and ends in $\pi$, then $\widehat{Q}_\ga$ is a linear isomorphism from
$
\{\pi'\}\times\De_\ga\times\Th_{\pi'}
$
onto
$
\{\pi\}\times\RR_+^\cA\times\Th_\ga
$.
\item
Thus for data $(\pi,\la,\tau)$ in
$
\{\pi\}\times\RR_+^\cA\times\Th_\ga
$
Equation (\ref{eq1s4ss2data}) determines an unique triple of data $(\pi',\la',\tau')$ in
$
\{\pi'\}\times\De_\ga\times\Th_{\pi'}
$.
\item
Let $\ga'$ be a path starting at $\pi'$ and ending in $\pi$ and $\ga$ be a path ending at $\pi$. If $\ga'$ is an end of $\ga$ then $\De_\ga\subset\De_{\ga'}$ and $\Th_\ga\subset\Th_{\ga'}$, therefore for data $(\pi,\la,\tau)$ in
$
\{\pi\}\times\RR_+^\cA\times\Th_\ga
$
the path $\ga'$ determines an unique triple of data $(\pi',\la',\tau')$ in
$
\{\pi'\}\times\De_{\ga'}\times\Th_{\pi'}
$
via Equation (\ref{eq1s4ss2data}).
\item
Data $(\pi,\la,\tau)$ and $(\pi',\la',\tau')$ satisfying Equation (\ref{eq1s4ss2data}) represent the same translation surface $X$.
\end{enumerate}

Recall that $\ga$ is a \emph{strongly complete} if it is concatenation of $d$ complete paths. Consider the function
$
(\pi,\la,\tau)\mapsto w(\pi,\la,\tau)
$
defined in \S \ref{s3ss1(formulastrata)} (see Definition \ref{def:w}). Proposition \ref{prop1s4ss2cords} below corresponds to Corollary \ref{cor1s4ss2} in the case of negative time.

\begin{prop}\label{prop1s4ss2cords}
Consider data $(\pi,\la,\tau)$ and a strongly complete path $\ga$ such that $(\pi,\la,\tau)$ belong to
$
\{\pi\}\times\RR_+^\cA\times\Th_\ga
$.
For any cord $v=\sum_{\al\in\cA(v)}\zeta_\al$ with respect to $(\pi,\la,\tau)$ there exists some end $\ga'$ of $\ga$ such that the data $(\pi',\la',\tau')$ determined by $\ga'$ via the relation (\ref{eq1s4ss2data}) satisfy
$$
\area(v)\geq w(\pi',\la',\tau').
$$
\end{prop}

The rest of the paragraph is devoted to the proof of Proposition \ref{prop1s4ss2cords}.

\medskip

\textbf{Assumption.} We suppose that $v$ is a top cord with respect to $(\pi,\la,\tau)$, the proof for bottom cords being the same. Thus in terms of the labelling in Equation (\ref{eqlabelling}) we assume
$\cA(v)=\{\al(l),\dots,\al(m)\}$
and
$v=\zeta_{\al(l)}+\dots+\zeta_{\al(m)}$. The letters $\al(l),\dots,\al(m)$ in $\cA(v)$ have of course a different position in the top line of combinatorial data different from $\pi$. Anyway, for simplicity, in \S \ref{s4ss2sss3controlfinitetime} and \S \ref{s4ss2sss4controlfinitetime} we keep their labelling $\al(i)$ with $l\leq i\leq m$ induced by $\pi$.

\subsubsection{Some special cases}\label{s4ss2sss3controlfinitetime}

In this subsection we fix two pairs of combinatorial-length-suspension data $(\pi,\la,\tau)$ and $(\pi',\la',\tau')$ and a Rauzy path $\ga$ as in Equation (\ref{eq1s4ss2data}). In the following Lemmas we treat special cases of paths $\ga$ satisfying different assumptions. In \S \ref{s4ss2sss4controlfinitetime} we show that the general case can be reduced to these special cases.

\begin{lem}\label{lem2s4ss2specialcase1}
Let $v$ be a top cord with respect to $(\pi,\la,\tau)$ of the simplest form $v=\zeta_\al$, that is $\cA(v)=\{\al\}$ for some $\al$ in $\cA$. The following holds.
\begin{enumerate}
\item
If $\ga=\nu_{\al}$ is a simple arrow with winner $\al$, then
$
\area(v)\geq w(\pi',\la',\tau')
$,
where $(\pi',\la',\tau')$ are the data determined by $\nu_\al$ via the relation (\ref{eq1s4ss2data}).
\item
Moreover $\nu_\al$ is a top arrow if and only if $\tau_\al<0$. Similarly $\nu_\al$ is a bottom arrow if and only if $\tau_\al>0$.
\end{enumerate}
\end{lem}

\begin{dimo}
Part (2) follows trivially from \S \ref{backgroundss3}, thus we just prove part (1). Suppose that $\nu_\al$ is a top arrow, the other case being the same. Let $\be$ be the looser of $\nu_\al$. The pre-image of $\pi$ under $\nu_\al$ has the form
$
\pi'=
\begin{pmatrix}
\dots&\al\\
\dots&\be
\end{pmatrix}
$,
thus we have
$
v=
\zeta_\al=\zeta'_\al-\zeta'_\be=
\langle
\zeta',
w_{\pi',\be,\al}
\rangle
$.
\end{dimo}

A \emph{neutral} path is a finite path $\ga$ that does not contain arrows with winner $\al\in\cA(v)$, moreover the only arrows in $\ga$ where some $\al\in\cA(v)$ loses (if they exist) are top arrows. The terminology is motivated by Lemma \ref{lem3s4ss2specialcase2}, which says that the Rauzy operations contained in a neutral path does not affect $v$. Recall that with our notation any letter $\al(i)\in\cA(v)$ is in $i$-th position in the top line of $\pi$, that is $\pi^t(\al(i))=i$ for any $l\leq i\leq m$.

\begin{lem}\label{lem3s4ss2specialcase2}
Let $\ga$ be a neutral path. Then the following holds.
\begin{enumerate}
\item
There exist some integer $k$ with $0\leq k< m-l$ such that $\al(i)$ is in $(i-k)$-th position in the top line of $\pi'$ for any $l\leq i\leq m$.
\item
In particular, if $\al(m)$ is in $d$-th position in the top line of $\pi'$, then $\ga$ is trivial.
\item
Finally, $v$ is a top cord with respect to $(\pi',\la',\tau')$ with
$
\cA_{(\pi',\la',\tau')}(v)=
\cA_{(\pi,\la,\tau)}(v)
$.
\end{enumerate}
\end{lem}

\begin{dimo}
We first prove part (1) of the Lemma. We suppose that $\ga$ is a simple arrow, the general case follows repeating the argument finitely many times. Let $\be$ be the winner of $\ga$, which does not belong to $\cA(v)$. If $\ga$ is a top arrow, then the statement obviously follows with $k=0$. If $\ga$ is a bottom arrow, then we have two possibilities. If $\pi^t(\be)>m$ then the statement follows with $k=0$. Otherwise we must have $\pi^t(\be)<l$, because $\be$ does not belong to $\cA(v)$. Moreover $\al(l)$ cannot be the loser of $\ga$ by assumption, thus $\pi^t(\be)\leq l-2$. In this case the statement follows with $k=1$. In order to prove part (2), observe that if $\al(m)$ is in $d$-th position in the top line of $\pi'$, then part (1) implies that $\al(m)$ is in $d$-th position also in the top line of $\pi$, and in such position $\al(m)$ can't neither win nor loose, since $\ga$ is neutral. Part (3) follows from part (1), indeed we have $
\zeta'_\al=\zeta_\al
$
for any $\al$ in $\cA(v)$, since $\ga$ does not contain arrows with winner in $\cA(v)$.
\end{dimo}

\begin{lem}\label{lem4s4ss2specialcase3}
Consider a path $\ga$ of the form
$
\ga=
\eta_{\al(m)}\ast\dots\ast\eta_{\al(l)}
\ast\ga^{(1)}
$,
where $\ga^{(1)}$ is a (possibly trivial) neutral path and
$
\eta_{\al(m)},\dots,\eta_{\al(l)}
$
are bottom arrows with losers respectively $\al(m),\dots,\al(l)$ and with the same winner $\be\in\cA\setminus\cA(v)$. Then any letter $\al(i)$ is in $(d-m+i)$-th position in the top line of $\pi'$ for any $l\leq i\leq m$, that is
$$
\pi'=
\begin{pmatrix}
\dots&\be&\dots&\al(l)&\dots&\al(m)\\
     &   &\dots&      &     &\be
\end{pmatrix}.
$$
Moreover $v$ is a top cord with respect to $(\pi',\la',\tau')$ with
$
\cA_{(\pi',\la',\tau')}(v)=
\cA_{(\pi,\la,\tau)}(v)
$.
\end{lem}

\begin{dimo}
Let $(\pi^{(1)},\la^{(1)},\tau^{(1)})$ be the data determined by $\ga^{(1)}$ via the relation (\ref{eq1s4ss2data}). Lemma \ref{lem3s4ss2specialcase2} implies that the letters $\al(l),\dots,\al(m)$ keep their reciprocal order in the top line of $\pi^{(1)}$ and they are possibly translated to the left. Moreover $\pi^{(1)}$ is the endpoint of the bottom arrow $\eta_{\al(l)}$ with winner $\be$ and looser $\al(l)$, therefore we have
$
\pi^{(1)}=
\begin{pmatrix}
\dots&\be&\al(l)&\dots&\al(m)&\dots\\
     &   & \dots&     &      &\be
\end{pmatrix}
$.
The pre-image of $\pi^{(1)}$ under the arrows
$
\eta_{\al(m)},\dots,\eta_{\al(l)}
$
is easily computable, hence the statement on $\pi'$ follows. Since $\pi'$ is as in the statement, then the sum $
\zeta'_{\al(l)}+\dots+\zeta'_{\al(m)}
$
is a top cord with respect to $(\pi',\la',\tau')$. Moreover $\ga$ does not contain arrows with winner in $\cA(v)$ therefore $\zeta'_\al=\zeta_\al$ for any $\al$ in $\cA(v)$. It follows that
$
v=
\zeta_{\al(l)}+\dots+\zeta_{\al(m)}=
\zeta'_{\al(l)}+\dots+\zeta'_{\al(m)}
$,
that is $v$ is a top cord with respect to $(\pi',\la',\tau')$.
\end{dimo}

The Corollary below follows immediately repeating $k$ times the argument in Lemma \ref{lem4s4ss2specialcase3}.

\begin{cor}\label{cor2s4ss2specialcase4}
Consider a path $\ga$ of the form $\ga=\ga_k\ast\dots\ast\ga_1$, where $\ga_1,\dots,\ga_k$ are as in Lemma \ref{lem4s4ss2specialcase3}. Then $\pi'$ and $v$ are as in the conclusion of Lemma \ref{lem4s4ss2specialcase3}.
\end{cor}

\begin{lem}\label{lem5s4ss2specialcase5}
Fix an integer $i$ such that $l\leq i\leq m$ and consider a path $\ga$ given by the concatenation
$
\ga=
\nu_{\al(i)}\ast
\eta_{\al(i)}\ast\dots\ast\eta_{\al(l)}\ast
\ga^{(2)}\ast
\ga^{(1)}
$,
where
\begin{enumerate}
\item
$\ga^{(1)}$ is a (possibly trivial) path as in Corollary \ref{cor2s4ss2specialcase4},
\item
$\ga^{(2)}$ is a (possibly trivial) neutral path,
\item
$\eta_{\al(i)},\dots,\eta_{\al(l)}$ are bottom arrows with looser respectively $\al(i),\dots,\al(l)$ and winner $\be$ in $\cA\setminus\cA(v)$,
\item
$\nu_{\al(i)}$ is a top arrow with winner $\al(i)$.
\end{enumerate}
Then, denoting by $(\pi',\la',\tau')$ the data determined by $\ga$ via (\ref{eq1s4ss2data}) we have
$
\area(v)\geq w(\pi',\la',\tau').
$
\end{lem}

\begin{dimo}
Denote
$
(\pi^{(1)},\la^{(1)},\tau^{(1)})
$
and
$
(\pi^{(2)},\la^{(2)},\tau^{(2)})
$
the triples of data determined via (\ref{eq1s4ss2data}) respectively by the paths $\ga^{(1)}$ and $\ga^{(2)}\ast\ga^{(1)}$. According to Corollary \ref{cor2s4ss2specialcase4} applied to the path $\ga^{(1)}$, for $\pi^{(1)}$ we have the same conclusion as for $\pi'$ in Lemma \ref{lem4s4ss2specialcase3}. Then Lemma \ref{lem3s4ss2specialcase2} applies to the path $\ga^{(2)}$ and moreover $\pi^{(2)}$ is the ending point of the bottom arrow $\eta_{\al(l)}$ with winner $\be\in\cA\setminus\cA(v)$, thus we have
$
\pi^{(2)}=
\begin{pmatrix}
\dots&\be&\al(l)&\dots&\al(m)&\dots\\
     &   &\dots &     &      &\be
\end{pmatrix}.
$
Observe that $\pi^{(2),b}(\al(i))<d$, since $\be\not\in\cA(v)$. Denote $\chi$ the letter in $\cA$ with
$
\pi^{(2),b}(\chi)=\pi^{(2),b}(\al(i))+1
$, where the case $\chi=\be$ is possible. It is easy to check that the pre-image of $\pi^{(2)}$ under the concatenation of arrows
$
\nu_{\al(i)}\ast
\eta_{\al(i)}\ast\dots\ast\eta_{\al(l)}
$
has the form
$$
\pi'=
\begin{pmatrix}
\dots&\be   & \dots&\al(l)&\dots&\al(i)\\
\dots&\al(i)&      &\dots & \be &  \chi
\end{pmatrix},
$$
where for any $l\leq j\leq i$ the letter $\al(j)$ is in $(d+l-j)$-th position in the top line of $\pi'$.
%We remark that in case $\chi=\be$ the letters $\al(i)$ and $\be$ are in $(d-1)$-th position and $d$-th position respectively in the bottom line of $\pi'$.
Since $\ga^{(1)}$ and $\ga^{(2)}$ do not contain arrows with winner in $\cA(v)$ then $\zeta^{(2)}_\al=\zeta_\al$ for any $\al$ in $\cA(v)$. Moreover the winner of the arrows $\eta_{\al(i)},\dots,\eta_{\al(l)}$ is $\be\in\cA\setminus\cA(v)$ and $\nu_{\al(i)}$ has winner $\al(i)$ and looser $\chi$, thus we have
$
\zeta^{(2)}_{\al(i)}=
\zeta'_{\al(i)}-\zeta'_\chi
$
and
$
\zeta'_\al=\zeta^{(2)}_\al
$
for any $\al$ in $\cA(v)\setminus\{\al(i)\}$. Setting
$
v':=\zeta_{\al(l)}+\dots+\zeta_{\al(i)}
$
we have
$$
v'=
\zeta_{\al(l)}+\dots+\zeta_{\al(i)}=
\zeta^{(2)}_{\al(l)}+\dots+\zeta^{(2)}_{\al(i)}=
\zeta'_{\al(l)}+\dots+\zeta'_{\al(i)}-\zeta'_\chi=
\langle\zeta',w_{\pi',\al(i),\chi}\rangle.
$$
Remark that in general $v'$ is not a period of the translation surface $X$, and thus nor a top cord with respect to $(\pi,\la,\tau)$. In order to prove that $\area(v')\leq\area(v)$, observe first that $\Im(v')<0$. Consider the data $(\pi,\la,\tau)$ and the corresponding broken line in the complex plane whose vertices are $\xi^t_{\al(j)}$ with $l\leq j\leq m$ and $\xi^t_{\al(m)}+\zeta_{\al(m)}$. The vertices $\xi^t_{\al(l)}$ and $\xi^t_{\al(i)}+\zeta_{\al(i)}$ are connected by  $v'=\zeta_{\al(l)}+\dots+\zeta_{\al(i)}$, whereas the period $v=\zeta_{\al(l)}+\dots+\zeta_{\al(m)}$ connects the extremal vertices $\xi^t_{\al(l)}$ and $\xi^t_{\al(m)}+\zeta_{\al(m)}$. Since $v$ is a top cord then the intermediate vertex $\xi^t_{\al(i)}+\zeta_{\al(i)}$ lies above the line trough
$\xi^t_{\al(l)}$ and $\xi^t_{\al(m)}+\zeta_{\al(m)}$, moreover such line has negative slope, because $\Im(v')<0$. It follows that
$
\Im(\xi^t_{\al(l)})
\geq
\Im(\xi^t_{\al(i)}+\zeta_{\al(i)})
\geq
\Im(\xi^t_{\al(m)}+\zeta_{\al(m)})
$
and therefore $|\Im(v')|\leq|\Im(v)|$. We have obviously
$
\Re(v')=
\la_{\al(l)}+\dots+\la_{\al(i)}\leq
\la_{\al(l)}+\dots+\la_{\al(m)}=\Re(v)
$,
thus $\area(v')\leq\area(v)$.
\end{dimo}

\begin{lem}\label{lem6s4ss2specialcase6}
Consider a path $\ga$ of the form $\ga=\nu_\al$ where $\nu_\al$ is a bottom arrow with winner $\al\in\cA(v)$ and let $(\pi',\la',\tau')$ be the data determined by $\nu_\al$ via (\ref{eq1s4ss2data}). We have the following two possibilities
\begin{enumerate}
\item
If $\al=\al(i)$ with $l\leq i<m$ then $v$ is a top cord with respect to $(\pi',\la',\tau')$ and moreover we have
$
\cA_{(\pi',\la',\tau')}(v)=
\cA_{(\pi,\la,\tau)}(v)\setminus\{\al(i+1)\}
$.
\item
If $\al=\al(m)$ then
$
\area(v)\geq w(\pi',\la',\tau')
$.
\end{enumerate}
\end{lem}

\begin{dimo}
The data $(\pi,\la,\tau)$ are the image of a bottom arrow with winner $\al$, thus $\al$ is in $d$-th position in the bottom line of $\pi$, moreover part (2) of Lemma \ref{lem2s4ss2specialcase1} implies $\tau_\al>0$ both in case (1) and in case (2).

If $\al=\al(i)$ with $l\leq i<m$ then $\al(i+1)$ belongs to $\cA(v)$ and the pre-image of $\pi$ under $\nu_{\al(i)}$ is
$$
\pi'=
\begin{pmatrix}
\dots&\al(l)&\dots&\al(i)&\al(i+2)&\dots&\al(m)&\dots&\al(i+1)\\
     &      &     &\dots &        &    &      &    &\al(i)
\end{pmatrix}.
$$
In this case we have
$
\zeta'_{\al(i)}=
\zeta_{\al(i)}+\zeta_{\al(i+1)}
$
and therefore
$
v=
\sum_{\al\in\cA(v)}\zeta_\al=
\sum_{\al\in\cA(v)\setminus\{\al(i+1)\}}\zeta'_\al
$.
It follows that $v$ is a top cord with respect to the data $(\pi',\la',\tau')$ with
$
\cA_{(\pi',\la',\tau')}(v)=
\cA_{(\pi,\la,\tau)}(v)\setminus\{\al(i+1)\}
$.

If $\al=\al(m)$ then we must have $\pi^t(\al(m))<d$, because $\pi^b(\al(m))=d$. If $\be\in\cA\setminus\cA(v)$ is the letter such that $\pi^t(\be)=\pi^t(\al(m))+1$, that is the letter which looses against $\al(m)$ in $\nu_{\al(m)}$, then the pre-image of $\pi$ under $\nu_{\al(m)}$ is
$$
\pi'=
\begin{pmatrix}
\dots&\al(l)&\dots&\al(m)&\dots&\be\\
     &      &\dots&      &     &\al(m)
\end{pmatrix},
$$
where the position of any letter $\al(i)\in\cA(v)$ in the top line of $\pi'$ is the same as its position in the top line of $\pi$. We have
$
\zeta_{\al(m)}=
\zeta'_{\al(m)}-\zeta'_\be=
\langle
\zeta',w_{\pi',\be,\al(m)}
\rangle
$
and it is enough to prove that $\area(\zeta_{\al(m)})\leq\area(v)$. Consider the data $(\pi,\la,\tau)$ and the corresponding broken line whose vertices are $\xi^t_{\al(j)}$ with $l\leq j\leq m$ and $\xi^t_{\al(m)}+\zeta_{\al(m)}$. The vertices $\xi^t_{\al(m)}$ and $\xi^t_{\al(m)}+\zeta_{\al(m)}$ are connected by the period $\zeta_{\al(m)}$, whereas the period $v=\zeta_{\al(l)}+\dots+\zeta_{\al(m)}$ connects the extremal vertices $\xi^t_{\al(l)}$ and $\xi^t_{\al(m)}+\zeta_{\al(m)}$. Since $v$ is a top cord then the intermediate vertex $\xi^t_{\al(m)}$ lies above the line trough
$\xi^t_{\al(l)}$ and $\xi^t_{\al(m)}+\zeta_{\al(m)}$, moreover such line has positive slope, because $\tau_{\al(m)}>0$. It follows that
$
\Im(\xi^t_{\al(l)})
\leq
\Im(\xi^t_{\al(m)})
\leq
\Im(\xi^t_{\al(m)}+\zeta_{\al(m)})
$
and therefore $|\Im(\zeta_{\al(m)})|\leq|\Im(v)|$. We have obviously
$
\Re(\zeta_{\al(m)})=
\la_{\al(m)}\leq
\la_{\al(l)}+\dots+\la_{\al(m)}=\Re(v)
$,
thus $\area(\zeta_{\al(m)})\leq\area(v)$.
\end{dimo}

\subsubsection{The general case: proof of Proposition \ref{prop1s4ss2cords}}\label{s4ss2sss4controlfinitetime}

Let $\ga$ be as in Proposition \ref{prop1s4ss2cords} and decompose it as
$$
\ga=\ga^{activ}\ast\ga^{neutral},
$$
where $\ga^{neutral}$ is the maximal neutral end of $\ga$ (recall that a finite Rauzy path is said neutral if it does not contain neither arrows with winner $\al\in\cA(v)$ not bottom arrows with looser $\al\in\cA(v)$, whereas top arrows with looser $\al\in\cA(v)$ are admitted). Let $(\pi^{(1)},\la^{(1)},\tau^{(1)})$ be the combinatorial-length-suspension data determined by $\ga^{neutral}$ via the relation (\ref{eq1s4ss2data}), thus in particular $\ga^{activ}$ ends in $\pi^{(1)}$. Since $\ga^{neutral}$ does not contain arrows with winner in $\cA(v)$, then $\ga^{activ}$ is the concatenation of $d$ paths $\ga_1\ast\dots\ast\ga_d$, where each letter in $\cA(v)$ wins at least once in each $\ga_i$. Moreover, the last arrow $\ga_{last}$ of $\ga^{activ}$ is as in cases A or B or C as below.

\begin{description}
\item[A]
We have $\ga_{last}=\nu_\al$, where $\nu_\al$ is a top arrow with winner $\al\in\cA(v)$. In this case part (2) of Lemma \ref{lem3s4ss2specialcase2} implies that $\al=\al(m)$ and $\ga^{neutral}$ is the trivial path. Therefore we have
$$
\pi^{(1)}=
\pi=
\begin{pmatrix}
\dots&\al(l)&\dots&\al(m)\\
\dots&\al(m)&\dots&
\end{pmatrix}.
$$
\item[B]
We have $\ga_{last}=\eta_\al$, where $\eta_\al$ is a bottom arrow with looser $\al\in\cA(v)$ and whose winner $\be$ does not belong to $\cA(v)$. In this case, part (1) of Lemma \ref{lem3s4ss2specialcase2} implies
$
\pi^{(1),t}(\be)=\pi^{(1),t}(\al(l))-1
$,
thus the loser of $\eta_\al$ is $\al=\al(l)$ and $\pi^{(1)}$ has the form
$$
\pi^{(1)}=
\begin{pmatrix}
\dots&\be&\al(l)&\al(l+1)&\dots&\al(m)&\dots\\
     &   &      & \dots  &     &      & \be
\end{pmatrix}.
$$
\item[C]
We have $\ga_{last}=\nu_\al$, where $\nu_\al$ is a bottom arrow with winner in $\cA(v)$. In this case we consider the following two subcases.
\begin{description}
\item[C1]
We have $\al=\al(m)$.
\item[C2]
We have $\al=\al(i)$ with $l\leq i<m$.
\end{description}
\end{description}

\begin{rem}\label{rem1appendix1ss1}
According to Lemma \ref{lem3s4ss2specialcase2}, $v$ is a top cord with respect to $(\pi^{(1)},\la^{(1)},\tau^{(1)})$ and we have
$
\cA_{(\pi^{(1)},\la^{(1)},\tau^{(1)})}(v)=
\cA_{(\pi,\la,\tau)}(v)
$.
Therefore, in order to prove Proposition \ref{prop1s4ss2cords} it is enough to replace $\ga$ by $\ga^{activ}$ and $(\pi,\la,\tau)$ by $(\pi^{(1)},\la^{(1)},\tau^{(1)})$. This amounts to consider data $(\pi,\la,\tau)$ in
$
\{\pi\}\times\RR_+^\cA\times\Th_\ga
$,
where $\ga$ is a path satisfying the conditions below.
\begin{enumerate}
\item
The path $\ga$ is the concatenation of $d$ paths $\ga_1\ast\dots\ast\ga_d$, where each letter in $\cA(v)$ wins in each $\ga_i$.
\item
The last arrow $\ga_{last}$ of $\ga$ is as in cases A or B or C above.
\end{enumerate}
\end{rem}

Case A. Let $(\pi',\la',\tau')$ be the data determined by $\nu_{\al(m)}$ via Equation (\ref{eq1s4ss2data}). If $\be$ is the letter such that $\pi^b(\be)=\pi^b(\al(m))+1$, we have
$
\zeta_{\al(m)}=\zeta'_{\al(m)}-\zeta'_{\be}
$.
Since the pre-image of $\pi$ under $\nu_{\al(m)}$ is
$
\pi'=
\begin{pmatrix}
\dots&\al(l)&\dots&\al(m)\\
\dots&\al(m)&\dots&\be
\end{pmatrix},
$
then Proposition \ref{prop1s4ss2cords} follows with $\ga'=\nu_{\al(m)}$, indeed we have
$$
v=
\zeta_{\al(l)}+\dots+\zeta_{\al(m)}=
\zeta'_{\al(l)}+\dots+\zeta'_{\al(m)}-\zeta'_{\be}=
\langle
\zeta',w_{\pi',\be,\al(m)}
\rangle.
$$

Case B. Decompose the path $\ga$ as $\ga=\ga^{start}\ast\ga^{end}$, where $\ga^{end}$ is minimal containing an arrow $\nu_\al$ with winner in $\cA(v)$. Observe that such decomposition exists, because each letter in $\cA(v)$ wins at least once in $\ga$. We consider separately two cases, according to the type (top or bottom) of $\nu_\al$.

If $\nu_\al$ is a top arrow, let $i$ be the integer with $l\leq i\leq m$ such that the winner of $\nu_\al$ is $\al=\al(i)$. It is not difficult to check that $\ga^{end}$ decomposes as in Lemma \ref{lem5s4ss2specialcase5}. Thus, according to the same Lemma, Proposition \ref{prop1s4ss2cords} follows with $\ga':=\ga^{end}$.

If $\nu_\al$ is a bottom arrow, it is easy to check that we have a decomposition
$
\ga^{end}=\nu_\al\ast\ga^{(1)}\ast\ga^{(2)}
$,
where $\ga^{(1)}$ is neutral and $\ga^{(2)}$ is as in Corollary \ref{cor2s4ss2specialcase4}. In this case, if $(\pi'',\la'',\tau'')$ are the data determined by $\ga^{(1)}\ast\ga^{(2)}$ via (\ref{eq1s4ss2data}), then Lemma \ref{lem3s4ss2specialcase2} and Corollary \ref{cor2s4ss2specialcase4} imply that $v$ is a top cord with respect to $(\pi'',\la'',\tau'')$ with
$
\cA_{(\pi'',\la'',\tau'')}(v)=
\cA_{(\pi,\la,\tau)}(v)
$.
Moreover the path $\ga^{start}\ast\nu_\al$ is a in Remark \ref{rem1appendix1ss1}, where its last arrow is as in case C, therefore in this case we can replace respectively $\ga$ by $\ga^{start}\ast\nu_\al$ and $(\pi,\la,\tau)$ by $(\pi'',\la'',\tau'')$ and then apply the discussion of Case C.

\medskip

Case C1. Let $(\pi',\la',\tau')$ be the data determined by $\nu_{\al(m)}$ via the relation (\ref{eq1s4ss2data}). Proposition \ref{prop1s4ss2cords} follows with $\ga:=\nu_{\al(m)}$, $\cA(\pi,\la',\tau'):=\{\al(m)\}$ and $v':=\zeta_{\al(m)}$, according to part (2) of Lemma \ref{lem6s4ss2specialcase6}.

\medskip

Case C2. Decompose $\ga$ as
$
\ga=\ga^{start}\ast\nu_{\al(i)}
$
and observe that $\ga^{start}$ is the concatenation of at least $d-1$ complete paths, because the initial path $\ga$ in Proposition \ref{prop1s4ss2cords} is assumed strongly complete. Part (1) of Lemma \ref{lem6s4ss2specialcase6} applies, therefore $v$ is a top cord with respect to the data $(\pi',\la',\tau')$ determined by $\nu_{\al(i)}$ via the relation (\ref{eq1s4ss2data}) and moreover we have
$
\cA_{(\pi',\la',\tau')}(v)=
\cA_{(\pi,\la,\tau)}(v)\setminus\{\al(i+1)\}
$.
We start an iterative procedure on the cardinality of $\cA(v)$, replacing the data $(\pi,\la,\tau)$ by $(\pi',\la',\tau')$, the path $\ga$ by $\ga^{start}$ and the set $\cA(v)$ by $\cA(v)\setminus\{\al(i+1)\}$ and then repeating the discussion of \S \ref{s4ss2sss4controlfinitetime} for these new data. At each step, either Proposition \ref{prop1s4ss2cords} is proved, or we end up in Case C2 and we apply part (1) of Lemma \ref{lem6s4ss2specialcase6}, so that the new set $\cA(v)$ looses an element. After at most $m-l-1$ steps, if Proposition \ref{prop1s4ss2cords} is not yet proved, then we end up with a top cord of the simple form $v=\zeta_\al$ for some $\al$ in $\cA$. In this case, Proposition \ref{prop1s4ss2cords} follows from part (1) of Lemma \ref{lem2s4ss2specialcase1}. The proof of Proposition \ref{prop1s4ss2cords} is completed.

\subsection{Geometric bounds for $(\pi,\la,\tau)$}\label{s4ss3}

This section is devoted to the proof of Proposition \ref{prop1s4ss1} and Proposition \ref{prop2s4ss1}. Fix combinatorial-length-suspension data $(\pi,\la,\tau)$ such that $\area(\pi,\la,\tau)=1$. In this section we give an estimate of the distortion $\De(T)$ of $T=(\pi,\la)$ in terms of the minimum of the values $\area(v)$ over all cords $v$ with respect to $(\pi,\la,\tau)$.

\subsubsection{Geometric bounds via the control on cords}

Let $k$ be an integer with $1\leq k<d$. Recall from \S \ref{backgroundss1} that since $\tau$ is a suspension datum for $\pi$ then we have
\beq\label{eq2s4ss3controloncords}
h_{\al(k+1)}\geq \tau_{\al(1)}+\dots+\tau_{\al(k)}>0
\textrm{ and }
-h_{\al(k+1)}\leq \tau_{\be(1)}+\dots+\tau_{\be(k)}<0.
\eeq
If we have
$
h_{\al(k)}<
\tau_{\al(1)}+\dots+\tau_{\al(k)}
$
we say that the letter $\al(k)$ is \emph{special in top line}. Observe that the condition is equivalent to have
$
\pi^b\big(\al(k)\big)=d
$,
that is $\al(k)=\be(d)$, and $\sum_\chi\tau_\chi>0$. Similarly, if we have
$
-h_{\al(k)}>
\tau_{\be(1)}+\dots+\tau_{\be(k)}
$
we say that the letter $\be(k)$ is \emph{special in bottom line} and we observe that this second condition is equivalent to have
$
\pi^t\big(\be(k)\big)=d
$,
that is $\be(k)=\al(d)$, and $\sum_\chi\tau_\chi<0$. We remark that for any triple of data $(\pi,\la,\tau)$ with $\sum_\chi\tau_\chi\not=0$ there always exists either an $\al(k)$ which is special in top line, or an $\be(k)$ which is special in bottom line. Moreover, a special letter in top line is necessarily unique, and it excludes the existence of special letters in bottom line. Similarly, a special letter in bottom line is unique, and it excludes the existence of special letters in top line.

\begin{lem}\label{lem1s4ss3controloncords1}
Consider an integer $k$ with $1\leq k<d$. The following holds.
\begin{enumerate}
\item
We always have
$
0<
\tau_{\al(1)}+\dots+\tau_{\al(k)}<
(\la_{\al(k+1)})^{-1}
$.
Moreover if $\al(k)$ is not special in top line then we have
$
0<
\tau_{\al(1)}+\dots+\tau_{\al(k)}<
(\la_{\al(k)})^{-1}
$.
\item
We always have
$
-(\la_{\be(k+1)})^{-1}<
\tau_{\be(1)}+\dots+\tau_{\be(k)}<0
$.
Moreover if $\be(k)$ is not special in top line then we have
$
-(\la_{\be(k)})^{-1}<
\tau_{\be(1)}+\dots+\tau_{\be(k)}<0
$.
\end{enumerate}
\end{lem}

\begin{dimo}
We just prove part (1), part (2) being the same. We have $\la_\chi\cdot h_\chi<1$ for any $\chi\in\cA$, thus the first estimate follows from Equation (\ref{eq2s4ss3controloncords}). Since $\al(k)$ is not special in top line, then
$
\sum_{j\leq k}\tau_{\al(j)}\leq h_{\al(k)}
$
and the second inequality follows trivially.
\end{dimo}

\begin{lem}\label{lem2s4ss3controloncords2}
Consider $\de>0$ such that any cord $v$ with respect to $(\pi,\la,\tau)$ satisfies $\area(v)>\de$. For any $1\leq k<d$ we have the following two conditions.
\begin{eqnarray*}
&&
(\la_{\al(1)}+\dots+\la_{\al(k)})\cdot
|\tau_{\al(1)}+\dots+\tau_{\al(k)}|\geq \de,\\
&&
(\la_{\be(1)}+\dots+\la_{\be(k)})\cdot
|\tau_{\be(1)}+\dots+\tau_{\be(k)}|\geq \de.
\end{eqnarray*}
\end{lem}

\begin{dimo}
We prove the first equation, the second being the same. Set
$
v_k:=\zeta_{\al(1)}+\dots+\zeta_{\al(k)}
$
for $1\leq k<d$ and observe that the $v_k$ for which $\area(v_k)$ is minimal is a top cord with respect to $(\pi,\la,\tau)$.
\end{dimo}

\begin{lem}\label{lem3s4ss3controloncords3}
Consider $\de>0$ such that any cord $v$ with respect to $(\pi,\la,\tau)$ satisfies $\area(v)>\de$. For any integer $k$ with $1\leq k<d$ there exists positive integers $j$ and $i$ with $j\leq k$ and $i\leq k$ such that
$$
\la_{\al(j)}>
\frac{\de}{k}\la_{\al(k+1)}
\textrm{ and }
\la_{\be(i)}>
\frac{\de}{k}\la_{\be(k+1)}
$$
\end{lem}

\begin{dimo}
Part (1) of \ref{lem1s4ss3controloncords1} and part (1) of Lemma \ref{lem2s4ss3controloncords2} imply
$$
\frac{1}{\la_{\al(k+1)}}>
\tau_{\al(1)}+\dots+\tau_{\al(k)}\geq
\frac{\de}{\la_{\al(1)}+\dots+\la_{\al(k)}},
$$
thus there exists $j\leq k$ with
$
\la_{\al(j)}>\de\la_{\al(k+1)}/k
$.
The existence of $\be(i)$ follows by a similar argument.
\end{dimo}

\begin{lem}\label{lem4s4ss3controloncords4}
Consider $\de>0$ such that any cord $v$ with respect to $(\pi,\la,\tau)$ satisfies $\area(v)>\de$. Then for any integers $m,l$ with $2\leq l+1<m\leq d$ both the following two properties hold.
\begin{enumerate}
\item
If $\al(l)$ is not special in top line, then there exists an integer $k$ with $l<k<m$ such that
$$
\la_{\al(k)}>\frac{\de}{m-l-1}\cdot
\min\{\la_{\al(l)},\la_{\al(m)}\}.
$$
\item
If $\be(l)$ is not special in bottom line, then there exists an integer $k$ with $l<k<m$ such that
$$
\la_{\be(k)}>\frac{\de}{m-l-1}\cdot
\min\{\la_{\be(l)},\la_{\be(m)}\}.
$$
\end{enumerate}
\end{lem}

\begin{dimo}
We just prove part (1), the second being the same. Fix $\epsilon>0$ such that $\la_{\be(l)}>\epsilon$ and $\la_{\be(m)}>\epsilon$. It is enough to show that
$$
\la_{\al(l+1)}+\dots+\la_{\al(m-1)}>
\de\cdot\epsilon.
$$
Setting
$
v:=\zeta_{\al(l+1)}+\dots+\zeta_{\al(m-1)}
$,
the statement is equivalent to $\Re(v)> \de\cdot\epsilon$. Assume by absurd that $0<\Re(v)\leq \de\cdot \epsilon$. We have $\la_{\al(l)}\cdot h_{\al(l)}<1$ and $\la_{\al(m)}\cdot h_{\al(m)}<1$, thus part (1) of Lemma \ref{lem1s4ss3controloncords1} implies
$
\tau_{\al(1)}+\dots+\tau_{\al(m-1)}<
1/\epsilon
$
and
$
\tau_{\al(1)}+\dots+\tau_{\al(l)}<
1/\epsilon
$,
where the second inequality holds because we assume that $\al(l)$ is not special in top line. Consider the broken line whose vertices are the points
$
\xi^t_{\al(k)}=
\zeta_{\al(1)}+\dots+\zeta_{\al(k-1)}
$
for $l<k\leq m$. If we assume $0<\Re(v)\leq \de\cdot \epsilon$, then the last two inequalities imply that any vertex $\xi^t_{\al(k)}$ with $l+1<k\leq m-1$ must lie above the line connecting the extremal vertices $\xi^t_{\al(l+1)}$ and $\xi^t_{\al(m)}$ indeed otherwise we would get a top cord $v'$ with
$
\area(v')<
\de\cdot \epsilon\cdot1/\epsilon=\de
$.
It follows that $v$ is a top cord, but then the assumption $0<\Re(v)\leq \de\cdot \epsilon$ implies $\area(v)<\de$, which is absurd.
\end{dimo}

Recall the notation
$
\|\la\|:=\sum_{\chi\in\cA}\la_\chi=1
$.

\begin{prop}\label{prop1s4ss3controloncords5}
Consider $\de>0$ such that any cord $v$ with respect to $(\pi,\la,\tau)$ satisfies $\area(v)>\de$. Then for any letter $\al$ in $\cA$ we have
$$
\la_\al>
\frac{\de^{d-1}}{d!}\cdot\|\la\|.
$$
\end{prop}

\begin{dimo}
There obviously exists a letter $\chi(1)\in\cA$ such that
$
\la_{\chi(1)}\geq \|\la\|/d
$.
We set $\cA_1:=\{\chi(1)\}$ and we define an increasing sequence of sub-alphabets
$
\cA_1\subset\dots
\subset\cA_i\subset\dots
\subset\cA_d:=\cA
$
such that for any $1\leq i\leq d$ and any $\chi\in\cA_i$ we have
$$
\la_\chi>
\de^{i-1}\cdot
\frac{(d-i)!}{d!}
\cdot\|\la\|.
$$
Consider $1\leq i<d$ and suppose that $\cA_i$ is defined. In order to define $\cA_{i+1}$ it is enough to find a letter
$
\chi(i+1)\in\cA\setminus\cA_i
$
such that
\beq\label{eqconditionBprop1s4ss3controloncords5}
\la_{\chi(i+1)}>
\frac{\de}{d-i}\cdot
\min_{\chi\in\cA_i}\la_\chi.
\eeq
The sequence of alphabets $\cA_i$ is defined according to two combinatorial procedures, the discriminant to pass from the first procedure to the second being the condition
\beq\label{eqconditionAprop1s4ss3controloncords5}
\{\al(1),\be(1)\}\subset \cA_i.
\eeq

Observe that $\cA_1$ obviously does not satisfy (\ref{eqconditionAprop1s4ss3controloncords5}). Fix $1\leq i<d$, suppose that $\cA_i$ is defined and moreover that the condition (\ref{eqconditionAprop1s4ss3controloncords5}) is not satisfied. Assume $\al(1)\not\in\cA_i$, the argument for the case $\be(1)\not\in\cA_i$ being the same. Let $\al(k)\in\cA_i$ be the letter such that
$
\pi^t(\chi)\geq \pi^t(\al(k))
$
for any other $\chi\in\cA_i$. According to Lemma \ref{lem3s4ss3controloncords3} there exists $j$ with $1\leq j<k$ and such that
$
\la_{\al(j)}>\de\la_{\al(k)}(k-1)^{-1}
$.
Observe that $k-1\leq d-i$, thus the letter
$
\chi(i+1):=\al(j)\in\cA\setminus\cA_i
$
satisfies the condition (\ref{eqconditionBprop1s4ss3controloncords5}).

Let $r$ be minimal such that $\cA_r$ is defined and moreover the condition (\ref{eqconditionAprop1s4ss3controloncords5}) is satisfied. Assume that $\cA_i$ is defined for $r\leq i<d$ and observe that such $\cA_i$ still satisfies (\ref{eqconditionAprop1s4ss3controloncords5}). Consider the conditions below.
\begin{enumerate}
\item
There exist integers $l,m$ with $2\leq l+1<m\leq d$ such that $\al(l)\in\cA_i$ and $\al(m)\in\cA_i$ and
$\al(j)\in\cA\setminus\cA_i$ for any $l<j<m$.
\item
There exist integers $l,m$ with $2\leq l+1<m\leq d$ such that $\be(l)\in\cA_i$ and $\be(m)\in\cA_i$ and
$\be(j)\in\cA\setminus\cA_i$ for any $l<j<m$.
\item
Condition (1) holds,
moreover $\al(l)$ is not special in top line.
\item
Condition (2) holds, moreover $\be(l)$ is not special in bottom line.
\end{enumerate}
We first prove that we have either condition (3) or condition (4), then we will define the required letter $\chi(i+1)$ using one of these two conditions. Since $\pi$ is admissible, then condition (\ref{eqconditionAprop1s4ss3controloncords5}) implies that either condition (1) or condition (2) hold. Suppose that condition (1) holds, but (3) is no true, that is $\al(l)$ is special in top line, thus in particular $\al(l)=\be(d)$. In this case, since $\cA_i$ has just $i$ elements, then there exists $\be(l')$ and $\be(m')$ satisfying condition (2), moreover $\be(l')$ cannot be special in bottom line, because $\al(l)$ is special in top line, therefore condition (4) holds. Conversely, suppose that condition (2) holds but (4) is not true. Then arguing as above, we get that condition (3) holds. Now we define the letter $\chi(i+1)$. If condition (3) holds, observing that $m-l\leq d-i$ and applying part (1) of Lemma \ref{lem4s4ss3controloncords4} we get $j$ with $l<j<m$ such that
$$
\la_{\al(j)}>
\frac{\de}{d-i}\cdot
\min\{\la_{\al(l)},\la_{\al(m)}\}.
$$
The letter $\chi(i+1):=\al(j)$ satisfies (\ref{eqconditionBprop1s4ss3controloncords5}). If condition (4) holds, then we apply part (2) of Lemma \ref{lem4s4ss3controloncords4} and the same argument.
\end{dimo}

\subsubsection{Proof of Proposition \ref{prop1s4ss1}}

Consider data $(\pi,\la,\tau)$ with $\area(\pi,\la,\tau)=1$. In order to prove Proposition \ref{prop1s4ss1} it is enough to prove that for any $\al$ in $\cA$ we have
$$
\frac{\la_\al}{\|\la\|}\geq
\frac{1}{d!}
m(\pi,\la,\tau)^{d-1}.
$$
Let $\de$ be the minimum of $\area(v)$ over all cords $v$ with respect to $(\pi,\la,\tau)$. According to Proposition \ref{prop1s4ss2cords} we have $\de\geq m(\pi,\la,\tau)$. Then the required estimate follows from Proposition \ref{prop1s4ss3controloncords5}. The proof of Proposition \ref{prop1s4ss1} is complete.

\subsubsection{Geometric bounds via the control on positive matrices}

Fix data $(\pi,\la,\tau)$ and let $X$ be the underlying translation surface. Let $T$ be the IET determined by the pair $(\pi,\la)$, acting on the interval
$I=(0,\sum_{\chi\in\cA}\la_\chi)$, which is embedded in $X$ along the horizontal direction and with the left endpoint in a conical singularity. The vertical flow of $X$ corresponds to the suspension flow $\phi^t$ over $T$ under the roof function $h=h(\pi,\tau)$ defined in \S \ref{backgroundss1}. For letters $\be$ and $\al$ with $\pi^t(\al)>1$ and $\pi^b(\be)>1$ consider the integer vectors $w^b_{\be,\pi}$, $w^t_{\al,\pi}$ and
$
w_{\be,\al,\pi}=
w^b_{\be,\pi}-w^t_{\al,\pi}
$,
defined in \S \ref{s3ss1sss1}. Recall that the singularities of $T$ and $T^{-1}$ are given by
$
u^t_\al=
\langle w^t_{\al,\pi},\la\rangle
$
and
$
u^b_\be=
\langle w^b_{\be,\pi},\la\rangle
$
respectively, and thus the values
$
|\langle\la,w_{\be,\al,\pi}\rangle|
$
represent the distance between singularities of $T$ and of $T^{-1}$. When
$
v:=\langle\zeta,w_{\be,\al,\pi}\rangle
$
is a period of $X$, then
$
|\langle\la,w_{\be,\al,\pi}\rangle|=
|\Re(v)|
$.
Denote
$
\|\la\|_\infty:=\max_\chi|\la_\chi|
$.

\begin{lem}\label{lem5s4ss4controlonpositivematrices1}
Let $\ga$ and $\ga'$ be two positive paths that can be concatenated and such that $\ga\ast\ga'$ ends in $\pi$, then consider $\la$ in
$\De_{\ga\ast\ga'}$. For any $\be$ and $\al$ as above we have
$$
|\langle\la,w_{\be,\al,\pi}\rangle|
\geq
\frac
{\|\la\|_\infty}
{\|B_\ga\|_\infty\cdot\|B_{\ga'}\|_\infty}.
$$
\end{lem}

\begin{dimo}
Let $(\pi',\la')$ be the image of $(\pi,\la)$ of the branch of the Rauzy map $Q_\ga$ determined by $\ga$. Since $\ga$ is positive, then Corollary \ref{cor1s4ss2} implies
$$
|\langle\la,w_{\be,\al,\pi}\rangle|
\geq
\min_{\chi\in\cA}\la'_\chi.
$$
On the other hand, $\la'$ belongs to $\De_{\ga'}$, that is the image under $^tB_{\ga'}$ of the cone $\RR^\cA_+$, thus we have
$$
\|\la'\|_\infty\leq
\|B_{\ga'}\|_\infty\cdot
\min_{\chi\in\cA}\la'_\chi.
$$
Finally, since $\la=^tB_\ga\la'$ then considering the norm of the linear operator $^tB_\ga$ we get
$$
\|\la\|_\infty
\leq
\cdot\|B_{\ga}\|_\infty\cdot \|\la'\|_\infty.
$$
The Lemma follows summing up the three results.
\end{dimo}

Let $\be$ and $\al$ be letters such that $\pi^t(\al)>1$ and $\pi^b(\be)>1$. Recall that
$
\langle w^t_{\al,\pi},\tau\rangle
$
equals the length of the vertical segment connecting $u^t_\al$ to the conical singularity where its orbit ends. Similarly,
$
\langle w^b_{\be,\pi},\tau\rangle
$
equals the length of the vertical segment connecting (in negative time) $u^b_\be$ to the conical singularity where its orbit starts. When
$
v:=\langle\zeta,w_{\be,\al,\pi}\rangle
$
is a period of $X$, then
$
|\langle\tau,w_{\be,\al,\pi}\rangle|=
|\Im(v)|
$.

\begin{lem}\label{lem6s4ss4controlonpositivematrices2}
Let $\ga$ and $\ga'$ be two positive paths that can be concatenated and such that $\ga\ast\ga'$ starts at $\pi$, then consider $\tau$ in
$\Th_{\ga\ast\ga'}$. For any $\be$ and $\al$ as above we have
$$
|\langle
\tau,w_{\be,\al,\pi}
\rangle|
\geq
2\cdot
\frac
{\|h\|_\infty}
{\|B_\ga\|_\infty\cdot\|B_{\ga'}\|_\infty}.
$$
\end{lem}

\begin{dimo}
Let $(\pi',\la',\tau')$ be the data determined by $\ga'$ via the relation (\ref{eq1s4ss2data}). Let $T'$ be the IET corresponding to the pair $(\pi',\la')$ and $I'$ be the interval where $T'$ acts. The interval $I'$ contains $I$ as subinterval (with the same left endpoint) and is horizontally embedded in $X$, so that $T'$ is the first return map to $I'$ of the vertical flow $\phi$ on $X$. Let $h'=h'(\pi',\tau')$ be the vector whose entries are the (piecewise constant) values of the height function for $T'$. Lemma \ref{lem2s4ss2} applies to $T$, since $\ga'$ is positive, thus $I$ does
not contain any singularity of $T'$ or of $(T')^{-1}$. Since $u^t_\al$ is a singularity of $T$, then its positive vertical orbit comes back to $I'$
before reaching a conical singularity, therefore we have
$
\langle w^t_{\al,\pi},\tau\rangle
\geq
\min_{\chi\in\cA}h'_\chi
$.
The symmetric argument in
negative time applies to $u^b_\be$, which is a singularity for $T^{-1}$. Combining the two estimates, we get
$$
|\langle
\tau,w_{\be,\al,\pi}
\rangle|
\geq
2\cdot\min_{\chi\in\cA}h'_\chi.
$$
The two height functions $h$ and $h'$ are related via the co-cycle by the relation $h=B_{\ga'}h'$, therefore we have
$$
\|h\|_\infty\leq
\|B_{\ga'}\|_\infty\cdot\|h'\|_\infty.
$$
Moreover the same relation holds for the matrix $B_\ga$, that is we have $h'=B_\ga h''$ for some vector $h''$ in $\RR^\cA_+$, hence positivity of $\ga$ implies that for any $\chi\in\cA$ we have
$$
\|h'\|_\infty\leq
\|B_\ga\|_\infty\cdot h'_\chi
$$
The Lemma follow summing up the three results above.
\end{dimo}

Lemma \ref{lem1backgroundss1} and Lemma
\ref{lem6s4ss4controlonpositivematrices2} immediately imply the following Corollary.

\begin{cor}\label{cor1s4ss4controlonpositivematrices3}
Let $\ga$ and $\ga'$ be two positive paths that can be concatenated and such that $\ga\ast\ga'$ start at $\pi$, then consider $\tau$ in
$\Th_{\ga\ast\ga'}$. For any $\be$ and $\al$ as above we have$$
|\langle\tau,w_{\be,\al,\pi}\rangle|
\geq
2\cdot\frac
{\|\tau\|_\infty}
{\|B_\ga\|_\infty\cdot\|B_{\ga'}\|_\infty}.
$$
\end{cor}

\subsubsection{Proof of Proposition \ref{prop2s4ss1}}

Let $\ga_1,\ga_2,\ga_3,\ga_4$ be positive paths that can be
concatenated so that $\pi$ is the ending point of $\ga_1\ast\ga_2$ and the starting point of $\ga_3\ast\ga_4$. Consider length-suspension data $(\la,\tau)$ for $\pi$ with
$
\area(\pi,\la,\tau)=\langle\la,h\rangle=1
$
and such that $\tau$ belongs to the sub-cone $\Th_{\ga_1\ast\ga_2}$ of $\Th_\pi$ and $\la$ belongs to $\De_{\ga_3\ast\ga_4}$. Fix any pair of letters $\be$ and $\al$ with
$\pi^t(\al)>1$ and $\pi^b(\be)>1$. Lemmas \ref{lem5s4ss4controlonpositivematrices1} and
\ref{lem6s4ss4controlonpositivematrices2} imply
$$
|\langle\la,w_{\be,\al,\pi}\rangle|
\cdot
|\langle\tau,w_{\be,\al,\pi}\rangle|
\geq
\frac
{\|\la\|_\infty}
{\|B_{\ga_3}\|_\infty\cdot\|B_{\ga_4}\|_\infty}
\cdot
\frac
{2\cdot\|h\|_\infty}
{\|B_{\ga_1}\|_\infty\cdot\|B_{\ga_2}\|_\infty}
\geq
$$
$$
\frac{2}{d}\cdot
\frac{1}
{\|B_{\ga_1}\|_\infty\cdot
\|B_{\ga_2}\|_\infty\cdot
\|B_{\ga_3}\|_\infty\cdot
\|B_{\ga_4}\|_\infty},
$$
where the second inequality follows combining the inequality
$
\|\cdot\|_2\leq \sqrt{d}\|\cdot\|_\infty
$
between the $L_\infty$ and $L_2$ norms and the Cauchy-Schwartz
inequality applied to the scalar product
$\langle\la,h\rangle=1$. Since the above inequality holds for any $\be$ and
$\al$, then it holds for $w(\pi,\la,\tau)$ too. Proposition \ref{prop2s4ss1} is proved.


\begin{thebibliography}{0000}

\bibitem[A,G,Y]{agy}
A.Avila, S.Gouezel, J.-C.Yoccoz:
\emph{"Exponential mixing for the Teichm\"uller
flow"}.
Publications math\'ematiques de l'IHES 104 (2006), 143-211.

\bibitem[B]{birkhoff}
G. Birkhoff:
\emph{"Extensions of Jentzsth's Theorem"}.
Trans. Amer. Math. Soc. 85 (1957), 219-227.

\bibitem[Boi,L]{boissylanneau}
C. Boissy, E. Lanneau:
\emph{"Dynamics and geometry of the Rauzy-Veech induction for quadratic differentials."}
Ergodic Theory Dynam. Systems 29 (2009), no. 3, 767-816.

\bibitem[Bos]{boshernitzan}
M. Boshernitzan:
\emph{"A condition for minimal interval exchange maps to be uniquely ergodic".}
Duke Math. J. 52 (1985) 723-752.

\bibitem[C,F]{cusick}
T. W. Cusick, M. E. Flahive:
\emph{"The Markoff and Lagrange Spectra"}.
Mathematicas Surveys and Monographs, 30, 1989.

\bibitem[E,Mi]{eskinmirza}
A. Eskin, M. Mirzakhani:
\emph{"Invariant and stationary measures for the $\slduer$ action on moduli space."}
Preliminary version of the preprint at http://www.math.uchicago.edu/~eskin/.

\bibitem[E,Ra,Mi]{EMR}
A. Eskin, K. Rafi, M. Mirzakhani:
\emph{"Counting closed geodesics in strata."}
ArXiv:1206.5574

\bibitem[F]{ferenczi}
S. Ferenczi:
\emph{"Dynamical generalizations of the Lagrange spectrum"}.
Arxiv:1108.3628.

\bibitem[H]{hall}
M. Hall:
\emph{"On the sum and products of continued fractions"}.
Annals of Math, 48, (1947), 966-993.

\bibitem[Ha1]{haas1}
A. Haas:
\emph{"Diophantine approximation on hyperbolic Riemann surfaces."}
Acta Math. 156 (1986), no. 1-2, 33-82.

\bibitem[Ha2]{haas2}
A. Haas:
\emph{"Diophantine approximation on hyperbolic orbifolds."} Duke Math. J. 56 (1988), no. 3, 531-547.

\bibitem[Ha,S]{haasseries}
A. Haas, C. Series:
\emph{"The Hurwitz constant and Diophantine approximation on Hecke groups"}.
J. London Math. Soc. 2, 34, 1986, 219-334.

\bibitem[Ham]{ursula}
U.~Hamenst{\"a}dt:
\emph{"Dynamics of the Teichmueller flow on compact invariant sets."}
J. Mod. Dynamics 4 (2010), 393-418.

\bibitem[He,P]{hersonskypulin}
S. Hersonsky, F. Paulin:
\emph{"On the almost sure spiraling of geodesics in negatively curved manifolds."}
J. Differential Geom. 85 (2010), no. 2, 271-314.

\bibitem[Hu,L]{hubertlelievre}
P. Hubert, S. Leli{\`e}vre:
\emph{"Prime arithmetic Teichm\"uller disc in $\cH(2)$"}. . Israel J. Math. 151 (2006), 281-321.

\bibitem[Ke]{keane}
M. Keane:
\emph{"Interval exchange transformations"}.
Math. Z., 141, 25-31, (1975).

\bibitem[K,M]{kimmarmi}
D. H. Kin, S. Marmi:
\emph{"Bounded type interval exchange maps."}
In preparation.

\bibitem[Ki]{kin}
Khinchin:
\emph{"Continued fractions"}.
English translation, P.Noordhoff, Groningen, (1963).

\bibitem[Mar]{luca1}
L. Marchese:
\emph{"The Khinchin theorem for interval exchange transformations"}.
J. Mod. Dyn., Volume 5, No.1, 123-183, (2011).

\bibitem[M,M,Y]{mmy}
S. Marmi, P. Moussa, J.-C. Yoccoz:
\emph{"The cohomological equation for Roth-type interval exchange maps."}
J. Amer. Math. Soc. 18 (2005), no. 4, 823-872.

\bibitem[Mau]{maucourant}
F. Maucourant:
\emph{"Sur les spectres de Lagrange et de Markoff des corps imaginaires quadratiques"} (French).
Ergodic Theory Dynam. Systems 23 (2003), no. 1, 193-205.

\bibitem[M]{masur}
H. Masur:
\emph{"Interval exchange transformation and measured foliations"}.
Annals of Mathematics, 115, 169--200, (1982).

\bibitem[Mor]{gugu}
G. Moreira:
\emph{"Introdu\c{c}ao \`a teoria dos n\'umeros"}. (Portuguese) [Introduction to number theory]
Monograf\'ias del Instituto de Matem\'atica y Ciencias Afines, 24. Instituto de Matem\'atica y Ciencias Afines, IMCA, Lima; Pontificia Universidad Cat\'olica del Per\'u, Lima, 2002. 110 pp. ISBN: 9972-899-01-2, 11-01

\bibitem[P,P1]{paulin1}
J. Parkkonen, F. Paulin:
\emph{"On the closedness of approximation spectra."}
J. Th\'eor. Nombres Bordeaux 21 (2009), no. 3, 701-710.

\bibitem[P,P2]{paulin2}
J. Parkkonen, F. Paulin:
\emph{"Prescribing the behaviour of geodesics in negative curvature."}
Geom. Topol. 14 (2010), no. 1, 277-392.

\bibitem[P,P3]{paulin3}
J. Parkkonen, F. Paulin:
\emph{"Spiraling spectra of geodesic lines in negatively curved manifolds."}
Math. Z. 268 (2011), no. 1-2, 101-142.

\bibitem[Pat]{patterson}
S. J. Patterson:
\emph{"Diophantine approximation in Fuchsian groups."} Philos. Trans. Roy. Soc. London Ser. A 282 (1976), no. 1309, 527-563,

\bibitem[R]{rauzy}
G. Rauzy:
\emph{"Echanges d'intervalles et transformations induites"}. (French)
Acta Arith.
34, no. 4, 315-328, (1979).

\bibitem[Ro]{ibarra}
S. A. Roma\~{n}a:
\emph{"Ph-D thesis"}.
In preparation.

\bibitem[S,S]{schmidtsheingorn}
T. A. Schmidt, M. Sheingorn:
\emph{"Riemann surfaces have Hall rays at each cusp"}. Illinois J. Math. 41 (1997), no. 3, 378-397.

\bibitem[Se0]{series0}
C. Series:
\emph{"The geometry of Markoff numbers"}.
Math. Intelligencer 7 (1985) 20–29.

\bibitem[Se1]{series1}
C. Series:
\emph{"The modular surface and continued fractions"}.
J. Lond. Math. Soc. 31 (1985) 69-80.

\bibitem[Se2]{series2}
C. Series:
\emph{"The Markoff spectrum in the Hecke group G5"},
J. London Math. Soc. 57, 1988, 151-181.

\bibitem[Sm,We]{smillieweiss}
J. Smillie, B. Weiss:
\emph{"Characterizations of lattice surfaces".}
Invent. Math. 180 (2010), no. 3, 535--557.

\bibitem[Ve1]{veech}
W. Veech:
\emph{"Gauss measures for transformations on the space of interval exchange maps"}.
Annals of Mathematics, 115, 201-242, (1982).


W. Veech
\bibitem[Ve2]{veechboshernitzancriterion}
W. Veech:
\emph{"Boshernitzan's criterion for unique ergodicity of an interval exchange transformation"}.
Ergodic Theory Dynam. Systems 7 (1987), 1, 27-48.


\bibitem[Ve3]{veechveechcurves}
W. Veech:
\emph{"Teichm{\"u}ller curves in moduli space, Eisenstein series and an application to triangular billiards"}.
Invent. Math.  97 (1989), 3, 553-583.

\bibitem[Vi]{viana}
M. Viana:
\emph{"Lecure notes on interval exchange transformations"}.

\bibitem[Vo]{vorobets}
Y. Vorobets:
\emph{"Planar structures and billiards in rational polygons: the Veech alternative". (Russian)}.
Uspekhi Mat. Nauk 51 (1996), no. 5(311), 3-42; translation in Russian Math. Surveys 51 (1996), no. 5, 779-817.

\bibitem[Vu1]{Vulakh:Triangle} 
L. Vulakh:
\emph{"The Markov Spectra for Triangle Groups"}.
J. Number Theory, Vol. 67 No. 1 (1997) 11-28.

\bibitem[Vu2]{Vulakh:Fuchsian}
L. Vulakh:
\emph{"The Markov spectra for Fuchsian groups."}  
Trans. Amer. Math. Soc. 352 (2000), 4067-4094.

\bibitem[Vu3]{vulakh3}
L. Vulakh:
\emph{"Diophantine Approximation in $R^n$."}
Trans. Amer. Math. Soc., Volume 347, No.  2 (1995) , 573--585.

\bibitem[Vu4]{Vulakh:Farey}  
L. Vulakh:
\emph{"Farey polytopes and continued
fractions associated with discrete hyperbolic groups."} Trans. Amer.
Math. Soc. 351 (1999), 2295-2323.

\bibitem[Vu5]{vulakh5}
L. Vulakh:
\emph{"Diophantine approximation on Bianchi groups"},
J. Number Theory, Vol. 54, No. 1 (1995) 73-80.

\bibitem[W]{Wright}
A. Wright:
\emph{"The field of definition of affine invariant submanifolds of the moduli space of abelian differentials"}. ArXiv:1210.4806
	
\bibitem[Y]{yoccoz}
J.-C. Yoccoz:
\emph{"Echanges d'intervalles"}.
Cours Coll\`ege de France, Janvier-Mars 2005.

\bibitem[Z1]{zorichuno}
A. Zorich:
\emph{"Finite Gauss measure on the space of interval exchange transformations. Lyapunov exponents"}.
Annales de l'Institut Fourier, 46:2, 325-370 (1996).

\bibitem[Z2]{zorichdue}
A. Zorich:
\emph{"Flat surfaces"}.
Frontiers in Number Theory, Physics and Geometry, Vol 1;
P.Cartier; B.Julia; P.Moussa; Vanhove (Editors), Springer-Verlag, 403-437,
(2006).

\end{thebibliography}
\end{document}